\documentclass[a4paper,11pt]{article}

\usepackage{authblk}

\usepackage{amssymb}
\let\Emptyset\emptyset

\usepackage{type1cm}        
\usepackage{makeidx}         
\usepackage{graphicx}        
\usepackage{multicol}        

\usepackage{newtxtext}       %
\usepackage{newtxmath}       

\makeindex             

\usepackage{mdframed}

\usepackage{tikz}
\usetikzlibrary{shapes,arrows}
\usetikzlibrary{patterns}

\usepackage{pgfplots}
\pgfplotsset{compat=1.13}
\usepgfplotslibrary{fillbetween}

\definecolor{PrimalOrange}{RGB}{255,153,85}
\definecolor{DualBlue}{RGB}{0,255,255}
\definecolor{SlabGreen}{RGB}{170,212,0}
\definecolor{HSUred}{RGB}{197,0,66}
\definecolor{mountainmeadow}{rgb}{0.19, 0.73, 0.56}
\definecolor{navyblue}{rgb}{0.0, 0.0, 0.5}

\definecolor{Goldenrod}{RGB}{255,185,15}
\definecolor{Eggshell}{RGB}{255,236,139}
\definecolor{LightOlive}{RGB}{204,255,153}

\usepackage{hyperref}
\usepackage{xcolor}

\usepackage{subfig}
\usepackage{graphicx}

\renewcommand{\emptyset}{\Emptyset}

\newcommand{\doi}[1]{\href{https://doi.org/#1}{doi:#1}}

\newcommand{\concentration}{{u}}
\newcommand{\convection}{{\boldsymbol v}}
\newcommand{\pressure}{{p}}
\newcommand{\dualz}{{z}}

\newcommand{\viscosity}{{\nu}}
\newcommand{\transportforce}{{g}}
\newcommand{\stokesforce}{{\boldsymbol f}}

\newtheorem{defi}{Definition}[section]
\newtheorem{theorem}[defi]{Theorem}

\newtheorem{remark}[defi]{Remark}

\newenvironment{mproof}{\paragraph{Proof.}}{\hfill$\blacksquare$}

\numberwithin{equation}{section}
\numberwithin{table}{section}
\numberwithin{figure}{section}

\usepackage{geometry}
\begin{document}

\title{\Large \textbf{On the implementation of an adaptive multirate framework for
coupled transport and flow}}
\author[M.\ P.\ Bruchh\"auser, U.\ K\"ocher, M.\ Bause]
{\large \textbf{Marius Paul Bruchh\"auser}\thanks{bruchhaeuser@hsu-hamburg.de ($^\ast$corresponding author)} 
$\,\boldsymbol{\cdot}$
\textbf{Uwe K\"ocher}\thanks{koecher@hsu-hamburg.de}
$\,\boldsymbol{\cdot}$
\textbf{Markus Bause}\thanks{bause@hsu-hamburg.de}\\
{\small Helmut Schmidt University, University of the Federal Armed Forces Hamburg,
Faculty of Mechanical Engineering, Chair of Numerical Mathematics\\ 
Holstenhofweg 85, 22043 Hamburg, Germany}
}
\date{}
\maketitle

\begin{abstract}
\noindent
In this work, a multirate in time approach resolving the different 
time scales of a convection-dominated transport and coupled fluid flow is 
developed and studied in view of goal-oriented error control by means of 
the Dual Weighted Residual (DWR) method. 
Key ingredients are an arbitrary degree discontinuous Galerkin time discretization
of the underlying subproblems, an a posteriori error representation for the 
transport problem coupled with flow and its implementation using space-time
tensor-product spaces.
The error representation allows the separation of the temporal and spatial 
discretization error which serve as local error indicators for adaptive mesh 
refinement.
The performance of the approach and its software implementation are studied
by numerical convergence examples as well as an example of physical interest for
convection-dominated transport.
\end{abstract}


\section{Introduction}
\label{sec:1:introduction}
In recent years, mathematical models of multi-physics coupling subproblems with 
characteristic time scales that differ by orders of magnitude have attracted
researchers' interest; cf., e.g., \cite{Jammoul2021,Ge2018,Gupta2016,Almani2016}. 
Their efficient numerical simulation with regard to the temporal discretization 
does not become feasible without using techniques adapted to these 
characteristic scales that resolve the solution components on their respective 
time length by an adaptation of the time steps sizes.
Such methods are referred to as multirate in time (for short, multirate) schemes. 
Firstly, they were introduced for the numerical approximation of systems of 
ordinary differential equations in \cite{Gear1984,Guenther1993}.
For a short review of multirate methods including a list of references we refer
to \cite{Gupta2016,Gander2013}.

In this work we focus on the multirate implementation of a fully space-time 
adaptive convection-dominated transport problem coupled with a time-dependent 
Stokes flow problem.
The implementation is based on our open-source code given by \cite{Koecher2019} 
for the \texttt{deal.II} finite element analysis library; 
cf.~\cite{dealiiReference93}.
With regard to our coupled model problem, we assume a highly time-dynamic 
process modeled by the transport equation such that the underlying temporal mesh
is discretized using smaller time step sizes compared to a slowly moving process
modeled by the viscous flow problem.
Our motivation comes through the definition of so-called characteristic times 
for the two subproblems that serve as quantities to measure the underlying 
dynamic in time and have their origin in the field of natural sciences and 
engineering sciences, cf., e.g., \cite{Gujer2008,Morgenroth2015}.
For the sake of physical realism, the transport problem is supposed to be
convection-dominated by assuming high P\'{e}clet numbers that are characterized 
by small diffusion relative to the convection, cf.~\cite{Burman2014,John2018}. 
The solution of these transport problems are typically characterized by the
occurrence of sharp moving fronts and layers. The key challenge for the
numerical approximation exists in the accurate and efficient solution while
avoiding non-physical oscillations or smearing effects.
The application of stabilization techniques is a typical approach to overcome
non-physical effects.
As shown in a comparative study for time-dependent convection-diffusion-reaction
equations in \cite{John2009}, stabilization techniques on globally refined meshes
fail to avoid these oscillations even after tuning stabilization parameters.
For a general review of stabilization techniques we refer to
\cite{Roos2008,John2018}.

For the efficient numerical simulation of multi-physics problems handling the 
challenges described above, it is indisputable that adaptive mesh refinement 
strategies in space and time are necessary. One possible technique for those 
adaptive strategies is goal-oriented a posteriori error control based on the 
Dual Weighted Residual method \cite{Becker2001,Bangerth2003}. 
For a general review of a posteriori error estimation we refer to
\cite{Ainsworth2000,Verfuerth1996}.

In this work we follow our approaches and implementations from
\cite{Bause2021,Koecher2019}.
An extension to our preceding work is that the flow problem now depends on time
and needs to be solved on a different time scale than the transport problem.
Precisely, this work is characterised by the following features.
\begin{itemize}
\item Development of a multirate concept with independent time scales for the
transport and flow problem, respectively.

\item Implementation of tensor-product space-time slabs for an arbitrary order
discontinuous Galerkin (dG) time discretization.

\item Implementation of coupling the Stokes flow velocity to the transport
problem using interpolation techniques between different finite element spaces
and meshes.

\end{itemize}
This work is organized as follows. In Sec. \ref{sec:2} we introduce the model
problem, the multirate decoupling of the transport and flow problems and their
space-time discretizations. In Sec. \ref{sec:3:dwr} we derive an a posteriori 
error representation for the transport problem. In Sec. \ref{sec:4:implementation} 
we explain the implementation of the space-time tensor-product spaces. The 
underlying algorithm and some related aspects are presented in 
Sec. \ref{sec:5:algorithm}. Numerical examples are given in
Sec. \ref{sec:6:examples} and in Sec. \ref{sec:7:conclusion}
we summarize with conclusions and give some outlook for future work.

\section{Model Problem, Multirate and Space-Time Discretization}
\label{sec:2}

In Subsection \ref{sec:2:1:modelproblem} we introduce the model problem 
of a convection-diffusion-reaction transport coupled with a time-dependent 
Stokes flow. For instance, such system is used to model species or heat 
transport in a creeping viscous fluid. 
Beyond that, such multi-physics systems of coupled flow and transport serve as 
prototype models for applications in several branches of natural and engineering 
scienes, for instance, contaminant transport and degradation in the subsurface, 
reservoir simulation, fluid-structure interaction, and thermal and mass 
transport in deformable porous media or  thermal expansion in solid mechanics; 
cf., e.g.,  \cite{Larson2007,Allaire1989,Wick2016,Odsaeter2019,Bengzon2010}.

In Section \ref{sec:2:2:multirate}, we explain our multirate in time approach for
the two subproblems, before we present the details of the space-time 
discretizations in Sections \ref{sec:2:3}-\ref{sec:2:5}.

\subsection{Model Problem}
\label{sec:2:1:modelproblem}

The time dependent convection-diffusion-reaction transport problem in 
dimensionless form is given by
\begin{equation}
\label{eq:1:transport_problem}
\begin{array}{rcl @{\,\,}l @{\,\,}l @{\,}l}
\partial_{t} \concentration
- \nabla \cdot (\varepsilon \nabla \concentration)
+ \convection \cdot \nabla \concentration
+ \alpha \concentration &=& \transportforce
& \text{in} & Q & = \Omega \times I\,,\\[.5ex]
\concentration &=& \concentration_D
& \text{on} & \Sigma_D & = \Gamma_D \times I\,,\\[.5ex]
\varepsilon \nabla \concentration \cdot \boldsymbol n &=& \concentration_N
& \text{on} & \Sigma_N & = \Gamma_N \times I\,,\\[.5ex]
\concentration &=&  \concentration_0
& \text{on} & \Sigma_0 & = \Omega \times \{ 0 \}\,,
\end{array}
\end{equation}
for a boundary partition $\partial\Omega = \Gamma_D \cup \Gamma_N$,
$\Gamma_D \neq \emptyset$ with outer unit normal vector $\boldsymbol n$.
The characteristic time $t_{\textnormal{transport}}$ of this transport equation
\eqref{eq:1:transport_problem} can be comprehended as a dimensionless time 
variable depending on the diffusive, convective as well as reactive part and is 
here defined by
\begin{equation}
\label{eq:2:characteristic-time-transport}
t_{\textnormal{transport}} := \min \Bigg\{\frac{L^2}{\varepsilon}\,;\, \frac{L}{V}\,;\,
\frac{1}{\alpha}\Bigg\}\,,
\end{equation}
where $0 < \varepsilon \ll 1$ is the diffusion coefficient, $\alpha > 0$ is the 
reaction coefficient, $L$ denotes the characteristic length of the domain 
$\Omega$, for instance, its diameter, and $V$ denotes a characteristic velocity of the flow field $\convection $, 
for instance, the mean inflow velocity given by 
$\frac{1}{T\, |\Gamma_{\text{inflow}}|} \int_I \int_{\Gamma_{\textnormal{inflow}}}
\convection_D\cdot (- \boldsymbol n) \;\mathrm{d} o\; \mathrm{d} t$;
cf.\cite{Gujer2008,Morgenroth2015} for more details.

The convection field $\convection$ in the transport problem \eqref{eq:1:transport_problem}
is determined by the dimensionless Stokes flow system
\begin{equation}
\label{eq:3:stokes_problem}
\begin{array}{rcl @{\,\,}l @{\,\,}l @{\,}l}
\partial_t \convection
- \nabla \cdot (2\viscosity\, \boldsymbol\epsilon(\convection)
+ \pressure \boldsymbol I) &=& \stokesforce
& \text{in} & Q & = \Omega \times I\,,\\[.5ex]
\nabla \cdot \convection &=& 0
& \text{in} & Q & = \Omega \times I\,,\\[.5ex]
\convection &=& \convection_D
& \text{on} & \Sigma_{\textnormal{inflow}}
& = \Gamma_{\textnormal{inflow}} \times I\,,\\[.5ex]
\convection &=& \boldsymbol 0
& \text{on} & \Sigma_{\textnormal{wall}}
& = \Gamma_{\textnormal{wall}} \times I\,,\\[.5ex]
%
(2\viscosity\, \boldsymbol\epsilon(\convection) + \pressure \boldsymbol I)
\boldsymbol n &=& \boldsymbol 0 & \text{on} & \Sigma_{\textnormal{outflow}}
& = \Gamma_{\textnormal{outflow}} \times I\,,\\[.5ex]
\convection &=& \convection_0
& \text{on} & \Sigma_0 & = \Omega \times \{ 0 \}\,,\\[.5ex]
\end{array}
\end{equation}
for a boundary partition
$\partial\Omega = \Gamma_{\textnormal{inflow}} \cup \Gamma_{\textnormal{wall}}
\cup \Gamma_{\textnormal{outflow}}$ which is (in general) independent from the
boundary partition of the transport problem. The appropriate choice for the
boundary partition and setting of the inflow profiles is standard and can be 
found in the literature \cite{John2016}.
The characteristic time $t_{\textnormal{flow}}$ of the 
Stokes flow equation is then defined by
\begin{equation}
\label{eq:4:characteristic-time-flow}
t_{\textnormal{flow}} := \frac{L}{V}\,,
\end{equation}
with $L$ and $V$ being chosen as in \eqref{eq:2:characteristic-time-transport}. With regard to the characteristic times of the two subproblems, we assume that $t_{\textnormal{transport}} \ll t_{\textnormal{flow}}$ 
such that we are using a finer temporal mesh to resolve the dynamics of a faster process 
given by the transport equation compared to the slower process of the viscous, 
creeping flow. This multirate in time approach is described in detail in the following section. 

In \eqref{eq:1:transport_problem}, \eqref{eq:3:stokes_problem}, we denote by
$\Omega \subset \mathbb{R}^{d}$, with $d=2,3$, a polygonal or polyhedral bounded
domain with Lipschitz boundary $\partial\Omega$ and $I=(0,T]$, $0 < T < \infty$,
is a finite time interval.
We assume that 
and
$\viscosity > 0$ is a viscosity coefficient.
Well-posedness of \eqref{eq:1:transport_problem}, \eqref{eq:3:stokes_problem}
and the existence of a sufficiently regular solution, such that all of the
arguments and terms used below are well-defined, are tacitly assumed without
mentioning explicitly all technical assumptions about the data and coefficients,
cf.~\cite{Roos2008} and \cite{Ern2021}.

\subsection{Multirate}
\label{sec:2:2:multirate}

\begin{figure}
\centering

\begin{tikzpicture}
\tikzstyle{ns1}=[line width=1.]

\node at (10.3,.4) {\small Transport};
\node at (10.2,-.2) {\small $t$};
\draw[->,ns1] (0,0) -- (10,0);

\draw[->,ns1] (0,-.5) -- (10,-.5);
\node at (10.5,-.9) {\small Stokes flow};

\node at (5.0,-.27) {\small $\ell=1$: Initialization loop};

\node at (0.2, 0.4) {\small $0$};
\draw[ns1] (0.2, -0.1) -- (0.2, 0.1);

\node at (1.2, 0.4) {\small $t_1^\textnormal{T}$};
\draw[ns1] (1.2, -0.1) -- (1.2, 0.1);

\node at (2.2, 0.4) {\small $t_2^\textnormal{T}$};
\draw[ns1] (2.2, -0.1) -- (2.2, 0.1);

\node at (3.2, 0.4) {\small $t_3^\textnormal{T}$};
\draw[ns1] (3.2, -0.1) -- (3.2, 0.1);

\node at (4.2, 0.4) {\small $t_4^\textnormal{T}$};
\draw[ns1] (4.2, -0.1) -- (4.2, 0.1);

\node at (5.2, 0.4) {\small $t_5^\textnormal{T}$};
\draw[ns1] (5.2, -0.1) -- (5.2, 0.1);

\node at (6.2, 0.4) {\small $t_6^\textnormal{T}$};
\draw[ns1] (6.2, -0.1) -- (6.2, 0.1);

\node at (7.2, 0.4) {\small $t_7^\textnormal{T}$};
\draw[ns1] (7.2, -0.1) -- (7.2, 0.1);

\node at (8.2, 0.4) {\small $t_8^\textnormal{T}$};
\draw[ns1] (8.2, -0.1) -- (8.2, 0.1);

\node at (9.2, 0.4) {\small $T$};
\draw[ns1] (9.2, -0.1) -- (9.2, 0.1);

\draw[ns1] (0.2, -0.6) -- (0.2, -0.4);
\node at (0.2, -0.9) {\small $0$};

\draw[ns1] (2.2, -0.6) -- (2.2, -0.4);
\node at (2.2, -0.9) {\small $t_1^\textnormal{F}$};

\draw[ns1] (5.2, -0.6) -- (5.2, -0.4);
\node at (5.2, -0.9) {\small $t_2^\textnormal{F}$};

\draw[ns1] (6.2, -0.6) -- (6.2, -0.4);
\node at (6.2, -0.9) {\small $t_3^\textnormal{F}$};

\draw[ns1] (9.2, -0.6) -- (9.2, -0.4);
\node at (9.2, -0.9) {\small $T$};

\end{tikzpicture}

\vskip3ex
\begin{tikzpicture}
\tikzstyle{ns1}=[line width=1.]

\node at (10.3,.4) {\small Transport};
\node at (10.2,-.2) {\small $t$};
\draw[->,ns1] (0,0) -- (10,0);

\draw[->,ns1] (0,-.5) -- (10,-.5);
\node at (10.5,-.9) {\small Stokes flow};

\node at (5.0,-.27) {\small $\ell=2$: First adaptively refined loop};

\node at (0.2, 0.4) {\small $0$};
\draw[ns1] (0.2, -0.1) -- (0.2, 0.1);

\node at (1.2, 0.4) {\small $t_1^\textnormal{T}$};
\draw[ns1] (1.2, -0.1) -- (1.2, 0.1);

\node at (1.7, 0.4) {\small $t_2^\textnormal{T}$};
\draw[ns1] (1.7, -0.1) -- (1.7, 0.1);

\node at (2.2, 0.4) {\small $t_3^\textnormal{T}$};
\draw[ns1] (2.2, -0.1) -- (2.2, 0.1);

\node at (3.2, 0.4) {\small $t_4^\textnormal{T}$};
\draw[ns1] (3.2, -0.1) -- (3.2, 0.1);

\node at (4.2, 0.4) {\small $t_5^\textnormal{T}$};
\draw[ns1] (4.2, -0.1) -- (4.2, 0.1);

\node at (5.2, 0.4) {\small $t_6^\textnormal{T}$};
\draw[ns1] (5.2, -0.1) -- (5.2, 0.1);

\node at (5.7, 0.4) {\small $t_7^\textnormal{T}$};
\draw[ns1] (5.7, -0.1) -- (5.7, 0.1);

\node at (6.2, 0.4) {\small $t_8^\textnormal{T}$};
\draw[ns1] (6.2, -0.1) -- (6.2, 0.1);

\node at (6.7, 0.4) {\small $t_9^\textnormal{T}$};
\draw[ns1] (6.7, -0.1) -- (6.7, 0.1);

\node at (7.2, 0.4) {\small $t_{10}^\textnormal{T}$};
\draw[ns1] (7.2, -0.1) -- (7.2, 0.1);

\node at (8.2, 0.4) {\small $t_{11}^\textnormal{T}$};
\draw[ns1] (8.2, -0.1) -- (8.2, 0.1);

\node at (9.2, 0.4) {\small $T$};
\draw[ns1] (9.2, -0.1) -- (9.2, 0.1);

\draw[ns1] (0.2, -0.6) -- (0.2, -0.4);
\node at (0.2, -0.9) {\small $0$};

\draw[ns1] (2.2, -0.6) -- (2.2, -0.4);
\node at (2.2, -0.9) {\small $t_1^\textnormal{F}$};

\draw[ns1] (5.2, -0.6) -- (5.2, -0.4);
\node at (5.2, -0.9) {\small $t_2^\textnormal{F}$};

\draw[ns1] (6.2, -0.6) -- (6.2, -0.4);
\node at (6.2, -0.9) {\small $t_3^\textnormal{F}$};

\draw[ns1] (9.2, -0.6) -- (9.2, -0.4);
\node at (9.2, -0.9) {\small $T$};

\end{tikzpicture}

\caption{Illustration of exemplary temporal meshes
for the initial loop and the first adaptively refined loop.
The temporal mesh of the transport solver is adaptively refined and
the mesh of the Stokes flow solver is fixed
as explained in Sec. \ref{sec:2:2:multirate}.}

\label{fig:1:multirate_time_scales}
\end{figure}
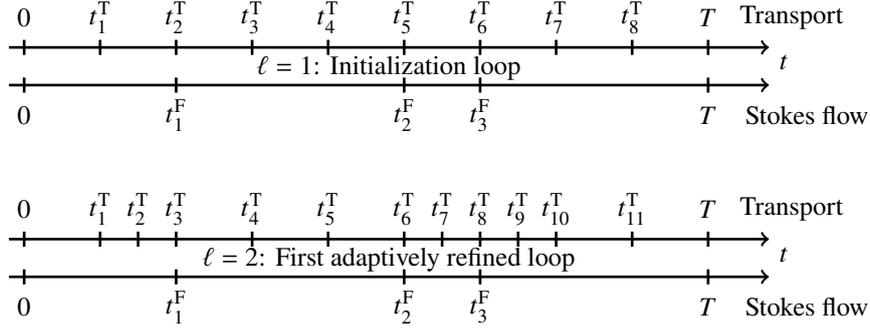

For the efficient approximation we use a multirate in time approach to mimic
the behaviour of a slowly moving fluid, that is approximated by a time-dependent
Stokes solver, and a faster convection-diffusion-reaction process. Precisely, the
problems given by \eqref{eq:1:transport_problem} and \eqref{eq:3:stokes_problem}
are considered on different time scales modeling the underlying physical processes.
We initialize the temporal mesh independently for the Stokes flow and the
transport problem with the following properties
\begin{itemize}
\item the Stokes flow temporal mesh is coarser or equal to that of the transport problem,
\item the endpoints in the temporal mesh of the Stokes solver must match with endpoints
in the temporal mesh of the transport problem.
\end{itemize}
We allow for adaptive time refinements of the temporal mesh of the transport problem
and for global temporal mesh refinements of the Stokes solver over the 
adaptation loops due to the lack of an error estimator for the flow problem.
An exemplary initialization and one manufactured refined temporal mesh are
illustrated in Fig. \ref{fig:1:multirate_time_scales}.

%

For the multirate decoupling of the transport problem,
we let $0 =: t_0^\textnormal{T} < t_1^\textnormal{T} < \dots < t_{N^\ell}^\textnormal{T} := T$ 
a set of time points for the partition of the
closure of the time domain $\bar{I}=[0,T]$ into left-open subintervals
$I_n:=(t_{n-1}^\textnormal{T},t_n^\textnormal{T}]$,
$n=1,\dots,N^\ell$. The number $N^\ell$ depends on the adaptivity loop $\ell$.
For the flow problem, we let $0 =: t_0^\textnormal{F} < t_1^\textnormal{F} < \dots 
< t_{N^{\textnormal{F},\ell}}^{\textnormal{F}} := T$ a set of time points for 
the partition of the closure of the time domain $\bar{I}=[0,T]$ into left-open 
subintervals $I_n^\textnormal{F}:=(t_{n-1}^\textnormal{F},t_n^\textnormal{F}]$,
$n=1,\dots,N^{\textnormal{F},\ell}$.
We approximate the solution $\{\convection, \pressure \}$ of the Stokes flow
problem on each $I_n^\textnormal{F}$
by means of a globally piecewise constant discontinuous Galerkin (dG($0$)) time
approximation.
For simplicity of the implementation, we ensure that each element of the set
$\{ t_0^\textnormal{F}, t_1^\textnormal{F}, \dots, t_{N^{\textnormal{F},\ell}} \}$
corresponds to an element of the set
$\{ t_0^\textnormal{T}, t_1^\textnormal{T}, \dots, t_{N^\ell} \}$.

Additionally, we approximate the solution of the transport problem with
an arbitrary degree $r \ge 0$ in time while we restrict the implementation here
for the time-dependent Stokes flow to a piecewise constant in time approximation.
This gives us for $r>0$ an additional level of the multirate in time character
between the two problems.

\subsection{Weak Formulation}
\label{sec:2:3}

In this section, we present the weak formulation of the transport and Stokes 
flow problem given by Eq.~\eqref{eq:1:transport_problem} and 
Eq.~\eqref{eq:3:stokes_problem}, respectively, to prepare the discretizations
in space and time following below.
Let
$X:=\{u \in L^2(0, T; H^1_0(\Omega)) \mid
\partial_t u \in L^2(0, T; H^{-1}(\Omega))\}$
and $Y_1 := \{ \boldsymbol v \in L^2(0, T; H^1_0(\Omega)^d) \mid
\partial_t \boldsymbol v \in L^2(0, T; H^{-1}(\Omega)^d) \}$.
Then, the weak formulation of \eqref{eq:1:transport_problem}
reads as follows:

\textit{For a given $\convection \in Y_1$ of \eqref{eq:6:weak_stokes}, 
find $\concentration \in X$ such that}
\begin{equation}
\label{eq:5:weak_transport_problem}
A(\concentration;\convection)(\varphi) = G(\varphi) \quad \forall \varphi \in X\,,
\end{equation}
\textit{where the bilinear form $A: \{X;Y_1\} \times X \rightarrow \mathbb{R}$
and the linear form
$G: L^2(0, T;$ $H^{-1}(\Omega))
\rightarrow \mathbb{R}$
are defined by}
\begin{displaymath}
\begin{array}{r@{\,}c@{\,}l}
A(\concentration;\convection)(\varphi) & := & \displaystyle
\int_{I} \big\{
(
\partial_{t} \concentration, \varphi)
+ a(\concentration;\convection)(\varphi) \big\}\; \mathrm{d} t
+ (\concentration(0), \varphi(0)\,,\\[3ex]
G(\varphi) & := & \displaystyle
\int_{I}
  (\transportforce, \varphi)\; \mathrm{d} t
+ (\concentration_{0},\varphi(0) ) \,,
\end{array}
\end{displaymath}
\textit{with the bilinear form}
\begin{displaymath}
a(\concentration;\convection)(\varphi) :=
  (\varepsilon \nabla \concentration, \nabla \varphi)
+ (\convection \cdot \nabla \concentration, \varphi)
+ (\alpha \concentration, \varphi)\,.
\end{displaymath}
%
%
Here, $(\cdot, \cdot)$ denotes the inner product of $L^2(\Omega)$ or duality
pairing of $H^{-1}(\Omega)$ with $H^1_0(\Omega)$, respectively. By $\|\cdot \|$
we denote the associated $L^2$-norm.

For the weak formulation of \eqref{eq:3:stokes_problem}
we additionally define
$Y_2 := \{p \in L^2(0,T;L_0^2(\Omega))\}$, with $L_0^2(\Omega):=\{p \in L^2(\Omega)
\mid \int_\Omega p\; \mathrm{d} \boldsymbol{x} = 0\}$. Then we get:

\textit{For $\stokesforce \in L^2(I;H^{-1}(\Omega)^d)$ and 
$\convection_0 \in L^2(\Omega)^d$, find $\{ \convection, \pressure\} \in Y_1 \times Y_2$, 
such that}
\begin{equation}
\label{eq:6:weak_stokes}
B(\convection,\pressure)(\boldsymbol{\psi},\chi) = F(\boldsymbol{\psi}) \quad
\forall \{\boldsymbol{\psi}, \chi\} \in Y_1 \times Y_2\,,
\end{equation}
\textit{where the bilinear form $B: \{ Y_1 \times Y_2 \} \times 
\{ Y_1 \times Y_2 \} \rightarrow \mathbb{R}$\,
as well as the linear form 
\\
$F: L^2(0,T;H^{-1}(\Omega)^d)~\rightarrow~\mathbb{R}$ are defined by}
\begin{displaymath}
\begin{array}{rcl}
B(\convection,\pressure)(\boldsymbol{\psi},\chi) &:=&
\displaystyle
\int_{I} \big\{
(\partial_{t} \convection,\boldsymbol{\psi})
+ (2 \viscosity \boldsymbol\epsilon(\convection), \boldsymbol\epsilon(\boldsymbol{\psi}))
- (\pressure, \nabla \cdot \boldsymbol{\psi})
- (\nabla \cdot \convection, \chi) \big\}\; \mathrm{d} t 
\\
& & 
+ (\convection(0),\boldsymbol{\psi}(0))\,,\\
F(\boldsymbol{\psi}) &:=&  
\displaystyle
\int_{I} \big\{
(\stokesforce, \boldsymbol{\psi})\big\}\; \mathrm{d} t
+ (\convection_0,\boldsymbol{\psi}(0))
\,.
\end{array}
\end{displaymath}

\subsection{Discretization in Time}
\label{sec:2:4}

The sets of time subintervals $I_n$ and $I_n^\textnormal{F}$ as introduced in 
Sec. \ref{sec:2:2:multirate}
are finite and countable. Therefore, the separation of the global space-time
cylinder $Q=\Omega \times I$ into a partition of space-time slabs
$\hat{Q}_n = \Omega \times I_n$ for the transport problem and 
$\hat{Q}_n^\textnormal{F} = \Omega \times I_n^\textnormal{F}$ for the Stokes flow 
problem, respectively, is reasonable.
The time domain of each space-time slab $\hat{Q}_n$ or $\hat{Q}_n^\textnormal{F}$ 
is then discretized using a one-dimensional triangulation $\mathcal{T}_{\tau,n}$ 
or $\mathcal{T}_{\sigma,n}$ for the closure of the subinterval 
$\bar{I}_n=[t_{n-1}^\textnormal{T}, t_n^\textnormal{T}]$ or
$\bar{I}_n^\textnormal{F}=[t_{n-1}^\textnormal{F}, t_n^\textnormal{F}]$, respectively.
This allows to have more than one cell in time on a slab $\hat{Q}_n$ or 
$\hat{Q}_n^\textnormal{F}$ and a different number of cells in time of pairwise 
different slabs $\hat{Q}_i$ and $\hat{Q}_j$ or
$\hat{Q}_i^\textnormal{F}$ and $\hat{Q}_j^\textnormal{F}$, $0 < i,j \le N^\ell,
N^{\textnormal{F},\ell}$, for the $\ell$-th adaptivity loop.
Furthermore, let $\mathcal{F}_\tau$ and $\mathcal{F}_\sigma$ be the sets of all 
interior time points given as
\begin{displaymath}
\begin{array}{rcl}
\mathcal{F}_\tau &:=& ( \{ t_1^\textnormal{T}, \dots, t_{N^\ell}^\textnormal{T} \}
\cup
\{ t \in \partial K_n \mid K_n \in \mathcal{T}_{\tau,n} \} ) \setminus \{ 0, T \}
\\[1.5ex]
\mathcal{F}_\sigma &:=& ( \{ t_1^\textnormal{F}, \dots, 
t_{N^{\textnormal{F},\ell}}^\textnormal{F} \}
\cup
\{ t \in \partial K_n \mid K_n \in \mathcal{T}_{\sigma,n} \} ) \setminus \{ 0, T \}
\end{array}
\end{displaymath}
with $1 \le n \le N^\ell,N^{\textnormal{F},\ell}$.
The commonly used time step size $\tau_K$ or $\sigma_K$ is here the diameter or 
length of the cell in time of $\mathcal{T}_{\tau,n}$ or $\mathcal{T}_{\sigma,n}$ 
and the global time discretization parameter $\tau$ or $\sigma$ is the maximum 
time step size $\tau_K$ or $\sigma_K$ of all cells in time of all slabs 
$\hat{Q}_n$ or $\hat{Q}_n^\textnormal{F}$, $0 < n \le N^\ell,N^{\textnormal{F},\ell}$.

For the discretization in time of the transport problem \eqref{eq:5:weak_transport_problem}
we use a discontinuous Galerkin method dG($r$) with an arbitrary polynomial degree
$r \ge 0$. Let $X_{\tau}^{\textnormal{dG}(r)}$ be the time-discrete function space
given as
\begin{equation}
\label{eq:7:Def_X_tau_dGr}
\begin{array}{rcl}
X_{\tau}^{\textnormal{dG}(r)} &:= \Big\{ &
 \concentration_{\tau} \in L^2(0, T; H_0^1(\Omega)) \,\,\big|\,\,
 \concentration_{\tau}|_{K_n} \in \mathcal{P}_r(K_n; H_0^1(\Omega))\,,\\[1.5ex]
 &~& K_n \in \mathcal{T}_{\tau,n}\,,\,\,
 n=1,\dots,N^\ell\,,\,\,
 \concentration_{\tau}(0)\in L^2(\Omega)
\Big\}\,,
\end{array}
\end{equation}
where $\mathcal{P}_{r}(K_n; H_0^1(\Omega))$ denotes the space of all polynomials
in time up to degree $r \ge 0$ on $K_n \in \mathcal{T}_{\tau,n}$
with values in $H_0^1(\Omega)$\,.
For some discontinuous in time function
$u_{\tau} \in X_{\tau}^{\textnormal{dG}(r)}$
we define the limits $u_{\tau}(t_F^\pm)$
from above and below of $u_{\tau}$ at $t_F$ as well as their jump at $t_F$ by
\begin{displaymath}
u_\tau(t_F^\pm) := \displaystyle \lim_{t \to t_F \pm 0} u_\tau(t)\,,\quad
[ u_\tau ]_{t_F} := u_\tau(t_F^+) - u_\tau(t_F^-) \,.
\end{displaymath}

The semidiscretization in time of the the transport problem
\eqref{eq:5:weak_transport_problem} then reads as follows:

\textit{For a given $\convection_{\sigma} \in Y_{\sigma}^{\textnormal{dG}(r)}$
of \eqref{eq:12:B_sigma_F_sigma},
find $\concentration_\tau \in X_{\tau}^{\textnormal{dG}(r)}$ such that
}
\begin{equation}
\label{eq:8:A_tau_u_phi_eq_F_phi}
A_{\tau}(\concentration_\tau; \convection)(\varphi_\tau) = G_\tau(\varphi_\tau)
\quad \forall \varphi_\tau \in X_{\tau}^{\text{dG}(r)}\,,
\end{equation}
where the semi-discrete bilinear form and linear form are given by
\begin{equation}
\label{eq:9:Def_A_tau_u_phi_and_Def_F_phi_tau}
\begin{array}{rcl}
A_{\tau}(\concentration_\tau; \convection)(\varphi_\tau) &:=& \displaystyle
\sum_{n=1}^{N^\ell}  \sum_{K_n \in \mathcal{T}_{\tau,n}}
  \int_{K_n} \big\{(
  \partial_{t} \concentration_\tau,\varphi_\tau)
  + a(\concentration_\tau; \convection_\sigma)(\varphi_\tau)
  \big\} \mathrm{d} t\;\\[1.5ex]
&~& + (
\concentration_{\tau}(0^+), \varphi_{\tau}(0^+))
+ \displaystyle \sum_{t_F \in \mathcal{F}_\tau}
(
[\concentration_\tau]_{t_F}, \varphi_\tau(t_F^+) )\,,\\[3ex]
G_\tau (\varphi_{\tau}) &:=&
\displaystyle \int_I (\transportforce, \varphi_{\tau})\;\mathrm{d}t
+ (\concentration_{0}, \varphi_{\tau}(0^+))\,,
\end{array}
\end{equation}
with the bilinear form $a(\cdot,\cdot)(\cdot)$ depending on the semi-discrete 
Stokes solution $\convection_\sigma$.
\begin{remark}
\label{rem:1:Galerkin_orthogonality_time}
For the error $e = \concentration-\concentration_{\tau}$
we get by subtracting Eq.~\eqref{eq:8:A_tau_u_phi_eq_F_phi} from
Eq.~\eqref{eq:5:weak_transport_problem} the identity
\begin{equation}
\label{eq:10:Galerkin_orthogonality_time}
\begin{aligned}
&\displaystyle
\sum_{n=1}^{N^\ell}  \sum_{K_n \in \mathcal{T}_{\tau,n}}
\int_{K_n}
\big\{
(
\partial_{t} e,\varphi_{\tau})
+ a(e, \convection_{\sigma})(\varphi_{\tau})
\big\}
\mathrm{d} t\;
\\
&  +  \displaystyle\sum_{t_F \in \mathcal{F}_\tau}
(
[e]_{t_F},\varphi_{\tau}(t_F^+))
+ (e(0^+),\varphi_{\tau}(0^+))
\\
&=
- \sum_{n=1}^{N^\ell}\sum_{K_n \in \mathcal{T}_{\tau,n}}\int_{K_n}
\big(
(\convection-\convection_{\sigma}) \cdot \nabla \concentration,
\varphi_{\tau}
\big)
\mathrm{d} t\,,
\end{aligned}
\end{equation}
with a non-vanishing right-hand side term depending on the stabilization and the
error in the approximation of the flow field.
Eq.~\eqref{eq:10:Galerkin_orthogonality_time} with the perturbation term on
the right-hand side replaces the standard Galerkin orthogonality of the
space-time finite element approximation.
\end{remark}

The discontinuous time-discrete function space for the Stokes flow problem is given by
\begin{equation}
\label{eq:11:Def_Y_sigma_dGr}
\begin{array}{rcl}
Y_{\sigma}^{\textnormal{dG}(r)} &:= \Big\{ &
\{ \convection_{\sigma},\pressure_\sigma \} \in 
L^2(0, T; H_0^1(\Omega)^d\times L_0^2(\Omega)) \,\,\big|\\[1.5ex]
&~& \convection_{\sigma}|_{K_n} \in \mathcal{P}_r(K_n; H_0^1(\Omega)^d)\,,
\convection_{\sigma}(0)\in L^2(\Omega)\,,\\[1.5ex]
 &~& \pressure_{\sigma}|_{K_n} \in \mathcal{P}_r(K_n; L_0^2(\Omega))\,, 
 K_n \in \mathcal{T}_{\sigma,n}\,,\,\,
 n=1,\dots,N^{\textnormal{F},\ell}
\Big\}\,.
\end{array}
\end{equation}
Then, the semidiscretization in time of the the Stokes flow problem
\eqref{eq:6:weak_stokes} reads as follows:

\textit{Find $\{ \convection_{\sigma},\pressure_\sigma \} \in 
Y_{\sigma}^{\textnormal{dG}(r)}$
such that
}
\begin{equation}
\label{eq:12:B_sigma_F_sigma}
B_{\sigma}(\convection_\sigma,\pressure_\sigma)(\boldsymbol{\psi}_\sigma,\chi_\sigma)
= F_\sigma(\boldsymbol{\psi}_\sigma)
\quad
\forall \{\boldsymbol{\psi}_{\sigma}, \chi_{\sigma h}\} \in Y_{\sigma}^{\text{dG}(r)}\,,
\end{equation}
where the semi-discrete bilinear form and linear form are given by
\begin{equation}
\label{eq:13:Def_B_sigma_F_sigma}
\begin{array}{rcl}
B_{\sigma}(\convection_\sigma,\pressure_\sigma)(\boldsymbol{\psi}_\sigma,\chi_\sigma)
&:=& \displaystyle
\sum_{n=1}^{N^{\textnormal{F},\ell}}  \sum_{K_n \in \mathcal{T}_{\sigma,n}}
\int_{K_n} \big\{(\partial_{t} \convection_\sigma,\boldsymbol{\psi}_\sigma)
+ (2 \viscosity \boldsymbol\epsilon(\convection_\sigma), 
\boldsymbol\epsilon(\boldsymbol{\psi_\sigma}))
\\[1.5ex]
&~& - (\pressure_\sigma, \nabla \cdot \boldsymbol{\psi_\sigma})
-(\nabla \cdot \convection_\sigma, \chi_\sigma)
  \big\} \mathrm{d} t 
\\[1.5ex]
&~&
+ (\convection_{\sigma}(0^+), \boldsymbol{\psi}_\sigma(0^+))
+ \displaystyle \sum_{t_F \in \mathcal{F}_\sigma}
([\convection_\sigma]_{t_F}, \boldsymbol{\psi}_\sigma(t_F^+) )\,,\\[3ex]
F_\sigma (\varphi_{\tau}) &:=&
\displaystyle \int_I (\stokesforce, \boldsymbol{\psi_\sigma})\;\mathrm{d}t
+ (\convection_{0}, \boldsymbol{\psi}_\sigma(0^+))\,.
\end{array}
\end{equation}

\subsection{Discretization in Space and SUPG Stabilization}
\label{sec:2:5}

Next, we describe the Galerkin finite element approximation in space of the
semi-discrete transport problem \eqref{eq:8:A_tau_u_phi_eq_F_phi} and
the flow problem \eqref{eq:12:B_sigma_F_sigma}, respectively.
%
We use Lagrange type finite element spaces of continuous functions that are
piecewise polynomials. For the discretization in space, we consider a
separation $Q_n=\mathcal{T}_{h,n}\times I_n$ or 
$Q_n^{\textnormal{F}}=\mathcal{T}_{h,n}^{\textnormal{F}}\times I_n^{\textnormal{F}}$, 
where $\mathcal{T}_{h,n}$ or $\mathcal{T}_{h,n}^{\textnormal{F}}$ build a
decomposition of the domain $\Omega$ into disjoint elements $K$ or 
$K^{\textnormal{F}}$, such that $\overline{\Omega}=\cup_{K\in\mathcal{T}_{h}}\overline{K}$ 
or $\overline{\Omega}=
\cup_{K^{\textnormal{F}}\in\mathcal{T}_{h}^{\textnormal{F}}}\overline{K^{\textnormal{F}}}$
for the transport and Stokes flow problem, respectively.
Here, we choose the elements $K\in\mathcal{T}_{h}$ or 
$K^{\textnormal{F}}\in\mathcal{T}_{h}^{\textnormal{F}}$ to be quadrilaterals
for $d=2$ and hexahedrals for $d=3$.
We denote by $h_{K}$ or $h_{K}^{\textnormal{F}}$ the diameter of the element 
$K$ or $K^{\textnormal{F}}$. The global space
discretization parameter $h$ or $h^{\textnormal{F}}$ is given by 
$h:=\max_{K\in\mathcal{T}_{h}}h_{K}$ or 
$h^{\textnormal{F}}:=\max_{K^{\textnormal{F}}\in
\mathcal{T}_{h}^{\textnormal{F}}}h_{K}^{\textnormal{F}}$, respectively.
Our mesh adaptation process yields locally refined cells, which is
enabled by using hanging nodes. We point out that
the global conformity of the finite element approach is preserved since the
unknowns at such hanging nodes are eliminated by interpolation between the
neighboring 'regular' nodes; cf.~\cite[Chapter 4.2]{Bangerth2003} and
\cite{Carey1984} for more details.
On $\mathcal{T}_{h}$ and $\mathcal{T}_{h}^{\textnormal{F}}$ we define the 
discrete finite element spaces by
$
V_{h}^{p,n}:=
\big\{v\in C(\overline{\Omega})\mid v_{|K}
\in Q_h^p(K)\,,\forall K\in\mathcal{T}_{h},
\big\}\,,
$
and
$
V_{h^{\textnormal{F}}}^{p,n}:=
\big\{v\in C(\overline{\Omega})\mid v_{|K^{\textnormal{F}}}
\in Q_h^p(K^{\textnormal{F}})\,,\forall K^{\textnormal{F}}\in
\mathcal{T}_{h}^{\textnormal{F}},
\big\}\,,
$
with $1 \le n \le N,N^{\textnormal{F}}$, where $Q_h^p(K)$ or 
$Q_{h^{\textnormal{F}}}^p(K^{\textnormal{F}})$ is the space defined on the 
reference element with maximum degree $p$ in each variable. 
By replacing $H_0^1(\Omega)$ in the
definition of the semi-discrete function space $X_{\tau}^{\textnormal{dG}(r)}$
in \eqref{eq:7:Def_X_tau_dGr} by $V_h^{p,n}$ and by replacing 
$H_0^1(\Omega)^d, L_0^2(\Omega)$ in the definition of the semi-discrete function 
space $Y_{\sigma}^{\textnormal{dG}(r)}$ in \eqref{eq:11:Def_Y_sigma_dGr} by 
$V_{h^{\textnormal{F}}}^{p,n}$, we obtain the fully discrete function spaces for 
the transport and Stokes flow problem, respectively, 
\begin{equation}
\begin{array}{rcl}
\label{eq:14:Def_X_tau_h_dGr_p_Y_sigma_h}
X_{\tau h}^{\text{dG}(r),p} := & \Big\{ &
\concentration_{\tau h}\in X_{\tau}^{\text{dG}(r)} \,\,\big|\,\,
\concentration_{\tau h}|_{K_n} \in \mathcal{P}_r(K_n;H_h^{p_{\concentration},n})
\,,\\
&~& \concentration_{\tau h}(0) \in H_h^{p_{\concentration},0},
K_n \in \mathcal{T}_{\tau,n}\,,\,\,n=1,\dots,N
\Big\}\,,
\\[1.5ex]
Y_{\sigma h^{\textnormal{F}}}^{\text{dG}(r),p} := & \Big\{ &
\{ \convection_{\sigma h},\pressure_{\sigma h} \} \in 
Y_{\sigma}^{\text{dG}(r)} \,\,\big|\,\,
\convection_{\sigma h}|_{K_n} \in \mathcal{P}_r(K_n; (H_h^{p_\convection,n})^d)\,,
\\[0.5ex]
 &~& \convection_{\sigma h}(0)\in (H_h^{p_\convection,0})^d\,,
 \pressure_{\sigma h}|_{K_n} \in \mathcal{P}_r(K_n; L_h^{p_\pressure,n})\,,\\ 
 &~& K_n \in \mathcal{T}_{\sigma,n}\,,\,\,
 n=1,\dots,N^{\textnormal{F},\ell}
\Big\}
\end{array}
\end{equation}
\begin{displaymath}
H_h^{p_{\concentration},n}:=V_h^{p_{\concentration},n}\cap H_0^1(\Omega), \quad
H_h^{p_{\convection},n}:=V_h^{p_{\convection},n}\cap H_0^1(\Omega), \quad
L_h^{p_{\pressure},n}:=V_h^{p_{\pressure},n}\cap L_0^2(\Omega).
\end{displaymath}
We note that the spatial finite
element space $V_h^{p,n}$ and $V_{h^{\textnormal{F}}}^{p,n}$ are allowed to be 
different on all subintervals $I_n$ and $I_n^{\textnormal{F}}$, respectively,
which is natural in the context of a discontinuous Galerkin approximation of the 
time variable and allows dynamic mesh changes in time.
Due to the conformity of $H_h^{p_{\concentration},n}$, $H_h^{p_{\convection},n}$
and $L_h^{p_{\pressure},n}$, we get
$X_{\tau h}^{\textnormal{dG}(r),p}\subseteq X_{\tau}^{\textnormal{dG}(r)}$ and
$Y_{\sigma h^{\textnormal{F}}}^{\textnormal{dG}(r),p}\subseteq 
Y_{\sigma}^{\textnormal{dG}(r)}$, respectively.

For convection-dominated transport, the finite element approximation needs to be
stabilized in order to avoid spurious and non-physical oscillations of
the discrete solution arising close to sharp fronts and layers. Here, we apply
the streamline upwind Petrov-Galerkin (SUPG) method introduced by Hughes and Brooks
\cite{Hughes1979,Brooks1982}.
With this in mind,the stabilized fully discrete discontinuous in time scheme for 
the transport problem reads as follows:

\textit{For a given $\convection_{\sigma h} \in 
Y_{\sigma h^{\textnormal{F}}}^{\text{dG}(r),p}$ of \eqref{eq:18:B_v_p_psi_chi_eq_F_psi},
find $\concentration_{\tau h} \in X_{\tau h}^{\textnormal{dG}(r),p}$ such that}
\begin{equation}
\label{eq:15:A_S_G_tauh}
 A_{S}(\concentration_{\tau h}; \convection_{\sigma h})(\varphi_{\tau h})
=
G_{\tau}(\varphi_{\tau h})
\quad \forall \varphi_{\tau h} \in X_{\tau h}^{\text{dG}(r),p}\,,
\end{equation}
\textit{where the linear form} $G_{\tau}(\cdot)$ \textit{is 
defined in \eqref{eq:9:Def_A_tau_u_phi_and_Def_F_phi_tau} and the stabilized 
bilinear form}
$A_{S}(\cdot;\cdot)(\cdot)$ \textit{is given by}
\begin{displaymath}
A_{S}(u_{\tau h};\convection_{\sigma h})(\varphi_{\tau h})
:=A_{\tau}(u_{\tau h};\convection_{\sigma h})(\varphi_{\tau h})
+S_A(u_{\tau h};\convection_{\sigma h})(\varphi_{\tau h})\,,
\end{displaymath}
\textit{with} $A_{\tau}(\cdot;\cdot)(\cdot)$ \textit{being defined in 
\eqref{eq:9:Def_A_tau_u_phi_and_Def_F_phi_tau}.
Here, the SUPG stabilized bilinear form} 
$S_A(\cdot;\cdot)(\cdot)$
\textit{is defined by}
\begin{equation}
\label{eq:16:Def_S_A}
\begin{array}{r@{\,}c@{\,}l@{\quad}}
\displaystyle
S_A(u_{\tau h};\convection_{\sigma h})(\varphi_{\tau h}) 
& := & 
\displaystyle
\sum_{n=1}^{N^\ell}  \sum_{K_n \in \mathcal{T}_{\tau,n}}
  \int_{K_n}
\sum\limits_{K\in \mathcal{T}_h}\delta_K\big( r(u_{\tau h}), 
\convection_{\sigma h} \cdot \nabla \varphi_{\tau h}\big)_K \,\mathrm{d} t 
\\[3.5ex]
& &
\displaystyle
+ \displaystyle\sum_{t_F \in \mathcal{F}_\tau}\sum\limits_{K\in\mathcal{T}_h} 
\delta_K
\big(
\left[u_{\tau h}\right]_{t_F}, 
\convection_{\sigma h} \cdot \nabla \varphi_{\tau h}(t_F^+)\big)_{K} 
\\[3.5ex]
& &
\displaystyle
+ \displaystyle\sum\limits_{K\in\mathcal{T}_h}
\delta_K \big(
u_{\tau h,0}^{+} - u_0,
\convection_{\sigma h} \cdot \nabla \varphi_{\tau h,0}^{+} \big)_{K}\,, 
\end{array}
\end{equation}
\textit{where} $\delta_K$ \textit{is the so-called stabilization parameter 
and the residual term} $r(\cdot;\cdot)$ \textit{is given by}
\begin{displaymath}
r(u_{\tau h}):=
\partial_{t} u_{\tau h} 
- \nabla\cdot\left(\varepsilon\nabla u_{\tau h}\right)  
+ \convection_{\sigma h} \cdot \nabla u_{\tau h} 
+ \alpha u_{\tau h} - g\,.
\end{displaymath}
We note that the bilinear form $a(\cdot;\convection_{\sigma h})(\cdot)$ occurring 
in $A_S(\cdot;\cdot)(\cdot)$ reads here as
\begin{displaymath}
a(\concentration_{\tau h}; \convection_{\sigma h})(\varphi_{\tau h}) =
(\varepsilon \nabla \concentration_{\tau h}, \nabla \varphi_{\tau h})
+(\convection_{\sigma h} \cdot \nabla \concentration_{\tau h}, \varphi_{\tau h})
+ (\alpha \concentration_{\tau h},\varphi_{\tau h})
\end{displaymath}
for the fully discrete solutions.
\begin{remark}
\label{rem:2:SUPGstabilization}
The proper choice of the stabilization parameter $\delta_K$ is an important
issue in the application of the SUPG approach; cf., e.g.,
\cite{John2011,John2009,John2018}
and the discussion therein. For time-dependent convection-diffusion-reaction
problems an optimal error estimate for $\delta_K=\mathrm{O}(h)$ is derived
in \cite{John2011}.
\end{remark}

\begin{remark}
\label{rem:3:GalerkinOrthogonality_space}
For the error $e = \concentration_{\tau}-\concentration_{\tau h}$
we get by subtracting Eq.~\eqref{eq:15:A_S_G_tauh} from
Eq.~\eqref{eq:8:A_tau_u_phi_eq_F_phi} the identity
\begin{equation}
\label{eq:17:Galerkin_orthogonality_space}
\begin{aligned}
&\displaystyle
\sum_{n=1}^{N^\ell}  \sum_{K_n \in \mathcal{T}_{\tau,n}}
\int_{K_n}
\big\{
(
\partial_{t} e,\varphi_{\tau h})
+ a(e; \convection_{\sigma h})(\varphi_{\tau h})
\big\}
\mathrm{d} t\;\\
&  +  \displaystyle\sum_{t_F \in \mathcal{F}_\tau}
(
[e]_{t_F},\varphi_{\tau h}(t_F^+)
+ (e(0^+),\varphi_{\tau h}(0^+))
\\
&=
S_A(\concentration_{\tau h}; \convection_{\sigma h})(\varphi_{\tau h})
- \sum_{n=1}^{N^\ell}\sum_{K_n \in \mathcal{T}_{\tau,n}}\int_{K_n}
\big(
(\convection_{\sigma}-\convection_{\sigma h}) \cdot \nabla \concentration_{\tau},
\varphi_{\tau h}
\big)
\mathrm{d} t\,,
\end{aligned}
\end{equation}
with a non-vanishing right-hand side term depending on the stabilization and the
error in the approximation of the flow field.
Eq.~\eqref{eq:17:Galerkin_orthogonality_space} with the perturbation term on
the right-hand side replaces the standard Galerkin orthogonality of the
space-time finite element approximation.
\end{remark}

Finally, the fully discrete discontinuous in time scheme for the Stokes flow 
problem reads as follows:

\textit{Find $\{ \convection_{\sigma h},\pressure_{\sigma h} \} \in 
Y_{\sigma h^{\textnormal{F}}}^{\textnormal{dG}(r),p}$
such that}
\begin{equation}
\label{eq:18:B_v_p_psi_chi_eq_F_psi}
B_{\sigma}(\convection_{\sigma h},\pressure_{\sigma h})
(\boldsymbol{\psi}_{\sigma h},\chi_{\sigma h})
= F_{\sigma}(\boldsymbol{\psi}_{\sigma h})
\quad
\forall \{\boldsymbol{\psi}_{\sigma h}, \chi_{\sigma h}\} \in 
Y_{\sigma h^{\textnormal{F}}}^{\textnormal{dG}(r),p}\,,
\end{equation}
\textit{with $B_{\sigma}(\cdot,\cdot)(\cdot,\cdot)$ and $F_{\sigma}(\cdot)$ 
being defined in \eqref{eq:13:Def_B_sigma_F_sigma}.}

\section{An A Posteriori Error Estimator for the Transport Problem}
\label{sec:3:dwr}

In this section we derive a DWR-based a posteriori error representation for the 
stabilized transport problem \eqref{eq:15:A_S_G_tauh} coupled with the 
flow problem via the convection tensor $\convection_{\sigma h}$ given by 
Eq.~\eqref{eq:18:B_v_p_psi_chi_eq_F_psi}.
Since the derivation is close to our work based on a coupling of a steady-state 
Stokes problem, we keep this section rather short by drawing attention only to 
the differences and refer to our work \cite{Bause2021} for a detailed version of
the proofs and further details.

Here, only goal quantities depending on the unknown $\concentration$ are
studied. For applications of practical interest, physical quantities in terms of
the transport quantity $u$ are typically of higher relevance than quantities in
the unknowns $\convection$ and $\pressure$ of the flow problem. In the sequel,
we introduce this goal quantity with the following properties.
\\
\textbf{Assumption (Target functional $J$)}

\emph{Let us assume $J:X\rightarrow \mathbb{R}$ to be a linear functional 
representing the goal quantity of physical interest. In general, 
this functional is given as}
\begin{equation}
\label{eq:19:Def_J_u}
J(u)=\int_0^T J_1(u(t))\;\mathrm{d}t + J_2(u(T))\,,
\end{equation}
\emph{where $J_1 \in L^2(I;H^{-1}(\Omega))$ and $J_2 \in H^{-1}(\Omega)$
are three times differentiable 
functionals defining the dual right-hand side and the dual initial at time
$t=T$, respectively, where each of them may be zero.}

Since we aim at controlling the respective errors due to the discretization in 
time as well as in space, we split the a posteriori error representation with 
respect to $J$ into the contributions
\begin{equation}
\label{eq:20:J_u_J_u_tauh_eq_J_u_J_u_tau_J_u_tauh}
 J(\concentration)-J(\concentration_{\tau h}) =
 J(\concentration)-J(\concentration_{\tau})
 + J(\concentration_{\tau})-J(\concentration_{\tau h})\,.
\end{equation}

For the respective error representations we define the Lagrangian
functionals
$\mathcal{L}: X\times X \rightarrow \mathbb{R}$,
$\mathcal{L}_\tau: X_{\tau}^{\textnormal{dG}(r)} \times X_{\tau}^{\textnormal{dG}(r)}
\rightarrow \mathbb{R}$, and
$\mathcal{L}_{\tau h}:
X_{\tau h}^{\textnormal{dG}(r),p} \times X_{\tau h}^{\textnormal{dG}(r),p}
\rightarrow \mathbb{R}$ by
\begin{subequations}
\label{eq:23:Def_L_u_z_Def_L_tau_u_z_Def_L_tau_h_u_z}
\begin{align}
\label{eq:21a:Def_L_u_z}
\mathcal{L}(\concentration,\dualz;\convection) & :=  J(\concentration)
+ G(\dualz)
- A(\concentration; \convection)(\dualz)\,,
\\
\label{eq:21b:Def_L_tau_u_z}
\mathcal{L}_{\tau}(\concentration_\tau,\dualz_\tau;\convection_{\sigma}) & :=
J(\concentration_{\tau}) + G_\tau(\dualz_{\tau})
- A_{\tau}(\concentration_{\tau}; \convection_{\sigma})(\dualz_{\tau})\,,
\\
\label{eq:21c:Def_L_tau_h_u_z}
\mathcal{L}_{\tau h}(\concentration_{\tau h},\dualz_{\tau h};\convection_{\sigma h}) 
& :=
J(\concentration_{\tau h})
+ G_\tau (\dualz_{\tau h})
- A_S(\concentration_{\tau h}; \convection_{\sigma h})(\dualz_{\tau h})\,.
\end{align}
\end{subequations}
Here, the Lagrange multipliers $\dualz$, $\dualz_\tau,$ and $\dualz_{\tau h}$
are called dual variables in contrast to the primal variables
$\concentration$, $\concentration_\tau,$ and
$\concentration_{\tau h}$; cf. \cite{Besier2012,Becker2001}.
Considering the directional derivatives of the Lagrangian functionals, also
known as
G\^{a}teaux derivatives, with respect to their first argument, i.e.\
\begin{displaymath}
\mathcal{L}^{\prime}_{\concentration}(\concentration,\dualz;\convection)(\varphi) :=
\lim_{t\neq0,t\rightarrow 0}
t^{-1}\big\{\mathcal{L}(\concentration + t\varphi,\dualz;\convection)
-\mathcal{L}(\concentration,\dualz;\convection)\big\},
\quad \varphi \in X\,,
\end{displaymath}
leads to the so-called dual problems: 
Find the continuous dual solution $\dualz \in X$, the semi-discrete
dual solution $\dualz_{\tau} \in X_{\tau}^{\textnormal{dG}(r)}$ and the fully
discrete dual solution $\dualz_{\tau h} \in X_{\tau h}^{\textnormal{dG}(r),p}$,
respectively, such that
\begin{subequations}
\label{eq:26:DualProblems}
\begin{align}
\label{eq:22a:DualProblems_continuous}
A^{\prime}(\concentration; \convection)(\varphi,\dualz)
&=
J^{\prime}(\concentration)(\varphi)
\quad 
\forall \varphi \in X\,,\\
\label{eq:22b:DualProblems_semidiscrete}
A_{\tau}^{\prime}(\concentration_{\tau}; \convection_{\sigma})(\varphi_{\tau},\dualz_{\tau})
& =
J^{\prime}(\concentration_{\tau})(\varphi_{\tau})
\quad 
\forall \varphi_{\tau}\in X_{\tau}^{\text{dG}(r)}\,,
\\
\label{eq:22c:DualProblems_fullydiscrete}
A_{S}^{\prime}(\concentration_{\tau h}; \convection_{\sigma h})(\varphi_{\tau h},\dualz_{\tau h})
& =
J^{\prime}(\concentration_{\tau h})(\varphi_{\tau h})
\quad \forall \varphi_{\tau h}\in X_{\tau h}^{\text{dG}(r),p}\,,
\end{align}
\end{subequations}
where we refer to our work \cite{Bause2021} for a detailed description of the
adjoint bilinear forms $A^{\prime},A_{\tau}^{\prime},A_S^{\prime}$ as well as 
the dual right hand side term $J^\prime$. 
\begin{remark}
\label{rem:4:PrimalProblems}
We note that the directional derivatives of the Lagrangian functionals with
respect to their second argument leads to the primal problems given by
Eqs.~\eqref{eq:5:weak_transport_problem},
\eqref{eq:8:A_tau_u_phi_eq_F_phi} and \eqref{eq:15:A_S_G_tauh},
respectively.
\end{remark}

In the following Thm.~\ref{thm:1:ErrorRepresentation} we derive error 
representation formulas in space and time for the transport problem depending on 
the residuals of the primal and dual problem as well as jump, stabilization and 
coupling terms due to a non vanishing Galerkin orthogonality described in 
Rem.~\ref{rem:1:Galerkin_orthogonality_time} and 
Rem.~\ref{rem:3:GalerkinOrthogonality_space}, respectively.
The primal and dual residuals based on the continuous and semi-discrete schemes
are defined by means of the G\^{a}teaux derivatives of the Lagrangian functionals
in the following way:
\begin{displaymath} 
\begin{array}{r@{\;}c@{\;}l@{\;}c@{\;}l@{\;}}
\rho(u;\convection)(\varphi) & := &  
\mathcal{L}_{\dualz}^{\prime}(\concentration,\dualz;\convection)(\varphi)
& = & G(\varphi)-A(\concentration;\convection)(\varphi)\,,
\\[1.5ex]
\rho^{\ast}(\concentration,\dualz;\convection)(\varphi) & := & 
\mathcal{L}_{\concentration}^{\prime}(\concentration,\dualz;\convection)(\varphi)
& = & J^{\prime}(\concentration)(\varphi)
-A^{\prime}(\concentration;\convection)(\varphi,\dualz)\,,
\\[1.5ex]
\rho_{\tau}(\concentration;\convection_{\sigma})(\varphi) & := &
\mathcal{L}_{\tau,\dualz}^{\prime}(\concentration,\dualz;\convection_{\sigma})(\varphi)
& = & G_{\tau}(\varphi)-A_{\tau}(\concentration;\convection_{\sigma})(\varphi)\,,
\\[1.5ex]
\rho_{\tau}^{\ast}(\concentration,\dualz;\convection_{\sigma})(\varphi) & := &
\mathcal{L}_{\tau,\concentration}^{\prime}(\concentration,\dualz;\convection_{\sigma})(\varphi)
& = & J^{\prime}(\concentration)(\varphi)
-A_{\tau}^{\prime}(\concentration;\convection_{\sigma})(\varphi,\dualz)\,.
\end{array}
\end{displaymath}
By using these residuals as well as the Galerkin orthogonality described in 
Rem.~\ref{rem:1:Galerkin_orthogonality_time} and 
Rem.~\ref{rem:3:GalerkinOrthogonality_space}, respectively, we get the following 
result for the DWR-based error representation in space and time for the 
transport problem.
\begin{theorem}
\label{thm:1:ErrorRepresentation}
Let $\{\concentration,\dualz\}\in X \times X$,
$\{\concentration_{\tau},\dualz_{\tau}\}
\in
X_{\tau}^{\textnormal{dG}(r)} \times X_{\tau}^{\textnormal{dG}(r)}$,
and
$\{\concentration_{\tau h},\dualz_{\tau h}\}
\in X_{\tau h}^{\textnormal{dG}(r),p} \times X_{\tau h}^{\textnormal{dG}(r),p}$
denote the stationary points of
$\mathcal{L}, \mathcal{L}_{\tau}$, and $\mathcal{L}_{\tau h}$
on the different levels of discretization, i.e.,
\begin{displaymath}
\begin{aligned}
\mathcal{L}^{\prime}(\concentration,\dualz;\convection)(\delta \concentration, \delta \dualz)
& = 0 \quad
\forall \{\delta \concentration,\delta \dualz\}\in X \times X\,,
\\
\mathcal{L}_{\tau}^{\prime}(\concentration_{\tau},\dualz_{\tau};\convection_{\sigma})
(\delta \concentration_{\tau}, \delta \dualz_{\tau})
& = 0
\quad \forall \{\delta \concentration_{\tau},\delta \dualz_{\tau}\}
\in X_{\tau}^{\text{dG}(r)} \times X_{\tau}^{\text{dG}(r)}\,,
\\
\mathcal{L}_{\tau h}^{\prime}(\concentration_{\tau h},\dualz_{\tau h};\convection_{\sigma h})
(\delta \concentration_{\tau h}, \delta \dualz_{\tau h})
& = 0
\quad \forall \{\delta \concentration_{\tau h},\delta \dualz_{\tau h}\}
\in X_{\tau h}^{\text{dG}(r),p} \times X_{\tau h}^{\text{dG}(r),p}\,.
\end{aligned}
\end{displaymath}
Additionally, for the errors $e = \concentration - \concentration_{\tau}$
and $e = \concentration_{\tau} - \concentration_{\tau h}$
we have the Eqs.~\eqref{eq:10:Galerkin_orthogonality_time} and 
\eqref{eq:17:Galerkin_orthogonality_space} of Galerkin
orthogonality type. Then, for the discretization errors in space and time we get
the representation formulas
\begin{subequations}
\label{eq:28:ERtransport}
\begin{align}
\label{eq:23a:ERtransport_time}
J(\concentration)-J(\concentration_{\tau}) & =
\frac{1}{2}\rho(\concentration_{\tau};\convection)(\dualz-\tilde{\dualz}_{\tau})
+ \frac{1}{2}\rho^{\ast}(\concentration_{\tau},\dualz_{\tau};\convection)
(\concentration-\tilde{\concentration}_{\tau})
\\
\nonumber
& \qquad
+ \frac{1}{2} \mathcal{D}_{\tau}^{\prime}(\concentration_{\tau},\dualz_{\tau})
(\tilde{\concentration}_{\tau}-\concentration_{\tau},\tilde{\dualz}_{\tau}-\dualz_{\tau})
\\
\nonumber
& \qquad
+ \mathcal{D}_{\tau}(\concentration_{\tau},\dualz_{\tau})
+ \mathcal{R}_{\tau}\,,
\\
\label{eq:23b:ERtransport_space}
J(\concentration_{\tau})-J(\concentration_{\tau h}) & =
\frac{1}{2}\rho(\concentration_{\tau h};\convection_{\sigma})(\dualz_{\tau}-\tilde{\dualz}_{\tau h})
+ \frac{1}{2}
\rho^{\ast}(\concentration_{\tau h},\dualz_{\tau h};\convection_{\sigma})
(\concentration_{\tau}-\tilde{\concentration}_{\tau h})
\\
\nonumber
& \qquad
+ \frac{1}{2} \mathcal{D}_{\tau h}^{\prime}(\concentration_{\tau h},\dualz_{\tau h})
(\tilde{\concentration}_{\tau h}-\concentration_{\tau h},\tilde{\dualz}_{\tau h}-\dualz_{\tau h})
\\
\nonumber
& \qquad
+ \mathcal{D}_{\tau h}(\concentration_{\tau h},\dualz_{\tau h}) + \mathcal{R}_{h}\,,
\end{align}
\end{subequations}
where $\mathcal{D}_{\tau}(\cdot,\cdot)$ and $\mathcal{D}_{\tau h}(\cdot,\cdot)$ 
are given by
\begin{equation}
\label{eq:24:Def_D_tau_h_phi_psi}
\begin{array}{rcl}
\mathcal{D}_{\tau}(\concentration_{\tau},\dualz_{\tau}) 
&=&
\displaystyle\sum_{t_F \in \mathcal{F}_\tau}
(
[\concentration_\tau]_{t_F},\dualz_{\tau}(t_F^+))
-\sum_{n=1}^{N^{\ell}}\sum_{K_n\in\mathcal{T}_{\tau,n}}\int_{K_n}\big(
(\convection-\convection_{\sigma})\cdot\nabla \concentration_\tau,\dualz_\tau
\big)\;\mathrm{d}t\,,\\
\mathcal{D}_{\tau h}(\varphi,\psi) &=&
\displaystyle
S_A(\concentration_{\tau h}; \convection_{\sigma h})(\dualz_{\tau h})
-\sum_{n=1}^{N^{\ell}}\sum_{K_n\in\mathcal{T}_{\tau,n}}\int_{K_n}\big(
(\convection_{\sigma}-\convection_{\sigma h})\cdot\nabla \concentration_{\tau h},
\dualz_{\tau h}
\big)\;\mathrm{d}t\,,
\end{array}
\end{equation}
and $\mathcal{D}^\prime_{\tau}(\cdot,\cdot)(\cdot,\cdot)$ and 
$\mathcal{D}^\prime_{\tau h}(\cdot,\cdot)(\cdot,\cdot)$ denoting the G\^{a}teaux 
derivatives with respect to the first and second argument and with 
$S_A(\cdot;\cdot)(\cdot)$ being defined in \eqref{eq:16:Def_S_A}.
Here,
$\{\tilde{\concentration}_{\tau},\tilde{\dualz}_{\tau}\}\in X_{\tau}^{\textnormal{dG}(r)}
\times X_{\tau}^{\textnormal{dG}(r)}$,
and
$\{\tilde{\concentration}_{\tau h},\tilde{\dualz}_{\tau h}\} \in
X_{\tau h}^{\textnormal{dG}(r),p} \times X_{\tau h}^{\textnormal{dG}(r),p}$
can be chosen arbitrarily and the remainder terms $\mathcal{R}_{\tau}$ and
$\mathcal{R}_{h}$ are of higher-order with respect to the errors
$\concentration-\concentration_\tau,\dualz-\dualz_\tau$ and
$\concentration_{\tau}-\concentration_{\tau h},\dualz_{\tau}-\dualz_{\tau h}$,
respectively.
\end{theorem}
%
\begin{remark}
\label{rem:3:3:AdditionalCouplingTerms}
We note that within the temporal error representation formula 
\eqref{eq:23a:ERtransport_time} additional terms due to the coupling 
occur, cf. Eq.~\eqref{eq:24:Def_D_tau_h_phi_psi}.
This is an extension of our previous results obtained in \cite{Bruchhaeuser2020}
and \cite{Bause2021}.
Furthermore, we indicate that the occurring differences $v-\tilde{v}_{\tau}$ and
$v_{\tau}-\tilde{v}_{\tau h}$ with regard to the primal and dual variables are 
called temporal and spatial weights, respectively.
\end{remark}
%
\begin{mproof}
The technique to prove the temporal error representation formula 
\eqref{eq:23a:ERtransport_time} is equivalent to the spatial counterpart that 
can be found in our work \cite[Thm. 3.1]{Bause2021} and was originally 
proved by Besier and Rannacher applied to the incompressible Navier-Stokes 
equations in \cite[Thm.~5.2]{Besier2012}.
More precisely, we are using a general result given in \cite[Lemma~5.1]{Besier2012}
with the following settings:
\begin{displaymath}
L = \mathcal{L}\,,\;\;
\tilde{L}=\mathcal{L}_{\tau}\,,\;\;
Y_1 = X \times X\,,\;\;
Y_2 = X_{\tau}^{\text{dG}(r)} \times X_{\tau}^{\text{dG}(r)}\,,\;\;
Y:=Y_1+Y_2\,,
\end{displaymath}
where $\mathcal{L},\mathcal{L}_{\tau}$ are the Lagrangian functional given by
Eq.~\eqref{eq:23:Def_L_u_z_Def_L_tau_u_z_Def_L_tau_h_u_z} and $Y,Y_1$ and $Y_2$
are function spaces defined in \cite[Lemma~5.1]{Besier2012}.
\end{mproof}

\section{Implementation of Tensor-Product Spaces}
\label{sec:4:implementation}

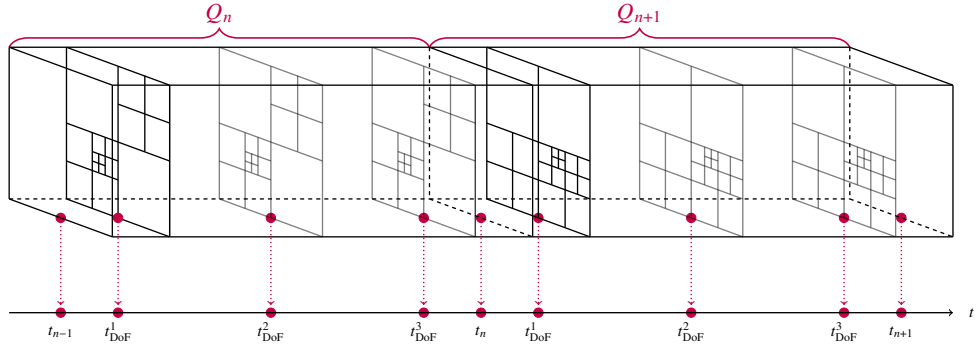
\begin{figure}[thbp]
\centering


\resizebox{.9\linewidth}{!}{%
\begin{tikzpicture}
\tikzstyle{ns1}=[line width=1.]
\tikzstyle{ns2}=[line width=1.,opacity=.5]
\tikzstyle{ns3}=[line width=1.,opacity=.5,dashed]
\tikzstyle{ns4}=[line width=0.9]

%

%
\draw[ns1] (-2.7,1) -- (0,0);
\draw[ns1] (0,0) -- (0,4);
\draw[ns1] (0,4) -- (-2.7,5);
\draw[ns1] (-2.7,5) -- (-2.7,1);
%
\draw[ns1] (-2.7,3) -- (0,2);
\draw[ns1] (-1.35,0.5) -- (-1.35,4.5);
\draw[ns1] (-2.7,2) -- (-1.35,1.5);
\draw[ns1] (-2.025,0.75) -- (-2.025,2.75);
\draw[ns1] (-2.025,2.25) -- (-1.35,2.);
\draw[ns1] (-1.6875,1.625) -- (-1.6875,2.625);
\draw[ns1] (-2.025,2.) -- (-1.6875,1.875);
\draw[ns1] (-1.8562,1.7) -- (-1.8562,2.2);
\draw[ns1] (-1.35,3.5) -- (0,3.);
\draw[ns1] (-0.65,2.25) -- (-0.65,4.25);

\node[] at (-1.35,-2.5) {\large $t_{\textnormal{DoF}}^{1}$};
\draw[HSUred,fill=HSUred](-1.35,0.5)circle(3.5pt);
\draw[->,ns1, HSUred, dotted] (-1.35,0.5) -- (-1.35,-1.8);
\draw[HSUred,fill=HSUred](-1.35,-2)circle(3.5pt);

\node[] at (-2.85,-2.5) {\large $t_{n-1}$};
\draw[HSUred,fill=HSUred](-2.85,0.5)circle(3.5pt);
\draw[->,ns1, HSUred, dotted] (-2.85,0.5) -- (-2.85,-1.8);
\draw[HSUred,fill=HSUred](-2.85,-2)circle(3.5pt);

\node[] at (8.15,-2.5) {\large $t_{n}$};
\draw[HSUred,fill=HSUred](8.15,0.5)circle(3.5pt);
\draw[->,ns1, HSUred, dotted] (8.15,0.5) -- (8.15,-1.8);
\draw[HSUred,fill=HSUred](8.15,-2)circle(3.5pt);

\draw[very thick,HSUred,decorate,decoration={brace,amplitude=14pt}]
(-4.2,5) -- (6.8,5) node[midway, above,yshift=14pt,]{\LARGE $Q_{n}$};

\draw[ns1] (-4.2,1) -- (-1.5,0);
\draw[ns1] (-1.5,0) -- (-1.5,4);
\draw[ns1] (-1.5,4) -- (-4.2,5);
\draw[ns1] (-4.2,5) -- (-4.2,1);

\draw[ns1] (-4.2,5) -- (-2.7,5);
\draw[ns1] (-1.5,4) -- (0,4);
\draw[ns1, dashed] (-4.2,1) -- (-2.7,1);
\draw[ns1] (-1.5,0) -- (0,0);

\draw[ns1] (0,0) -- (8,0);
\draw[ns1] (8,4) -- (0,4);
\draw[ns1] (-2.7,5) -- (5.3,5); 
\draw[ns1, dashed] (-2.7,1) -- (5.3,1); 

\draw[ns1,dashed] (6.8,1) -- (9.5,0);
\draw[ns1] (9.5,0) -- (9.5,4);
\draw[ns1] (9.5,4) -- (6.8,5);
\draw[ns1,dashed] (6.8,5) -- (6.8,1);

\draw[ns1] (5.3,5) -- (6.8,5);
\draw[ns1] (8,4) -- (9.5,4);
\draw[ns1, dashed] (5.3,1) -- (6.8,1);
\draw[ns1] (8,0) -- (9.5,0);

\node[] at (2.65,-2.5) {\large $t_{\textnormal{DoF}}^{2}$};
\draw[HSUred,fill=HSUred](2.65,0.5)circle(3.5pt);
\draw[->,ns1, HSUred, dotted] (2.65,0.5) -- (2.65,-1.8);
\draw[HSUred,fill=HSUred](2.65,-2)circle(3.5pt);

%
\draw[ns2] (1.3,1) -- (4,0);
\draw[ns2] (4,0) -- (4,4);
\draw[ns2] (4,4) -- (1.3,5);
\draw[ns2] (1.3,5) -- (1.3,1);
\draw[ns2] (1.3,3) -- (4,2);
\draw[ns2] (2.65,0.5) -- (2.65,4.5);
\draw[ns2] (1.3,2) -- (2.65,1.5);
\draw[ns2] (1.975,0.75) -- (1.975,2.75);
\draw[ns2] (1.975,2.25) -- (2.65,2.);
\draw[ns2] (2.3125,1.625) -- (2.3125,2.625);
\draw[ns2] (1.975,2.) -- (2.3125,1.875);
\draw[ns2] (2.1438,1.7) -- (2.1438,2.2);
\draw[ns2] (2.65,3.5) -- (4,3.);
\draw[ns2] (3.35,2.25) -- (3.35,4.25);

\node[] at (6.65,-2.5) {\large $t_{\textnormal{DoF}}^{3}$};
\draw[HSUred,fill=HSUred](6.65,0.5)circle(3.5pt);
\draw[->,ns1, HSUred, dotted] (6.65,0.5) -- (6.65,-1.8);
\draw[HSUred,fill=HSUred](6.65,-2)circle(3.5pt);

\draw[ns2] (8,0) -- (8,4);
\draw[ns2] (0,4) -- (0,0);
\draw[ns2] (5.3,1) -- (8,0);
\draw[ns2] (8,0) -- (8,4);
\draw[ns2] (8,4) -- (5.3,5);
\draw[ns2] (5.3,5) -- (5.3,1);
\draw[ns2] (5.3,3) -- (8,2);
\draw[ns2] (6.65,0.5) -- (6.65,4.5);
\draw[ns2] (5.3,2) -- (6.65,1.5);
\draw[ns2] (5.975,0.75) -- (5.975,2.75);
\draw[ns2] (5.975,2.25) -- (6.65,2.);
\draw[ns2] (6.3125,1.625) -- (6.3125,2.625);
\draw[ns2] (5.975,2.) -- (6.3125,1.875);
\draw[ns2] (6.1438,1.7) -- (6.1438,2.2);
\draw[ns2] (6.65,3.5) -- (8,3.);
\draw[ns2] (7.35,2.25) -- (7.35,4.25);

%


\node[] at (9.65,-2.5) {\large $t_{\textnormal{DoF}}^{1}$};
\draw[HSUred,fill=HSUred](9.65,0.5)circle(3.5pt);
\draw[->,ns1, HSUred, dotted] (9.65,0.5) -- (9.65,-1.8);
\draw[HSUred,fill=HSUred](9.65,-2)circle(3.5pt);

%
\draw[ns1] (8.3,1) -- (11,0);
\draw[ns1] (11,0) -- (11,4);
\draw[ns1] (11,4) -- (8.3,5);
\draw[ns1] (8.3,5) -- (8.3,1);
\draw[ns1] (8.3,3) -- (11,2);
\draw[ns1] (9.65,0.5) -- (9.65,4.5);
\draw[ns1] (9.65,1.5) -- (11,1.0);
\draw[ns1] (10.35,0.25) -- (10.35,2.25);
\draw[ns1] (9.65,2.0) -- (10.325,1.75);
\draw[ns1] (10,1.375) -- (10,2.375);
\draw[ns1] (10.325,1.75) -- (11,1.5);
\draw[ns1] (10.675,1.125) -- (10.675,2.125);
\draw[ns1] (10,2.125) -- (10.35,2.0);
\draw[ns1] (10.175,1.8125) -- (10.175,2.3125);
\draw[ns1] (8.3,2) -- (9.65,1.5);
\draw[ns1] (8.975,0.75) -- (8.975,2.75);

\node[] at (19.15,-2.5) {\large $t_{n+1}$};
\draw[HSUred,fill=HSUred](19.15,0.5)circle(3.5pt);
\draw[->,ns1, HSUred, dotted] (19.15,0.5) -- (19.15,-1.8);
\draw[HSUred,fill=HSUred](19.15,-2)circle(3.5pt);

\draw[very thick,HSUred,decorate,decoration={brace,amplitude=14pt}]
(6.8,5) -- (17.8,5) node[midway, above,yshift=14pt,]{\LARGE $Q_{n+1}$};

\draw[ns1] (6.8,5) -- (8.3,5);
\draw[ns1] (9.5,4) -- (11,4);
\draw[ns1,dashed] (6.8,1) -- (8.3,1);
\draw[ns1] (9.5,0) -- (11,0);

\draw[ns1] (11,0) -- (19,0);
\draw[ns1] (19,4) -- (11,4);
\draw[ns1] (8.3,5) -- (16.3,5); 
\draw[ns1, dashed] (8.3,1) -- (16.3,1); 

\draw[ns1] (16.3,5) -- (17.8,5);
\draw[ns1] (19,4) -- (20.5,4);
\draw[ns1,dashed] (16.3,1) -- (17.8,1);
\draw[ns1] (19,0) -- (20.5,0);

\draw[ns1,dashed] (17.8,1) -- (20.5,0);
\draw[ns1] (20.5,0) -- (20.5,4);
\draw[ns1] (20.5,4) -- (17.8,5);
\draw[ns1,dashed] (17.8,5) -- (17.8,1);

\node[] at (13.65,-2.5) {\large $t_{\textnormal{DoF}}^{2}$};
\draw[HSUred,fill=HSUred](13.65,0.5)circle(3.5pt);
\draw[->,ns1, HSUred, dotted] (13.65,0.5) -- (13.65,-1.8);
\draw[HSUred,fill=HSUred](13.65,-2)circle(3.5pt);

%
\draw[ns2] (12.3,1) -- (15,0);
\draw[ns2] (15,0) -- (15,4);
\draw[ns2] (15,4) -- (12.3,5);
\draw[ns2] (12.3,5) -- (12.3,1);
\draw[ns2] (12.3,3) -- (15,2);
\draw[ns2] (13.65,0.5) -- (13.65,4.5);
\draw[ns2] (13.65,1.5) -- (15,1.0);
\draw[ns2] (14.35,0.25) -- (14.35,2.25);
\draw[ns2] (13.65,2.0) -- (14.325,1.75);
\draw[ns2] (14,1.375) -- (14,2.375);
\draw[ns2] (14.325,1.75) -- (15,1.5);
\draw[ns2] (14.675,1.125) -- (14.675,2.125);
\draw[ns2] (14,2.125) -- (14.35,2.0);
\draw[ns2] (14.175,1.8125) -- (14.175,2.3125);
\draw[ns2] (12.3,2) -- (13.65,1.5);
\draw[ns2] (12.975,0.75) -- (12.975,2.75);

\node[] at (17.65,-2.5) {\large $t_{\textnormal{DoF}}^{3}$};
\draw[HSUred,fill=HSUred](17.65,0.5)circle(3.5pt);
\draw[->,ns1, HSUred, dotted] (17.65,0.5) -- (17.65,-1.8);
\draw[HSUred,fill=HSUred](17.65,-2)circle(3.5pt);

%
\draw[ns2] (16.3,1) -- (19,0);
\draw[ns2] (19,0) -- (19,4);
\draw[ns2] (19,4) -- (16.3,5);
\draw[ns2] (16.3,5) -- (16.3,1);
\draw[ns2] (16.3,3) -- (19,2);
\draw[ns2] (17.65,0.5) -- (17.65,4.5);
\draw[ns2] (17.65,1.5) -- (19,1.0);
\draw[ns2] (18.35,0.25) -- (18.35,2.25);
\draw[ns2] (17.65,2.0) -- (18.325,1.75);
\draw[ns2] (18,1.375) -- (18,2.375);
\draw[ns2] (18.325,1.75) -- (19,1.5);
\draw[ns2] (18.675,1.125) -- (18.675,2.125);
\draw[ns2] (18,2.125) -- (18.35,2.0);
\draw[ns2] (18.175,1.8125) -- (18.175,2.3125);
\draw[ns2] (16.3,2) -- (17.65,1.5);
\draw[ns2] (16.975,0.75) -- (16.975,2.75);


\node[] at (21,-2) {\large $t$};
\draw[->,ns1] (-4.2,-2) -- (20.5,-2);

\end{tikzpicture}}

\caption{Two consecutive space-time slabs, exemplary for a discontinuous Galerkin
$\textnormal{dG}(2)$ time discretization generated with three Gaussian quadrature
points. The three degrees of freedom (DoF) time points on each slab
are the support points for the temporal basis functions.
Each of the illustrated slabs here has one temporal cell and
an independent and adaptively refined spatial triangulation.}
\label{fig:2:STCylinderSlabs}
\end{figure}

In this section we analyse the implementation of space-time tensor-product
spaces in detail. An exemplary illustration of a space-time cylinder that is
distributed into space-time tensor-product slabs is given in
Fig. \ref{fig:2:STCylinderSlabs}.
Precisely, we explain the details here for the scalar-valued transport
equation with primal and dual finite element spaces.
The implementation for the primal vector-valued Stokes flow problem is very
similar with the difference that the spatial finite element has
$d$+$1$ components for the velocity and pressure variables.
We denote the number of spatial degrees of freedom by
$N_{\textnormal{DoF}}^{\textnormal{s,n}}$ for one degree of freedom in time and
the number of temporal degrees of freedom by $N_{\textnormal{DoF}}^{\textnormal{t,n}}$
on the $n$-th slab.

To implement the space-time tensor-product space,
as illustrated in Fig. \ref{fig:2:STCylinderSlabs},
we start with the usual discretization for the finite element method in space
having only one degree of freedom in time in an adaptive time marching process,
but here we do this for each slab.
Therefore, we generate the geometrical triangulation, i.e. a spatial mesh,
and colourize the boundaries.
Boundary colours can mark for instance Dirichlet type boundary conditions,
Neumann type boundary conditions, etc.
Next, we initialize each slab by creating an independent copy of the generated
spatial triangulation.

\begin{figure}[tbh]
\centering

\resizebox{0.85\linewidth}{!}{%
\begin{tikzpicture}

\draw[thick,rounded corners,fill=lightgray,opacity=0.5] (-0.5,-0.3) rectangle (14.8,1.8);
\draw[thick,rounded corners,fill=lightgray] (-0.2,0) rectangle (4,1.5);
\draw (1.9,1) node[]{\large Geometrical};
\draw (1.9,0.4) node[]{\large Mesh Description};
\draw[thick,rounded corners,fill=lightgray] (4.5,0) rectangle (7.5,1.5);
\draw (6,1) node[]{\large Mapping};
\draw (6,0.4) node[]{\large (Mesh Cell)};
\draw[thick,rounded corners,fill=lightgray] (8,0) rectangle (11,1.5);
\draw (9.5,1) node[]{\large Boundary};
\draw (9.5,0.4) node[]{\large Manifold};
\draw[thick,rounded corners,fill=lightgray] (11.5,0) rectangle (14.5,1.5);
\draw (13,1) node[]{\large Boundary};
\draw (13,0.4) node[]{\large Colorisation};

\draw[thick,rounded corners,fill=Eggshell] (-0.2,-2.6) rectangle (4.0,-1.1);
\draw (1.90,-1.6) node[]{\large Spatial};
\draw (1.90,-2.2) node[]{\large Triangulation $\mathcal{T}_{h,n}$};

\draw[ultra thick, stealth reversed-stealth](1.9,-0.07) -- (1.9,-1.40);

\draw[thick,rounded corners,fill=white,opacity=1., draw=gray] (-0.2,-8.5) rectangle (4.6,-4.6);
\draw[gray] (2.2,-8.1) node[]{\large Dual FE Space};

\draw[ultra thick, stealth reversed-stealth, draw=gray](2.8,-2.7) -- (2.8,-4.9);

\draw[thick,rounded corners,fill=white,opacity=1.] (-0.5,-7.6) rectangle (4.3,-3.7);
\draw (1.9,-7.2) node[]{\large Primal FE Space};

\draw[ultra thick, stealth reversed-stealth](1.9,-2.7) -- (1.9,-3.95);

\draw[thick,rounded corners] (-0.2,-4.7) rectangle (4.0,-4.0);
\draw (1.90,-4.35) node[]{\large Constraints};

\draw[thick,rounded corners] (-0.2,-5.7) rectangle (4.0,-5.0);
\draw (1.90,-5.37) node[]{\large Fin. Elem. System};

\draw[thick,rounded corners,fill=Eggshell] (-0.2,-6.7) rectangle (4.0,-6.0);
\draw (1.90,-6.38) node[]{\large Degrees of Freedom};


\draw[thick,rounded corners,fill=lightgray] (10.3,-2.6) rectangle (14.5,-1.1);
\draw (12.4,-1.6) node[]{\large Temporal};
\draw (12.4,-2.2) node[]{\large Endpoints $(t_{n-1}, t_n)$};

\draw[thick,rounded corners,fill=Eggshell] (1.+4.5,-2.6) rectangle (1.+8.7,-1.1);
\draw (1.+6.60,-1.6) node[]{\large Temporal};
\draw (1.+6.60,-2.2) node[]{\large Triangulation $\mathcal{T}_{\tau,n}$};

\draw[ultra thick, stealth reversed-stealth](10.6,-1.9) -- (9.4,-1.9);

\draw[thick,rounded corners,fill=white,opacity=1., draw=gray]
	(5.6-0.2,-8.5) rectangle (5.6+4.6,-4.6);
\draw[gray] (5.6+2.2,-8.1) node[]{\large Dual FE Space};

\draw[ultra thick, stealth reversed-stealth, draw=gray]
	(5.6+2.8,-2.7) -- (5.6+2.8,-4.9);

\draw[thick,rounded corners,fill=white,opacity=1.]
	(5.6-0.5,-7.6) rectangle (5.6+4.3,-3.7);
\draw (5.6+1.9,-7.2) node[]{\large Primal FE Space};

\draw[ultra thick, stealth reversed-stealth]
	(5.6+1.9,-2.7) -- (5.6+1.9,-3.95);


\draw[thick,rounded corners] (5.6-0.2,-5.7) rectangle (5.6+4.0,-5.0);
\draw (5.6+1.90,-5.37) node[]{\large Fin. Elem. System};

\draw[thick,rounded corners,fill=Eggshell] (5.6+-0.2,-6.7) rectangle (5.6+4.0,-6.0);
\draw (5.6+1.90,-6.38) node[]{\large Degrees of Freedom};

\draw[thick,rounded corners,fill=SlabGreen] (-0.5,-11.7) rectangle (14.8,-9.6); 

\draw[ultra thick, stealth reversed-, gray](4.0,-8.2) -- (4.0,-9.3);
\draw[ultra thick, -stealth, gray](0.7+4.0,-9.3) -- (0.7+4.0,-10.0);
\draw[ultra thick, stealth reversed-, gray] (5.9,-8.2) -- (5.9,-9.3);
\draw[ultra thick, gray] (5.9,-9.3) -- (4.0,-9.3);

\draw[ultra thick, stealth reversed-](3.7,-7.3) -- (3.7,-9.0);
\draw[ultra thick, -stealth](0.7+3.7,-9.0) -- (0.7+3.7,-10.0);
\draw[ultra thick, stealth reversed-] (5.6,-7.3) -- (5.6,-9.0);
\draw[ultra thick] (5.6,-9.0) -- (3.7,-9.0);

\draw[thick,rounded corners,fill=LightOlive,opacity=0.8] (-0.2,-3.2-8.2) rectangle (4.0,-3.2-6.7);
\draw (1.9,-3.2-7.20) node[]{\large Space-Time};
\draw (1.9,-3.2-7.80) node[]{\large Degrees of Freedom};

\draw[thick,rounded corners,fill=LightOlive,opacity=0.8]
	(.5+4.5,-3.2-8.2) rectangle (.5+8.7,-3.2-6.7);
\draw (.5+6.60,-3.2-7.20) node[]{\large Space-Time};
\draw (.5+6.60,-3.2-7.80) node[]{\large Constraints};

\draw[thick,rounded corners,fill=LightOlive,opacity=0.8] (10.3,-3.2-8.2) rectangle (14.4,-3.2-6.7);
\draw (12.4,-3.2-7.20) node[]{\large Space-Time};
\draw (12.4,-3.2-7.80) node[]{\large Sparsity Patterns};

\end{tikzpicture}
}

\caption{Illustration of the generation of space-time elements for
tensor-product spaces on a slab.
}
\label{fig:3:DTMslab}
\end{figure}
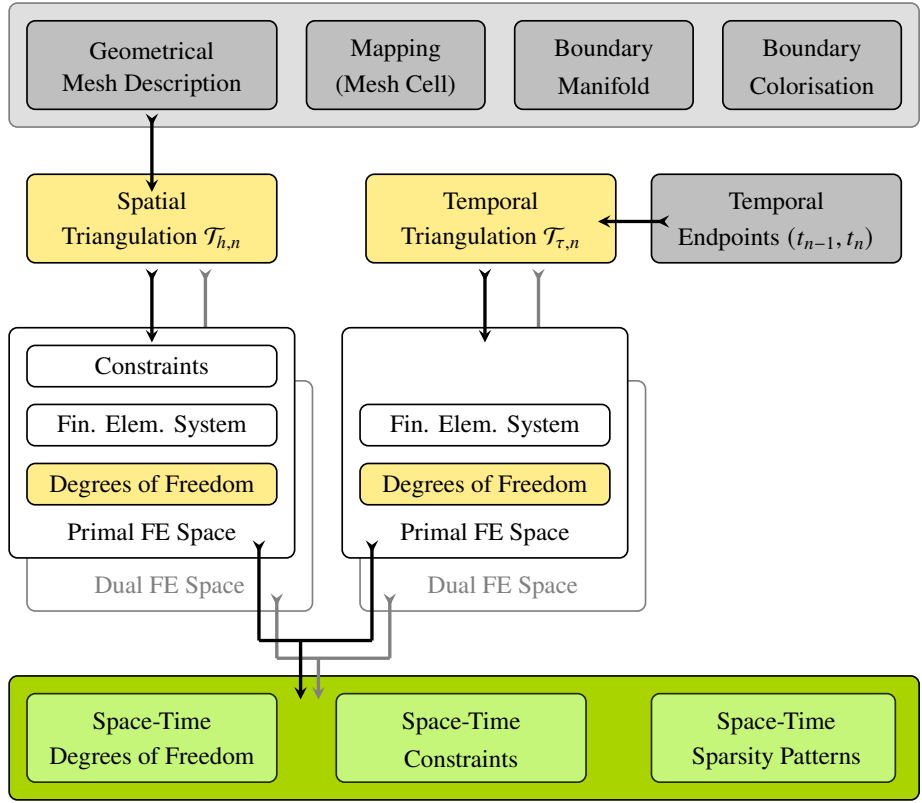

Then, for one degree of freedom in time on each slab,
we distribute the spatial degrees of freedom and generate affine constraints
objects.
Remark that an affine constraints object may include information on
handling degrees of freedom on hanging nodes or on Dirichlet type boundary nodes.
The sparsity pattern for a sparse matrix is now generated with
the geometric triangulation, the spatial degree of freedom (DoF) handler and
the constraints object for one degree of freedom in time.

Next, the space-time tensor-product degrees of freedom on a slab are aligned
by their local degree of freedom in time on a slab. Precisely, the first
degree of freedom in time has the global number $0$ and the last one has the
number $N_{\textnormal{DoF}}^{\textnormal{t,n}}$-$1$.
The numbering of the local temporal degrees of freedom is increasingly ordered
by their temporal mesh cell index.
Remark that we have an one-dimensional additional triangulation (temporal mesh)
for the time subinterval $(t_{n-1}, t_n)$ corresponding to the $n$-th slab;
refer to the Fig. \ref{fig:3:DTMslab} for details.
Overall, we have $N_{\textnormal{DoF}}^{\textnormal{t,n}}$ times
$N_{\textnormal{DoF}}^{\textnormal{s,n}}$ degrees of freedom on the $n$-th slab.

Next, the space-time tensor-product constraints are created by taking the original
constraints object and shifting all entries accordingly such that the
$N_{\textnormal{DoF}}^{\textnormal{t,n}}$ are represented.
Precisely, the spatial degrees of freedom from
$0$ to $N_{\textnormal{DoF}}^{\textnormal{s}}$-$1$ are associated to the first
local temporal degree of freedom on a slab.
If there are more than one temporal degrees of freedom on a slab,
the corresponding spatial degrees of freedom are shifted by the number
$N_{\textnormal{DoF}}^{\textnormal{s,n}}$ times the local temporal degree of
freedom index.

\begin{figure}[thbp]
\centering


\resizebox{.4\linewidth}{!}{%
\input{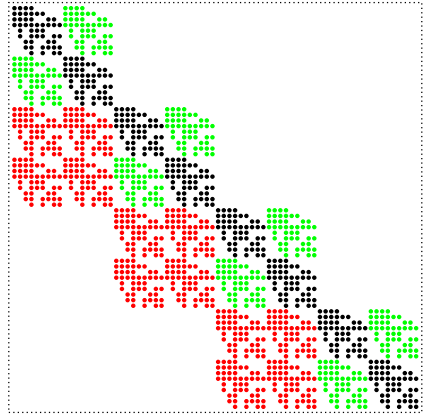}
}

\caption{
Sparsity pattern for the primal operator on a slab for cG(1)-dG(1)
with 4 mesh cells in space and 4 mesh cells in time
before condensing the Dirichlet nodes.
The blocks having black dots correspond to a classical sparsity pattern
for one degree of freedom in time.
The blocks with green dots are additional couplings between the two time basis
functions on a temporal cell.
The blocks with red dots are additional couplings from the temporal jump trace
operator between two time cells.}

\label{fig:4:sparsitypattern:cG1dG1}
\end{figure}

\begin{figure}[tbhp]
\centering

\includegraphics[width=.35\linewidth]{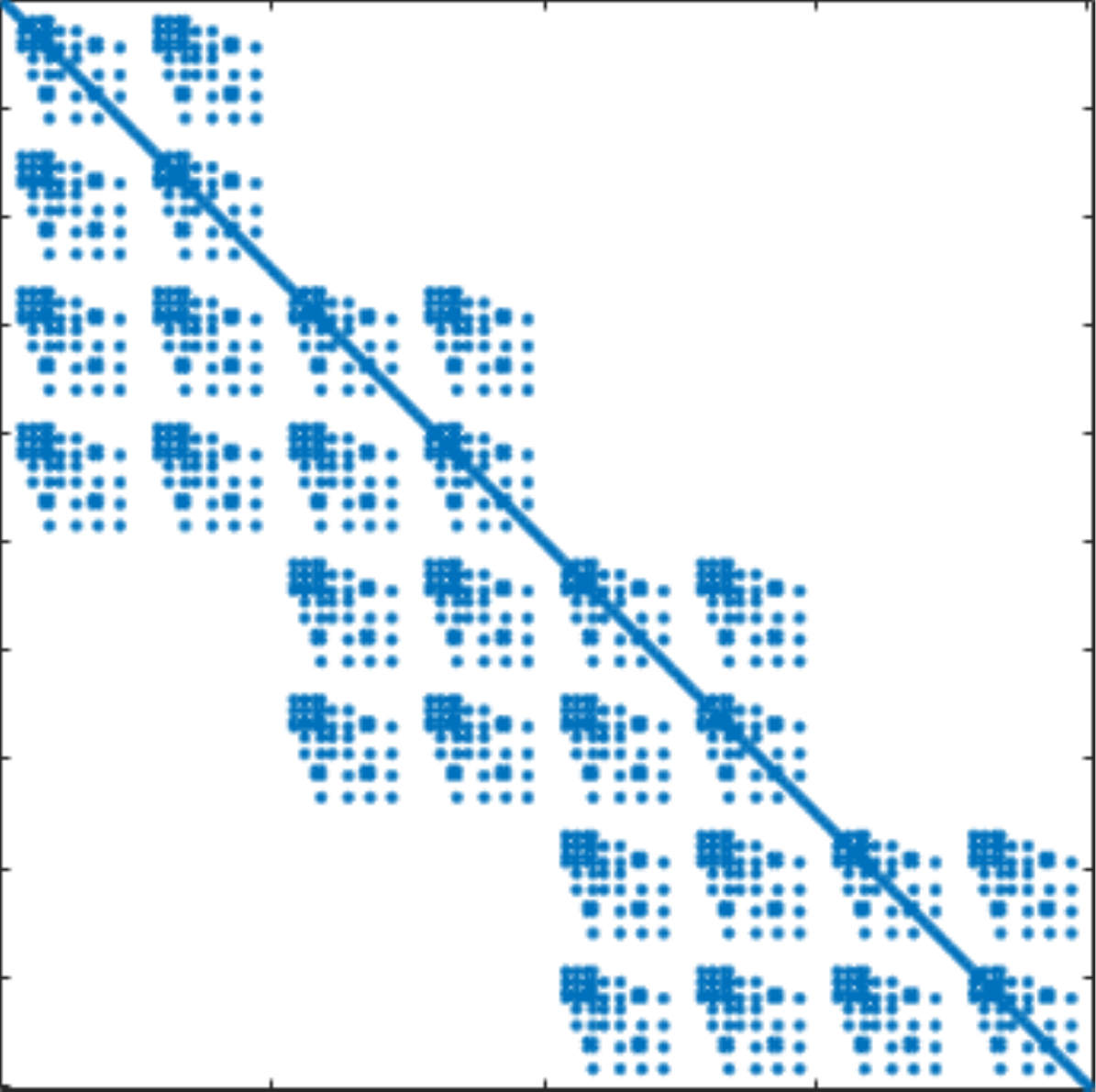}
~
\includegraphics[width=.35\linewidth]{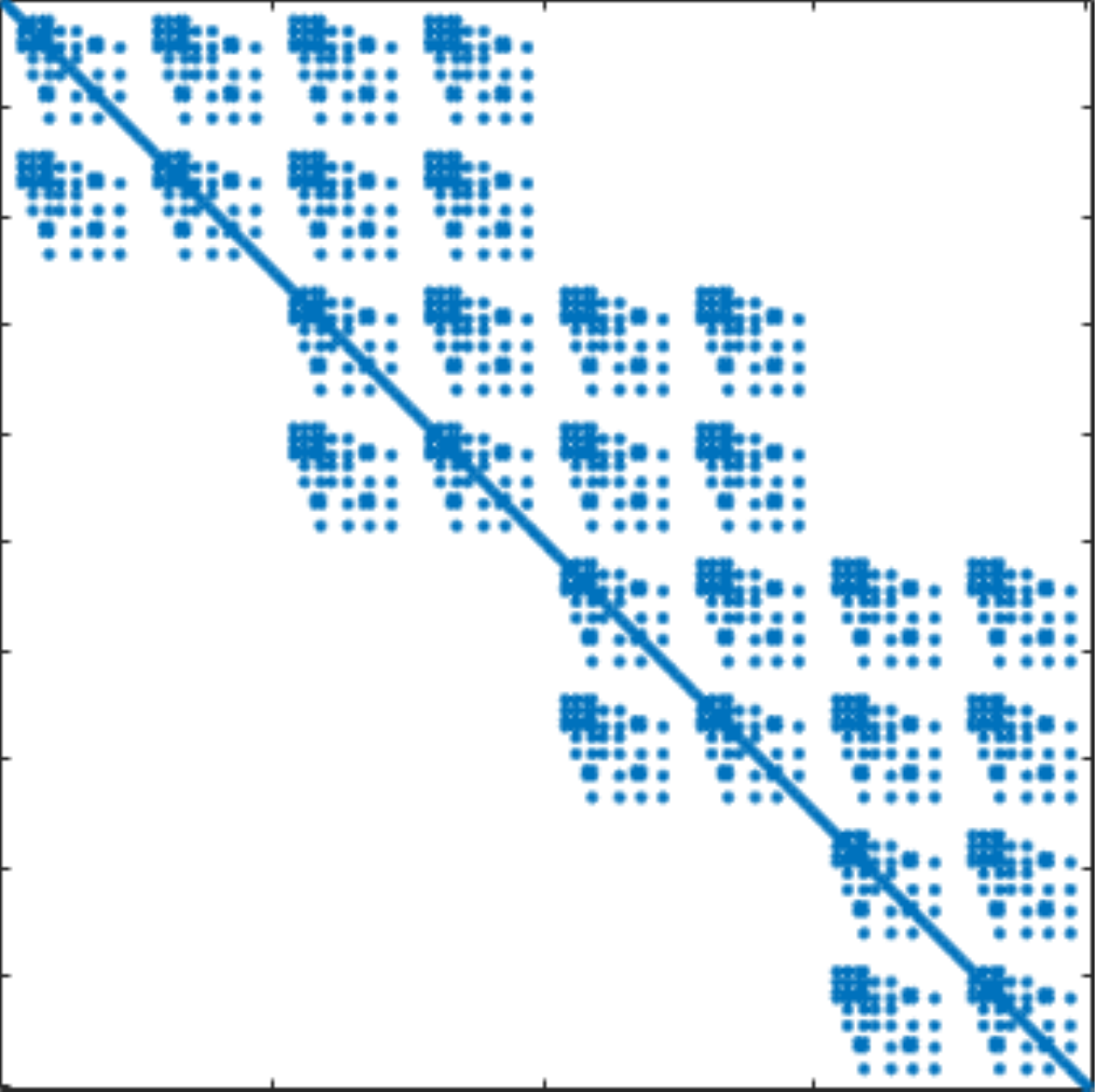}

\caption{%
Primal and dual operator matrix sparsity patterns on a slab
for cG(1)-dG(1) on 16 cells in space and 4 cells in time
after condensing the Dirichlet nodes.
}

\label{fig:5:sparsitypatterns}
\end{figure}

\begin{figure}[tbhp]
\centering

\includegraphics[width=.35\linewidth]{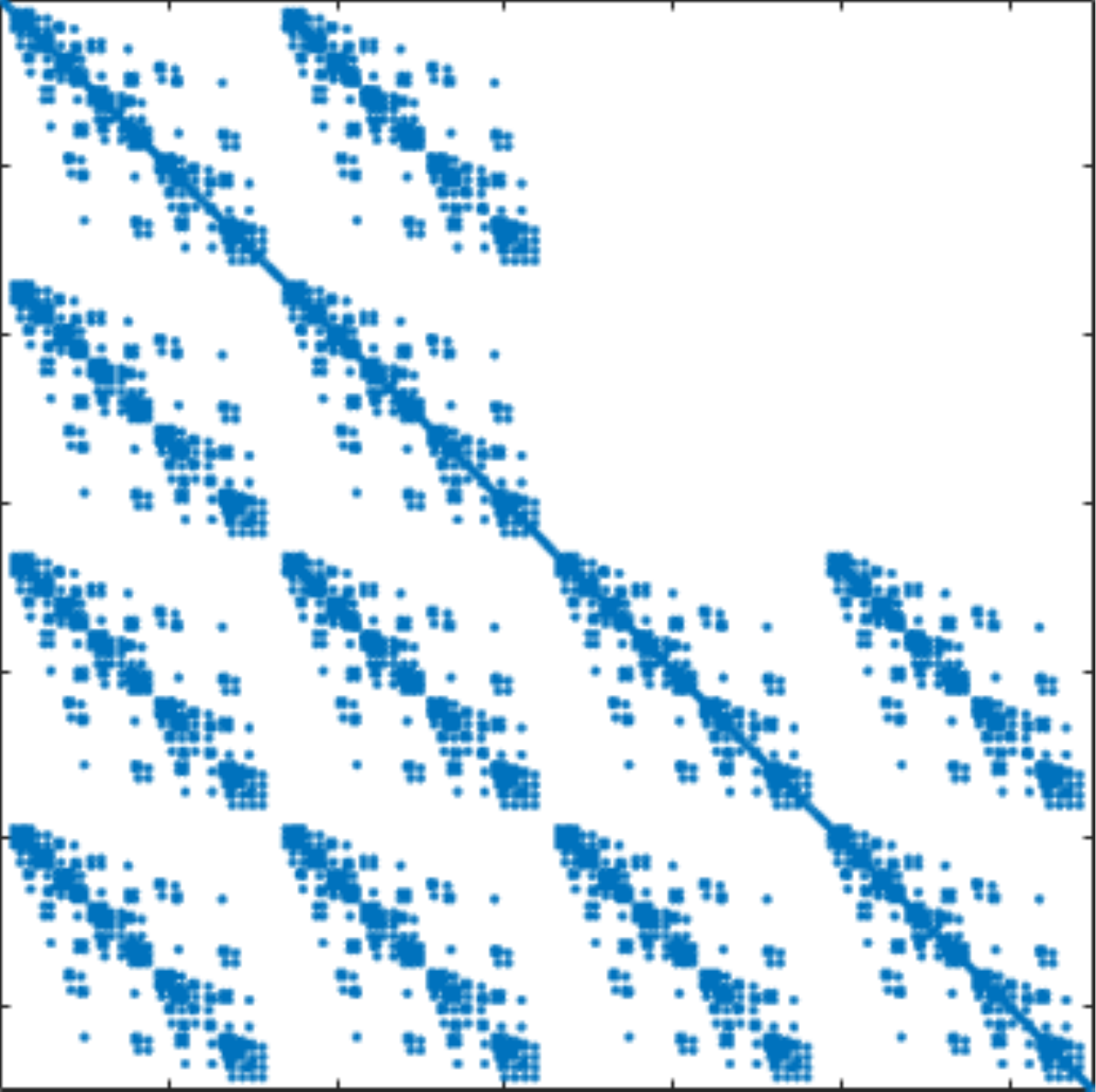}
~
\includegraphics[width=.35\linewidth]{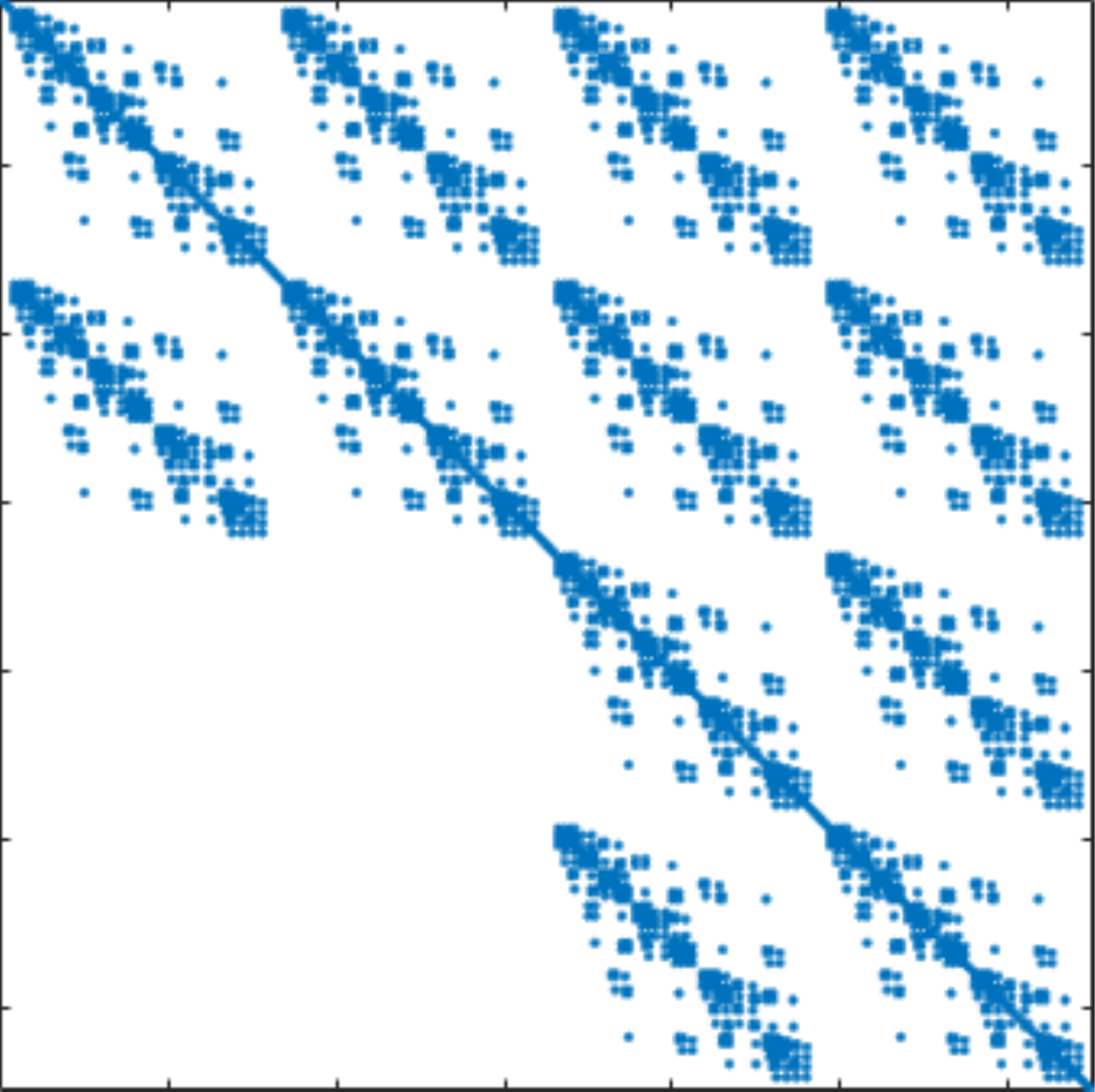}

\caption{%
Primal and dual operator matrix sparsity patterns on a slab
for cG(1)-dG(1) on 64 cells in space and 2 cells in time
after condensing the Dirichlet nodes.
}

\label{fig:6:sparsitypatterns}
\end{figure}

For each degree of freedom in time, the sparsity pattern is now copied into the
diagonal blocks for the space-time tensor product sparsity pattern.
A higher-order polynomial degree in time introduces couplings between the temporal
basis functions resulting in additional coupling blocks.
For the case of more than one time cell per slab, additional couplings appear for
temporal derivatives between the time basis functions of two consecutive
time cells. For the primal problem, the evolution is forward in time and therefore
these couplings appear in the left lower part.
For the dual problem, the evolution is backward in time and therefore the
coupling diagonals appear in the right upper part.
Exemplary sparsity patterns are given in Fig. \ref{fig:4:sparsitypattern:cG1dG1},
Fig. \ref{fig:5:sparsitypatterns} and Fig. \ref{fig:6:sparsitypatterns}.

For the assembly process we can use the basis functions and their derivatives in
time similar to the classical finite element approach in space.
But the distribution of the local contributions must respect the order of the 
temporal basis functions.
First, we take the mapping from a local to a global degree of freedom in space.
To respect the temporal basis functions, we shift the local to global mapping
accordingly by the factor of local degrees of freedom in space on a spatial cell.
This results in a shift of each global degree of freedom by the factor of
$N_{\textnormal{DoF}}^{\textnormal{s,n}}$ times the global degree of freedom of the
respective basis function in time.
The local matrix has therefore the size of the local degrees of freedom on a
spatial mesh cell times the local degrees of freedom in time on a temporal mesh cell.

In the case of more than on time cell per slab, an additional local matrix is
assembled for the coupling of the trial basis functions of the previous time cell
and the test basis functions of the current time cell.
This implements the negative part of the jump trace operator in time
which is transferred to the right-hand side in a classical time marching approach.

Finally, the space-time constraints of the slab have to be applied to
the system matrix, the solution vector and the right-hand side vector.
The space-time hanging node constraints have to be condensed in the solution vector
after solving the linear system for all degrees of freedom on the slab.

\clearpage
\section{Algorithm}
\label{sec:5:algorithm}

Here we present the multirate in time adaptive algorithm,
give the definition of the (localized) error indicators and explain the 
approximation techniques used for the weights occurring within these indicators.
Our space-time adaptivity strategy uses the following
algorithm.

\noindent\rule{\textwidth}{1pt}
  \begin{center}
   \textbf{Algorithm: goal-oriented multirate space-time adaptivity}
  \end{center}
\vspace{-0.3cm}
\noindent\rule{\textwidth}{0.5pt}
\textbf{Initialization:}
Generate the initial space-time slabs $Q_{n}^{1}=\mathcal{T}_{h,n}^{1}\times 
\mathcal{T}_{\tau,n}^{1}\,,n=1,\dots,N^1\,,$  as well as
$Q_{n}^{\textnormal{F},1}=\mathcal{T}_{h,n}^{\textnormal{F},1}\times 
\mathcal{T}_{\sigma,n}^{\textnormal{F},1}\,,n=1,\dots,N^{\textnormal{F},1}\,,$ 
$N^{\textnormal{F},1} \leq N^1\,,$ for the transport and Stokes flow problem, 
respectively, where we restrict $\mathcal{T}_{\sigma,n}^{1},
\mathcal{T}_{\sigma,n}^{\textnormal{F},1}$ to consist of only one cell in time
for each slab.

\noindent\rule{\textwidth}{0.5pt}
\textbf{DWR-loop $\ell=1,\dots$:}
\begin{enumerate}
\item %
  \textbf{Find the solutions} $\{\convection_{\sigma h},\pressure_{\sigma h}\} \in 
  Y_{\sigma h^{\textnormal{F}}}^{\textnormal{dG}(0),p}$
  of the Stokes flow problem (\ref{eq:18:B_v_p_psi_chi_eq_F_psi}).

\item \textbf{Find the primal solution}
  $\concentration_{\tau h} \in X_{\tau h}^{\textnormal{dG}(r),p}$
  of the stabilized transport problem (\ref{eq:15:A_S_G_tauh}).

\item \textbf{Break if the goal yields convergence}.

\item \textbf{Find the dual solution}
  $\dualz_{\tau h} \in X_{\tau h}^{\textnormal{dG}(r),q}, q > p,$
  of the dual transport problem (\ref{eq:22c:DualProblems_fullydiscrete}).

\item \textbf{Evaluate the localized a posteriori space-time error indicators}
  $\eta_{h}$ and $\eta_{\tau}$ given by Eq. \eqref{eq:26:eta_h} and
   \eqref{eq:25:eta_tau}, respectively, for the transport problem.

\item Refine the temporal and spatial meshes of the \textbf{transport} problem 
as follows:

\begin{enumerate}
\item[(i)] \textbf{If} 
$|\eta_{\tau}^{\ell}|>
\omega\,|\eta_{h}^{\ell}|\,,\omega \geq 1$: 
  
\textbf{Mark the slabs} $Q_{\tilde{n}}^{\ell}$, $\tilde{n}\in\{1,\dots,N^\ell\}$,
\textbf{for temporal refinement} if the corresponding 
$\eta_{\tau}^{\tilde{n},\ell}$ is in 
the set of $\theta_\tau^\textnormal{top}\,,0\leq\theta_\tau^\textnormal{top}\leq 1\,,$ 
percent of the worst indicators.
  
\item[(ii)] \textbf{Else if} 
$|\eta_{h}^{\ell}|>
\omega\,|\eta_{\tau}^{\ell}|$: 
  
\textbf{Mark the cells\,} $\tilde{K}\in\mathcal{T}_{h,n}^{\ell}$
\textbf{\,for spatial refinement} if the corresponding 
$\eta_{h}^{n,\ell}|_{\tilde{K}}$ 
is in the set of $\theta_{h,1}^\textnormal{top}$ or 
$\theta_{h,2}^\textnormal{top}$
(for a slab that is or is not marked for temporal refinement),
$0\leq\theta_{h,2}^\textnormal{top}
\leq
\theta_{h,1}^\textnormal{top}\leq 1\,,$ 
percent of the worst indicators,
\textbf{or}, respectively,
mark \textbf{for spatial coarsening} if 
$\eta_{h}^{n,\ell}|_{\tilde{K}}$  is in the 
set of $\theta_h^\textnormal{bottom}\,,
0\leq\theta_h^\textnormal{bottom}\leq 1\,,$ percent of the best 
indicators.

\item[(iii)] \textbf{Else}: 

\textbf{Mark the slabs} $Q_{\tilde{n}}^{\ell}$
\textbf{for temporal refinement} as well as \textbf{mark the cells\,} 
$\tilde{K}\in\mathcal{T}_{h,n}^{\ell}$ \textbf{\,for spatial coarsening and refinement}
as described in Step~6(i) and Step~6(ii), respectively.
  
\item[(iv)] \textbf{Execute spatial adaptations} on all slabs
of the transport problem under the use of mesh smoothing operators.
  
\item[(v)] \textbf{Execute temporal refinement} on all slabs of the transport 
problem. 
\end{enumerate}

\item \textbf{If} $\|\convection-\convection_{\sigma h}\|_{(0,T)\times\Omega} >
\|u-u_{\tau h}\|_{(0,T)\times\Omega}$:

Refine the spatial and temporal mesh of the \textbf{Stokes flow} problem 
globally. 

\item Increase $\ell$ to $\ell+1$ and return to Step~1.
\end{enumerate}
\vspace{-0.3cm}
\noindent\rule{\textwidth}{0.5pt}

Regarding this algorithm, we note the following issues.
\begin{remark}
\label{rem:?:algorithm}
~\\
\vspace{-0.6cm}
\begin{itemize}
\item For the spatial discretization of the Stokes flow problem we are using 
Taylor-Hood elements $Q_p/Q_{p-1}\,, p \geq 2$.

\item Within the Steps 2, 4 and 5 of the algorithm, the computed convection field
$\convection_{\sigma h}$ of the Stokes problem is interpolated to the adaptively 
refined spatial and temporal triangulation of the space-time slabs.

\item Our simulation tools of the \texttt{DTM++} project are frontend solvers
for the \texttt{deal.II} library; cf. \cite{dealiiReference93}.

\item Technical details of the implementation are given in
\cite{Koecher2019,Bause2021}.
\end{itemize}
\end{remark}

In the following, we give some details regarding the localization of the error 
representations that are derived in Thm.~\ref{thm:1:ErrorRepresentation}. 
Their practical realization and the definition of error indicators $\eta_{\tau}$ 
and $\eta_{h}$ is obtained by neglecting the remainder terms $\mathcal{R}_{\tau}$ 
and $\mathcal{R}_{h}$ of the result given in Thm.~\ref{thm:1:ErrorRepresentation}
and splitting the resulting quantities into elementwise contributions.

\begin{align}
\label{eq:25:eta_tau}
\nonumber
J(\concentration)-J(\concentration_{\tau}) & \doteq
\frac{1}{2}\rho_{\mathrm{t}}^{n,\ell}(\concentration_{\tau})(\dualz-\tilde{\dualz}_{\tau})
+ \frac{1}{2}\rho_{\mathrm{t}}^{\ast,n}(\concentration_{\tau},\dualz_{\tau})
(\concentration-\tilde{\concentration}_{\tau})
\\
\nonumber
&
+ \frac{1}{2} \mathcal{D}_{\tau}^{\prime,n,\ell}(\concentration_{\tau},\dualz_{\tau})
(\tilde{\concentration}_{\tau}-\concentration_{\tau},\tilde{\dualz}_{\tau}-\dualz_{\tau})
+ \mathcal{D}_{\tau h}^{n,\ell}(\concentration_{\tau h},\dualz_{\tau h})
\\
& =: \eta_{\tau}^\ell = \displaystyle \sum_{n=1}^{N^\ell} 
\eta_{\tau}^{n,\ell}\,,
\end{align}

\begin{align}
\label{eq:26:eta_h}
\nonumber
J(\concentration_{\tau})-J(\concentration_{\tau h}) & \doteq
\frac{1}{2}\rho_{\mathrm{t}}^{n,\ell}(\concentration_{\tau h})(\dualz_{\tau}-\tilde{\dualz}_{\tau h})
+ \frac{1}{2}
\rho_{\mathrm{t}}^{\ast,n,\ell}(\concentration_{\tau h},\dualz_{\tau h})
(\concentration_{\tau}-\tilde{\concentration}_{\tau h})
\\
\nonumber
&
+ \frac{1}{2} \mathcal{D}_{\tau h}^{\prime,n,\ell}(\concentration_{\tau h},\dualz_{\tau h})
(\tilde{\concentration}_{\tau h}-\concentration_{\tau h},\tilde{\dualz}_{\tau h}-\dualz_{\tau h})
+ \mathcal{D}_{\tau h}^{n,\ell}(\concentration_{\tau h},\dualz_{\tau h})
\\
& =: \eta_h^\ell = \displaystyle \sum_{n=1}^{N^\ell} \eta_{h}^{n,\ell}
= \displaystyle \sum_{n=1}^{N^\ell}
 \sum\limits_{K\in\mathcal{T}_h^{n,\ell}} \eta_{h,K}^{n,\ell}\,.
\end{align}

To compute the error indicators $\eta_{\tau}$ and $\eta_{h}$ we replace
all unknown solutions by the approximated fully discrete solutions
$\concentration_{\tau h} \in X_{\tau h}^{\textnormal{dG}(r), p}$,
$\dualz_{\tau h} \in X_{\tau h}^{\textnormal{dG}(r), q}$,
with $p < q$, and $\convection_{\sigma h} \in 
Y_{\sigma h^{\textnormal{F}}}^{\textnormal{dG}(0),p_\convection}$, 
$p_\convection \geq 2$,
whereby the arising weights are approximated in the following way.

\begin{itemize}
\itemsep1.5ex

\item Approximate the temporal weights $u-\tilde{u}_{\tau}$ and $z-\tilde{z}_{\tau}$,
respectively, by means of a higher-order extrapolation using Gauss-Lobatto 
quadrature points given by
\begin{displaymath}
\begin{array}{r@{\,}c@{\,}l}
u-\tilde{u}_{\tau} & \approx & \operatorname{E}_{\tau}^{(r+1)}u_{\tau h}-u_{\tau h}\,,
\\[1.0ex]
z-\tilde{z}_{\tau} & \approx & \operatorname{E}_{\tau}^{(r+1)}z_{\tau h}-z_{\tau h}\,,
\end{array}
\end{displaymath}
using an extrapolation in time operator $\operatorname{E}_{\tau}^{(r+1)}$ 
thats acts on a time cell of length $\tau_K$ and lifts the solution to a piecewise 
polynomial of degree ($r$+$1$) in time.
This approximation technique is a new approach compared to out previous work 
\cite[Sec. 4]{Bause2021}, where a higher-order finite element approximation was 
used, and is done for the purpose to reduce numerical costs solving the dual 
problem. 
\item Approximate the spatial weights $u_{\tau}-\tilde{u}_{\tau h}$ and 
$z_{\tau}-\tilde{z}_{\tau h}$ by means of a patch-wise higher-order interpolation
and a higher-order finite elements approach, respectively, given by
\begin{displaymath}
\begin{array}{r@{\,}c@{\,}l}
u_{\tau}-\tilde{u}_{\tau h} & \approx & \operatorname{I}_{2h}^{(2p)}u_{\tau h}-u_{\tau h}\,,
\\[1.0ex]
z_{\tau}-\tilde{z}_{\tau h} & \approx & z_{\tau h}-\operatorname{R}_{h}^{p}z_{\tau h}\,,
\end{array}
\end{displaymath}
using an interpolation in space operator $\operatorname{I}_{2h}^{(2p)}$ and an
restriction in space operator $\operatorname{R}_{h}^{p}$ that are described in 
detail in our work \cite[Sec. 4]{Bause2021}.

\end{itemize}


\section{Numerical Examples}
\label{sec:6:examples}

In the following section we study the convergence, computational efficiency
and stability of the introduced goal-oriented DWR based adaptivity approach for
the coupled transport and Stokes flow problem.
The first example given in Sec.~\ref{sec:6:1} is an academic test problem with 
given analytical solutions to study the convergence behavior of the two 
subproblems and, in particular, the coupling between them.  
The second example given in Sec.~\ref{sec:6:2} serves to demonstrate the 
performance properties of the algorithm with regard to adaptive mesh refinement 
in space and time.
Finally, the third example in Sec.~\ref{sec:6:3} is motivated by problem of physical
relevance in which we simulate a convection-dominated transport with goal-oriented
adaptivity of a species through a channel with a constraint.

\subsection{Example 1 (Higher-order space-time convergence studies)}
\label{sec:6:1}

In a first numerical example we study the space-time higher-order convergence 
behavior to validate the correctness of the higher-order implementations in space 
and time.
Therefore, we consider the two cases of a solely solved Stokes flow problem 
as well as a non-stabilized solved convection-diffusion-reaction transport problem 
coupled with this Stokes equation via the convection field solution 
$\convection_{\sigma h}$.
The latter may be compared to the results of a solely solved transport equation 
with a constant convection field $\convection=(2,3)^\top$ published in our 
work \cite[Example 1]{Bause2021}.
For this purpose, we investigate problem \eqref{eq:3:stokes_problem} with the given
analytical solution
\begin{equation}
\label{eq:27:BRdynamic}
\begin{array}{r@{\,}c@{\,}l}
\convection(\boldsymbol x, t) &:=&
\begin{pmatrix}
\sin(t)\sin^2(\pi x_1)\sin(\pi x_2)\cos(\pi x_2) \\
-\sin(t)\sin(\pi x_1)\cos(\pi x_1)\sin^2(\pi x_2)
\end{pmatrix}\,,
\\[2.5ex]
\pressure(\boldsymbol x, t) &:=&
\sin(t)\sin(\pi x_1)\cos(\pi x_1)\sin(\pi x_2)\cos(\pi x_2)\,,
\end{array}
\end{equation}
with $\boldsymbol x = (x_1, x_2)^\top \in \mathbb{R}^2\,, t \in \mathbb{R}$ and
$\nabla \cdot \convection = 0$.
The viscosity is set to $\viscosity=0.5$.
The problem is defined on $Q=\Omega\times I:=(0,1)^2\times (0,1]$. The initial 
and boundary conditions are given as
\begin{displaymath}
\convection = 0 \,\text{ on }\, \Sigma_0 = \Omega \times \{ 0 \}\,,
\quad 
\convection = 0 \,\text{ on }\, \Sigma_D = \partial\Omega \times (0,1)\,,
\end{displaymath}
and the volume force term $\stokesforce$ 
is calculated from the given analytical solution \eqref{eq:27:BRdynamic} and 
Eq.~\eqref{eq:3:stokes_problem}. 
This example is a typical test problem for time-dependent incompressible 
flow and can be found, for instance, in \cite[Example 1]{Besier2012}. 

For the following test settings, the solution $\{\convection,\pressure\}$ is 
approximated with the space-time higher-order methods \{cG(2)-dG(2),cG(1)-dG(2)\} 
and \{cG(3)-dG(3),cG(2)-dG(3)\}, respectively.
Due to the same polynomial orders of the spatial and temporal discretizations
with respect to the convection field $\convection$, and lower polynomial order in
space compared to in time with respect to the pressure variable $\pressure$,
we expect experimental orders of convergence 
(EOC $:= -\log_2(|| e ||_\ell / || e ||_{\ell-1})$) 
for the convection field $\convection$ of
$\textnormal{EOC}^{2,2} \approx 3$ for the cG(2)-dG(2) method and
$\textnormal{EOC}^{3,3} \approx 4$ for the cG(3)-dG(3) method, as well as  
experimental orders of convergence for the pressure variable $\pressure$ of
$\textnormal{EOC}^{1,2} \approx 2$ for the cG(1)-dG(2) method and
$\textnormal{EOC}^{2,3} \approx 3$ for the cG(2)-dG(3) method
for a global refinement convergence test.
The results are given by Tab.~\ref{tab:1:BRdynamic:global} and nicely confirm
our expected results for the respective spatial and temporal discretizations.
%
\begin{table}[ht]
\centering

\resizebox{0.75\linewidth}{!}{%
\begin{tabular}{c rrr | cc | cc}
\hline
\hline
$\ell$ & $N$ & $N_K$ & $N_{\text{DoF}}^{\text{tot}}$ &
$||\convection-\convection_{\sigma h}^{2,2}||$ & EOC & 
$||\pressure-\pressure_{\sigma h}^{1,2}||$ & EOC\\
\hline
1 & 4   & 16    & 2244      & 3.7974e-03 & ---  & 2.0593e-02 & --- \\
2 & 8   & 64    & 15816     & 4.4945e-04 & 3.08 & 2.5898e-03 & 2.99\\
3 & 16  & 256   & 118416    & 5.5129e-05 & 3.03 & 5.5087e-04 & 2.23\\
4 & 32  & 1024  & 915744    & 6.8603e-06 & 3.01 & 1.3405e-04 & 2.04\\
5 & 64  & 4096  & 7201344   & 8.5697e-07 & 3.00 & 3.3291e-05 & 2.01\\
6 & 128 & 16384 & 57115776  & 1.0713e-07 & 3.00 & 8.3039e-06 & 2.00\\
7 & 256 & 65536 & 454953216 & 1.3394e-08 & 3.00 & 2.0740e-06 & 2.00\\
\hline
\end{tabular}
}
%
%
\vspace{1.5ex}

\resizebox{0.75\linewidth}{!}{%
\begin{tabular}{c rrr | cc | cc}
\hline
\hline
$\ell$ & $N$ & $N_K$ & $N_{\text{DoF}}^{\text{tot}}$ &
$||\convection-\convection_{\sigma h}^{3,3}||$ & EOC & 
$||\pressure-\pressure_{\sigma h}^{2,3}||$ & EOC\\
\hline
1 & 4   & 16    & 6704       & 3.0085e-04 & ---  & 3.6349e-03 & --- \\
2 & 8   & 64    & 49248      & 1.9122e-05 & 3.98 & 4.9603e-04 & 2.87\\
3 & 16  & 256   & 377024     & 1.2073e-06 & 3.99 & 4.6902e-05 & 3.40\\
4 & 32  & 1024  & 2949504    & 7.5829e-08 & 3.99 & 4.6449e-06 & 3.34\\
5 & 64  & 4096  & 23331584   & 4.7494e-09 & 4.00 & 5.3029e-07 & 3.13\\
6 & 128 & 16384 & 185599488  & 2.9712e-10 & 4.00 & 6.4815e-08 & 3.03\\
\hline
\end{tabular}
}

\caption{Global convergence for 
$\convection_{\sigma h}^{2,2},\pressure_{\sigma h}^{1,2}$ in a 
cG(2)-dG(2),cG(1)-dG(2) and $\convection_{\sigma h}^{3,3},\pressure_{\sigma h}^{2,3}$ 
in a cG(3)-dG(3), cG(2)-dG(3) primal approximation for a time-dependent Stokes 
flow problem for Sec.~\ref{sec:6:1}.
$\ell$ denotes the refinement level, $N$ the total cells in time,
$N_K$ the cells in space on a slab,
$N_{\text{DoF}}^{\text{tot}}$ the total space-time degrees of freedom,
$||\cdot||$ the global $L^2(L^2)$-norm error and
EOC the experimental order of convergence.
}
\label{tab:1:BRdynamic:global}
\end{table}

The second part of the first example now serves to verify the higher-order 
implementation of the coupled problem. Therefore, the convection-diffusion-reaction 
transport problem \eqref{eq:1:transport_problem} is coupled with the 
Stokes flow problem \eqref{eq:3:stokes_problem} via the convection 
field solution $\convection_{\sigma h}$, using the exact solution $\convection$
given by Eq.~\eqref{eq:27:BRdynamic}. More precisely, we study problem 
\eqref{eq:1:transport_problem} with the given analytical solution
\begin{equation}
\label{eq:28:KB2}
\begin{array}{l@{\,}c@{\,}l}
u(\boldsymbol x, t) &:=&
u_1 \cdot u_2\,,\,\,
\boldsymbol x = (x_1, x_2)^\top \in \mathbb{R}^2 \text{ and }
t \in \mathbb{R}\,,\\[.5ex]
u_1(\boldsymbol x, t) &:=&
  (1 + a \cdot ( (x_1 - m_1(t))^2 + (x_2 - m_2(t))^2  ) )^{-1}\,,\\[.5ex]
u_2(t) &:=& \nu_1(t) \cdot s \cdot \arctan( \nu_2(t) )\,,
\end{array}
\end{equation}
with $m_1(t) := \frac{1}{2} + \frac{1}{4} \cos(2 \pi t)$ and
$m_2(t) := \frac{1}{2} + \frac{1}{4} \sin(2 \pi t)$, and,
$\nu_1(\hat t) := -1$,
$\nu_2(\hat t) := 5 \pi \cdot (4 \hat t - 1)$,
for $\hat t \in [0, 0.5)$ and
$\nu_1(\hat t) := 1$,
$\nu_2(\hat t) := 5 \pi \cdot (4 (\hat t-0.5) - 1)$,
for $\hat t \in [0.5, 1)$, $\hat t = t - k$,
$k \in \mathbb{N}_0$, and,
scalars $a = 50$ and $s=-\frac{1}{3}$.
The (analytic) solution \eqref{eq:28:KB2} mimics a counterclockwise rotating
cone which additionally changes its height and orientation over the period
$T=1$. Precisely, the orientation of the cone switches from negative to positive
while passing $t=0.25$ and from positive to negative while passing $t=0.75$.
The inhomogeneous Dirichlet boundary condition, the inhomogeneous initial condition
and the right-hand side forcing term $\transportforce$, are calculated from 
the given analytic solution \eqref{eq:28:KB2} and Eq. \eqref{eq:1:transport_problem}, 
where the latter uses the exact Stokes solution $\convection$ given by 
Eq.~\eqref{eq:27:BRdynamic}. Moreover, we note that the assembly of the transport
system matrix uses the approximated fully-discrete Stokes solution 
$\convection_{\sigma h}$ of \eqref{eq:18:B_v_p_psi_chi_eq_F_psi} that has to be
transferred to the spatial and temporal mesh of the transport problem, cf.
Rem.~\ref{rem:?:algorithm} in Sec.~\ref{sec:5:algorithm}.

Since we study the global space-time refinement behavior here, we restrict the 
convection-diffusion-reaction transport problem to a non-stabilized case, i.e. we set 
$\delta_0:=0$ within the local SUPG stabilization parameter 
$\delta_K = \delta_0 \cdot h_K$, where $h_K$ denotes the cell diameter of the
spatial mesh cell $K$.
Moreover, we set the diffusion coefficient $\varepsilon=1$ and choose a constant 
reaction coefficient $\alpha=1$.
\begin{table}[ht]
\centering

\resizebox{0.75\linewidth}{!}{%
\begin{tabular}{c rr | cc || rr | cc}
\hline
\hline
$\ell$ & $N$ & $N_K$ & $||\concentration-\concentration_{\tau h}^{1,1}||$ & EOC &
$N$ & $N_K$ & $||\convection-\convection_{\sigma h}^{2,0}||$ & EOC \\
\hline
1 & 4   & 4     & 8.4766e-02 & ---  & 4   & 16    & 4.5783e-03 & ---  \\ 
2 & 8   & 16    & 2.7780e-02 & 1.61 & 8   & 64    & 1.7196e-03 & 1.41 \\ 
3 & 16  & 64    & 9.1450e-03 & 1.60 & 16  & 256   & 9.9383e-04 & 0.79 \\ 
4 & 32  & 256   & 3.0372e-03 & 1.59 & 32  & 1024  & 5.6228e-04 & 0.82 \\ 
5 & 64  & 1024  & 7.7372e-04 & 1.97 & 64  & 4096  & 3.0374e-04 & 0.89 \\ 
6 & 128 & 4096  & 1.9407e-04 & 2.00 & 128 & 16384 & 1.5868e-04 & 0.94 \\ 
7 & 256 & 16384 & 4.9596e-05 & 1.97 & 256 & 65536 & 1.5868e-04 & 0.94 \\ 
\hline
\end{tabular}
}
%
%
\vspace{1.5ex}

\resizebox{0.75\linewidth}{!}{%
\begin{tabular}{c rr | cc | rr | cc}
\hline
\hline
$\ell$ & $N$ & $N_K$ & $||\concentration-\concentration_{\tau h}^{2,2}||$ & EOC &
$N$ & $N_K$ & $||\convection-\convection_{\sigma h}^{3,0}||$ & EOC \\
\hline
1 & 4   & 1    & 9.9045e-02 & ---  & 4   & 4     & 7.6009e-03 & ---  \\ 
2 & 8   & 4    & 4.8261e-02 & 1.04 & 8   & 16    & 1.7488e-03 & 2.12 \\
3 & 16  & 16   & 6.0634e-03 & 2.99 & 16  & 64    & 9.9237e-04 & 0.82 \\
4 & 32  & 64   & 1.0858e-03 & 2.48 & 32  & 256   & 5.6213e-04 & 0.82 \\
5 & 64  & 256  & 1.5131e-04 & 2.84 & 64  & 1024  & 3.0373e-04 & 0.89 \\
6 & 128 & 1024 & 2.0859e-05 & 2.86 & 128 & 4096  & 1.5868e-04 & 0.94 \\
7 & 256 & 4096 & 2.7041e-06 & 2.95 & 256 & 16384 & 8.1226e-05 & 0.97 \\
\hline
\end{tabular}
}

\caption{Global convergence for 
$\convection_{\sigma h}^{2,2},\pressure_{\sigma h}^{1,2}$ in a 
cG(2)-dG(2),cG(1)-dG(2) and $\convection_{\sigma h}^{3,3},\pressure_{\sigma h}^{2,3}$ 
in a cG(3)-dG(3), cG(2)-dG(3) primal approximation for a time-dependent Stokes 
flow problem for Sec.~\ref{sec:6:1}.
$\ell$ denotes the refinement level, $N$ the total cells in time,
$N_K$ the cells in space on a slab,
$N_{\text{DoF}}^{\text{tot}}$ the total space-time degrees of freedom,
$||\cdot||$ the global $L^2(L^2)$-norm error and
EOC the experimental order of convergence.
}
\label{tab:2:KB2BRdynamic:global}
\end{table}
The global space-time refinement behavior is illustrated by 
Tab.~\ref{tab:2:KB2BRdynamic:global} and nicely confirms
our results with respect to the expected EOCs for the solely solved transport 
problem obtained in \cite[Example 1]{Bause2021}, cf. columns four and five of 
Tab.~\ref{tab:2:KB2BRdynamic:global}.
Furthermore, with regard to the EOCs of the Stokes solution, we note that both
approximations cG(2)-dG(0) as well as cG(3)-dG(0) are restricted through the 
lowest order approximation in time, cf. columns eight and nine of 
Tab.~\ref{tab:2:KB2BRdynamic:global}.

\subsection{Example 2 (Space-time adaptivity studies for the coupled problem)}
\label{sec:6:2}
The second example serves to study the goal-oriented space-time adaptivity 
behavior of our algorithm introduced in Sec.~\ref{sec:5:algorithm}. 
More precisely, the transport problem is adaptively refined in space and time 
using an approximated Stokes solution $\convection_{\sigma h}$ on a coarser 
global refined mesh in space and time. 
The initial space-time meshes of the transport problem are once more refined 
compared to the initial meshes of the Stokes flow problem, cf. the first row of 
Tab.~\ref{tab:3:KB2BRdynamic:ada}.
The temporal and spatial mesh of the Stokes flow problem is refined globally if
the global $L^2(L^2)$-error $||\convection-\convection_{\sigma h}^{2,0}||$ is 
larger than its counterpart $||u-u_{\tau h}^{1,0}||$ or rather $||u-u_{\tau h}^{1,1}||$ 
for the transport problem (cf. columns five and nine of Tab.~\ref{tab:3:KB2BRdynamic:ada}
and Tab.~\ref{tab:4:KB2BRdynamic:ada1121}, respectively.).

We study problem \eqref{eq:1:transport_problem} and \eqref{eq:3:stokes_problem}
with the given analytical solutions \eqref{eq:28:KB2} and \eqref{eq:27:BRdynamic},
respectively, with the same settings as given in Sec.~\ref{sec:6:1}.
Our target quantity for the transport problem is chosen to control the global 
$L^2(L^2)$-error of $e$, $e = \concentration - \concentration_{\tau h}$, in 
space and time, given by
\begin{equation}
\label{eq:29:L2Goal}
J(\varphi)= \frac{1}{\|e\|_{(0,T)\times\Omega}}
\displaystyle\int_I(\varphi,e)\mathrm{d}t\,,
\quad \mathrm{with} \;\; \|\cdot\|_{(0,T)\times\Omega}
=\left(\int_I(\cdot,\cdot)\;\mathrm{d}t\right)^{\frac{1}{2}}\,.
\end{equation}
The tuning parameters of the goal-oriented adaptive Algorithm given in
Sec.~\ref{sec:5:algorithm} are chosen here in a way to balance
automatically the potential misfit of the spatial and temporal errors as
\begin{displaymath}
\theta_h^\textnormal{top} = 0.5 \cdot
\left| \frac{\eta_h}{|\eta_h| + |\eta_\tau|} \right|\,,\quad
\theta_h^\textnormal{bottom} = 0\quad
\textnormal{and}\quad
\theta_\tau^\textnormal{top} = 0.5 \cdot
\left| \frac{\eta_\tau}{|\eta_h| + |\eta_\tau|} \right|\,.
\end{displaymath}
For measuring the accuracy of the error estimator, we will study the so-called 
effectivity index given by
\begin{equation}
\label{eq:30:Ieff}
\mathcal{I}_{\textnormal{eff}} = \left|\frac{\eta_{\tau}+\eta_{h}}{J(\concentration)
-J(\concentration_{\tau h})}\right|
\end{equation}
as the ratio of the estimated error over the exact error.
Desirably, the index $\mathcal{I}_{\textnormal{eff}}$ should be close to one.

In Tab.~\ref{tab:3:KB2BRdynamic:ada} and Tab.~\ref{tab:4:KB2BRdynamic:ada1121} 
we present the development of the total discretization error 
$J(e)=\|e\|_{(0,T)\times\Omega}$ for \eqref{eq:29:L2Goal},
the approximated spatial and temporal error estimators $\eta_h$ and $\eta_{\tau}$ 
as well as the effectivity index $\mathcal{I}_{\mathrm{eff}}$ during an adaptive 
refinement process for two different primal and dual solution pairings 
$\{u_{\tau h},z_{\tau h}\}$:
cG(1)-dG(0)/cG(2)-dG(0), cG(1)-dG(1)/cG(2)-dG(1)
of the transport problem. Moreover, the development of the total 
discretization error $||\convection-\convection_{\sigma h}^{2,0}||$ for the Stokes 
flow solution on a global refined mesh in space and time and the corresponding
number of slabs and spatial cells is displayed.
Thereby, $\convection_{\sigma h}^{2,0}$ corresponds to a Stokes solution 
approximation in a cG($2$)-dG($0$) discretization.
We use an approximation of the temporal weights by a higher-order extrapolation 
strategy using Gauss-Lobatto quadrature points.
Here and in the following, $\ell$ denotes the refinement level or DWR loop,
$N$ or $N^{\text{F}}$ the total cells in time, $N_{K}^{\text{max}}$ or
$N_{K}^{\text{F},\text{max}}$ the number of spatial cells on the finest spatial 
mesh within the current loop, and $N_{\text{DoF}}^{\text{tot}}$ or
$N_{\text{DoF}}^{\text{F},\text{tot}}$ the total space-time degrees of freedom
of the transport or Stokes flow problem, respectively. 

Regarding the accuracy of the underlying error estimator, as given by the last 
column of Tab.~\ref{tab:3:KB2BRdynamic:ada} or Tab.~\ref{tab:4:KB2BRdynamic:ada1121},
respectively, we observe a good quantitative 
estimation of the discretization error as the respective effectivity index 
increases getting close to one. 
With regard to efficiency reasons for a space-time adaptive algorithm, it is 
essential to ensure an equilibrated reduction of the temporal as well as spatial
discretization error, cf. \cite[Sec. 3.3]{Besier2012}.
Referring to this, we point out a good equilibration of the spatial and 
temporal error indicators $\eta_{h}$ and $\eta_{\tau}$ in the course of the 
refinement process (columns ten and eleven of Tab.~\ref{tab:3:KB2BRdynamic:ada}
and Tab.~\ref{tab:4:KB2BRdynamic:ada1121}).
\begin{table}[h!]
\centering

\resizebox{0.95\linewidth}{!}{%
\begin{tabular}{c | rcrc | rrr | c | ccc | c}
\hline
\hline
DWR & \multicolumn{4}{c|}{Stokes Flow} & \multicolumn{8}{c}{Transport}
\\
\hline
$\ell$ & $N^{\text{F}}$ & $N_K^{\text{F},\text{max}}$ & $N_{\text{DoF}}^{\text{F},\text{tot}}$ & 
$||\convection-\convection_{\sigma h}^{2,0}||$ &
$N$ & $N_K^{\text{max}}$ & $N_{\text{DoF}}^{\text{tot}}$ & $\|e^{1,0,2,0}\|$ & 
${\eta}_h$ & ${\eta}_\tau$ & ${\eta}_{\tau h}$ &
$\mathcal{I}_{\textnormal{eff}}$
\\
\hline
1  &  5 &  4 &   295 & 1.966e-02 &  10 &  16 &   250 & 5.245e-02 & 4.270e-03 & 3.335e-04 & 4.604e-03 & 0.09\\ 
2  &    &    &       & 1.966e-02 &  10 &  40 &   438 & 4.681e-02 & 7.122e-04 & 3.547e-03 & 4.260e-03 & 0.09\\ 
3  &    &    &       & 1.966e-02 &  14 &  40 &   616 & 1.729e-02 & 3.587e-03 & 1.594e-03 & 5.180e-03 & 0.30\\ 
4  & 10 & 16 &  1870 & 4.099e-03 &  19 &  88 &  1913 & 1.100e-02 & 1.607e-03 & 3.784e-03 & 5.391e-03 & 0.49\\ 
5  &    &    &       & 4.099e-03 &  26 & 160 &  4074 & 7.139e-03 & 6.659e-04 & 3.638e-03 & 4.304e-03 & 0.60\\ 
6  &    &    &       & 4.099e-03 &  36 & 160 &  5534 & 5.036e-03 & 9.436e-04 & 2.585e-03 & 3.528e-03 & 0.70\\ 
7  &    &    &       & 4.099e-03 &  50 & 268 & 11954 & 3.439e-03 & 4.237e-04 & 2.146e-03 & 2.570e-03 & 0.75\\ 
8  & 20 & 64 & 13180 & 9.684e-04 &  70 & 268 & 16752 & 2.585e-03 & 5.376e-04 & 1.453e-03 & 1.991e-03 & 0.77\\ 
9  &    &    &       & 9.684e-04 &  98 & 448 & 37900 & 1.844e-03 & 2.343e-04 & 1.183e-03 & 1.417e-03 & 0.77\\ 
10 & 20 & 64 & 13180 & 9.684e-04 & 137 & 448 & 52937 & 1.402e-03 & 2.908e-04 & 8.384e-04 & 1.129e-03 & 0.81\\ 
\hline
\end{tabular}
}

\vskip-1ex
\caption{Adaptive refinement in the transport problem including effectivity indices
for goal quantity \eqref{eq:29:L2Goal},
with $\varepsilon = 1$,$\delta_0=0$, and $\omega = 3$ for Sec.~\ref{sec:6:2}
using a Stokes solution $\convection_{\sigma h}^{2,0}$ corresponding to a
cG(2)-dG(0) approximation on a global refined mesh in space and time.
$e^{1,0,2,0}$ corresponds to
the adaptive solution approximation
$\concentration_{\tau h}^{1,0}$ in cG(1)-dG(0)
and dual solution approximation
$\dualz_{\tau h}^{2,0}$ in cG(2)-dG(0).}

\label{tab:3:KB2BRdynamic:ada}
\end{table}

\begin{table}[h!]
\centering

\resizebox{0.95\linewidth}{!}{%
\begin{tabular}{c | rcrc | rrr | c | ccc | c}
\hline
\hline
DWR & \multicolumn{4}{c|}{Stokes Flow} & \multicolumn{8}{c}{Transport}
\\
\hline
$\ell$ & $N^{\text{F}}$ & $N_K^{\text{F},\text{max}}$ & $N_{\text{DoF}}^{\text{F},\text{tot}}$ & 
$||\convection-\convection_{\sigma h}^{2,0}||$ &
$N$ & $N_K^{\text{max}}$ & $N_{\text{DoF}}^{\text{tot}}$ & $\|e^{1,1,2,1}\|$ & 
${\eta}_h$ & ${\eta}_\tau$ & ${\eta}_{\tau h}$ &
$\mathcal{I}_{\textnormal{eff}}$
\\
\hline
1  &  5 &   4 &   295 & 1.966e-02 &  20 &  16 &   1000 & 2.547e-02 & 2.139e-02 & 1.231e-04 & 2.126e-02 & 0.83\\ 
2  &    &     &       & 1.966e-02 &  20 &  28 &   1640 & 1.158e-02 & 1.279e-02 & 4.406e-03 & 1.719e-02 & 1.48\\ 
3  & 10 &  16 &  1870 & 4.099e-03 &  20 &  76 &   2956 & 7.730e-03 & 3.986e-03 & 5.234e-03 & 9.220e-03 & 1.19\\ 
4  &    &     &       & 4.099e-03 &  28 & 124 &   6468 & 4.340e-03 & 3.634e-03 & 4.309e-03 & 7.943e-03 & 1.83\\ 
5  &    &     &       & 4.099e-03 &  39 & 172 &  11694 & 2.840e-03 & 2.314e-03 & 3.446e-03 & 5.760e-03 & 2.02\\ 
6  & 20 &  64 & 13180 & 9.684e-04 &  54 & 232 &  20348 & 1.945e-03 & 8.625e-04 & 3.172e-03 & 4.035e-03 & 2.07\\ 
7  &    &     &       & 9.684e-04 &  64 & 232 &  24600 & 1.889e-03 & 2.207e-04 & 2.273e-03 & 2.052e-03 & 1.08\\ 
8  &    &     &       & 9.684e-04 &  75 & 232 &  28022 & 1.871e-03 & 5.367e-05 & 2.367e-03 & 2.313e-03 & 1.23\\ 
9  &    &     &       & 9.684e-04 & 147 & 316 &  66570 & 1.469e-03 & 2.111e-03 & 6.047e-04 & 1.506e-03 & 1.02\\ 
10 & 40 & 256 & 98680 & 4.691e-04 & 283 & 532 & 220298 & 8.416e-04 & 9.537e-04 & 8.517e-05 & 8.685e-04 & 1.03\\ 
\hline
\end{tabular}
}

\vskip-1ex
\caption{Adaptive refinement in the transport problem including effectivity indices
for goal quantity \eqref{eq:29:L2Goal},
with $\varepsilon = 1$,$\delta_0=0$, and $\omega = 3$ for Sec.~\ref{sec:6:2}
using a Stokes solution $\convection_{\sigma h}^{2,0}$ corresponding to a
cG(2)-dG(0) approximation on a global refined mesh in space and time.
$e^{1,1,2,1}$ corresponds to
the adaptive solution approximation
$\concentration_{\tau h}^{1,0}$ in cG(1)-dG(1)
and dual solution approximation
$\dualz_{\tau h}^{2,1}$ in cG(2)-dG(1).}

\label{tab:4:KB2BRdynamic:ada1121}
\end{table}

\begin{figure}[h!]
\begin{minipage}{\linewidth}
\centering
\begin{tikzpicture}
\begin{axis}[%
width=4.in,
height=1.in,
scale only axis,
xlabel={t},
ylabel={\textcolor{navyblue}{$\tau_K(I_n^{1})$}, \textcolor{HSUred}{$\sigma_K(I_n^{\textnormal{F},1})$}},
xmin=0.0,
xmax=1.0,
ymin=0.,
ymax=0.22,
yminorticks=true,
yticklabels={0.0,0, 0.05 , 0.1, 0.15 , 0.2},
]
\addplot [
color=navyblue,
solid,
line width=1.5pt,
mark=*,
mark size = 1.,
only marks,
mark options={solid,navyblue}
]
table[row sep=crcr]{
0.1 0.1 \\
0.2 0.1 \\
0.3 0.1 \\
0.4 0.1 \\
0.5 0.1 \\
0.6 0.1 \\
0.7 0.1 \\
0.8 0.1 \\
0.9 0.1 \\
1 0.1 \\
};
\addplot [
color=HSUred,
solid,
line width=0.5pt,
mark=*,
mark size = 1.,
only marks,
mark options={fill=HSUred}
]
table[row sep=crcr]{
0.2 0.2 \\
0.4 0.2 \\
0.6 0.2 \\
0.8 0.2 \\
1 0.2 \\
};
\end{axis}
\end{tikzpicture}
\end{minipage}

\begin{minipage}{\linewidth}
\centering
\begin{tikzpicture}
\begin{axis}[%
width=4.in,
height=1.in,
scale only axis,
/pgf/number format/.cd, 1000 sep={},
xlabel={t},
ylabel={\textcolor{navyblue}{$\tau_K(I_n^{5})$}, \textcolor{HSUred}{$\sigma_K(I_n^{\textnormal{F},5})$}},
xmin=0.0,
xmax=1.0,
ymin=0.,
ymax=0.12,
yminorticks=true,
yticklabels={0.0,0, 0.05 , 0.1},
]

\addplot [
color=navyblue,
solid,
line width=1.0pt,
mark=*,
mark size = 1.5,
only marks,
mark options={solid,navyblue}
]
table[row sep=crcr]{
0.1 0.1 \\
0.15 0.05 \\
0.2 0.05 \\
0.225 0.025 \\
0.2375 0.0125 \\
0.25 0.0125 \\
0.2625 0.0125 \\
0.275 0.0125 \\
0.3 0.025 \\
0.35 0.05 \\
0.4 0.05 \\
0.5 0.1 \\
0.55 0.05 \\
0.6 0.05 \\
0.65 0.05 \\
0.7 0.05 \\
0.725 0.025 \\
0.7375 0.0125 \\
0.75 0.0125 \\
0.7625 0.0125 \\
0.775 0.0125 \\
0.8 0.025 \\
0.85 0.05 \\
0.9 0.05 \\
0.95 0.05 \\
1 0.05 \\
};

\addplot [
color=HSUred,
solid,
line width=0.5pt,
mark=*,
mark size = 1.,
only marks,
mark options={fill=HSUred}
]
table[row sep=crcr]{
0.1 0.1 \\
0.2 0.1 \\
0.3 0.1 \\
0.4 0.1 \\
0.5 0.1 \\
0.6 0.1 \\
0.7 0.1 \\
0.8 0.1 \\
0.9 0.1 \\
1 0.1 \\
};
\end{axis}
\end{tikzpicture}
\end{minipage}
\begin{minipage}{\linewidth}
\centering
\begin{tikzpicture}
\begin{axis}[%
width=4.in,
height=1.in,
scale only axis,
/pgf/number format/.cd, 1000 sep={},
xlabel={t},
ylabel={\textcolor{navyblue}{$\tau_K(I_n^{10})$}, \textcolor{HSUred}{$\sigma_K(I_n^{\textnormal{F},10})$}},
xmin=0.0,
xmax=1.0,
yminorticks=true,
yticklabels={0.0,0, \phantom{-,-}2 , \phantom{-,-}4 },
]

\addplot [
color=navyblue,
solid,
line width=1.0pt,
mark=*,
mark size = 1.5,
only marks,
mark options={solid,navyblue}
]
table[row sep=crcr]{
0.0125 0.0125 \\
0.01875 0.00625 \\
0.025 0.00625 \\
0.0375 0.0125 \\
0.05 0.0125 \\
0.0625 0.0125 \\
0.075 0.0125 \\
0.0875 0.0125 \\
0.1 0.0125 \\
0.10625 0.00625 \\
0.1125 0.00625 \\
0.125 0.0125 \\
0.1375 0.0125 \\
0.15 0.0125 \\
0.15625 0.00625 \\
0.1625 0.00625 \\
0.175 0.0125 \\
0.1875 0.0125 \\
0.2 0.0125 \\
0.20625 0.00625 \\
0.2125 0.00625 \\
0.21875 0.00625 \\
0.225 0.00625 \\
0.228125 0.003125 \\
0.23125 0.003125 \\
0.234375 0.003125 \\
0.2375 0.003125 \\
0.240625 0.003125 \\
0.242188 0.0015625 \\
0.24375 0.0015625 \\
0.245312 0.0015625 \\
0.246875 0.0015625 \\
0.248438 0.0015625 \\
0.25 0.0015625 \\
0.251563 0.0015625 \\
0.253125 0.0015625 \\
0.254687 0.0015625 \\
0.25625 0.0015625 \\
0.257812 0.0015625 \\
0.259375 0.0015625 \\
0.2625 0.003125 \\
0.265625 0.003125 \\
0.26875 0.003125 \\
0.271875 0.003125 \\
0.275 0.003125 \\
0.28125 0.00625 \\
0.2875 0.00625 \\
0.29375 0.00625 \\
0.3 0.00625 \\
0.3125 0.0125 \\
0.325 0.0125 \\
0.3375 0.0125 \\
0.35 0.0125 \\
0.3625 0.0125 \\
0.375 0.0125 \\
0.3875 0.0125 \\
0.4 0.0125 \\
0.4125 0.0125 \\
0.425 0.0125 \\
0.43125 0.00625 \\
0.4375 0.00625 \\
0.44375 0.00625 \\
0.45 0.00625 \\
0.4625 0.0125 \\
0.475 0.0125 \\
0.48125 0.00625 \\
0.4875 0.00625 \\
0.5 0.0125 \\
0.5125 0.0125 \\
0.525 0.0125 \\
0.5375 0.0125 \\
0.55 0.0125 \\
0.5625 0.0125 \\
0.575 0.0125 \\
0.5875 0.0125 \\
0.6 0.0125 \\
0.6125 0.0125 \\
0.625 0.0125 \\
0.6375 0.0125 \\
0.65 0.0125 \\
0.65625 0.00625 \\
0.6625 0.00625 \\
0.675 0.0125 \\
0.6875 0.0125 \\
0.69375 0.00625 \\
0.7 0.00625 \\
0.70625 0.00625 \\
0.7125 0.00625 \\
0.71875 0.00625 \\
0.725 0.00625 \\
0.728125 0.003125 \\
0.73125 0.003125 \\
0.734375 0.003125 \\
0.7375 0.003125 \\
0.740625 0.003125 \\
0.742188 0.0015625 \\
0.74375 0.0015625 \\
0.745313 0.0015625 \\
0.746875 0.0015625 \\
0.748437 0.0015625 \\
0.75 0.0015625 \\
0.751563 0.0015625 \\
0.753125 0.0015625 \\
0.754687 0.0015625 \\
0.75625 0.0015625 \\
0.757812 0.0015625 \\
0.759375 0.0015625 \\
0.760937 0.0015625 \\
0.7625 0.0015625 \\
0.765625 0.003125 \\
0.76875 0.003125 \\
0.775 0.00625 \\
0.78125 0.00625 \\
0.7875 0.00625 \\
0.79375 0.00625 \\
0.8 0.00625 \\
0.8125 0.0125 \\
0.825 0.0125 \\
0.83125 0.00625 \\
0.8375 0.00625 \\
0.84375 0.00625 \\
0.85 0.00625 \\
0.85625 0.00625 \\
0.8625 0.00625 \\
0.875 0.0125 \\
0.8875 0.0125 \\
0.9 0.0125 \\
0.9125 0.0125 \\
0.925 0.0125 \\
0.9375 0.0125 \\
0.95 0.0125 \\
0.9625 0.0125 \\
0.975 0.0125 \\
0.98125 0.00625 \\
};

\addplot [
color=HSUred,
solid,
line width=0.5pt,
mark=*,
mark size = 1.,
only marks,
mark options={fill=HSUred}
]
table[row sep=crcr]{
0.05 0.05 \\
0.1 0.05 \\
0.15 0.05 \\
0.2 0.05 \\
0.25 0.05 \\
0.3 0.05 \\
0.35 0.05 \\
0.4 0.05 \\
0.45 0.05 \\
0.5 0.05 \\
0.55 0.05 \\
0.6 0.05 \\
0.65 0.05 \\
0.7 0.05 \\
0.75 0.05 \\
0.8 0.05 \\
0.85 0.05 \\
0.9 0.05 \\
0.95 0.05 \\
1 0.05 \\
};
\end{axis}
\end{tikzpicture}
\end{minipage}

\caption{Distribution of the temporal step size $\tau_K$ of the transport problem
(adaptive) and $\sigma_K$ of the Stokes flow problem (global) over the time 
interval $I=(0,T]$ for the initial (1) and after 5 and 10 DWR-loops, 
corresponding to Tab.~\ref{tab:3:KB2BRdynamic:ada}.}
\label{fig:7:DistributionTauSigmaEx2}
\end{figure}
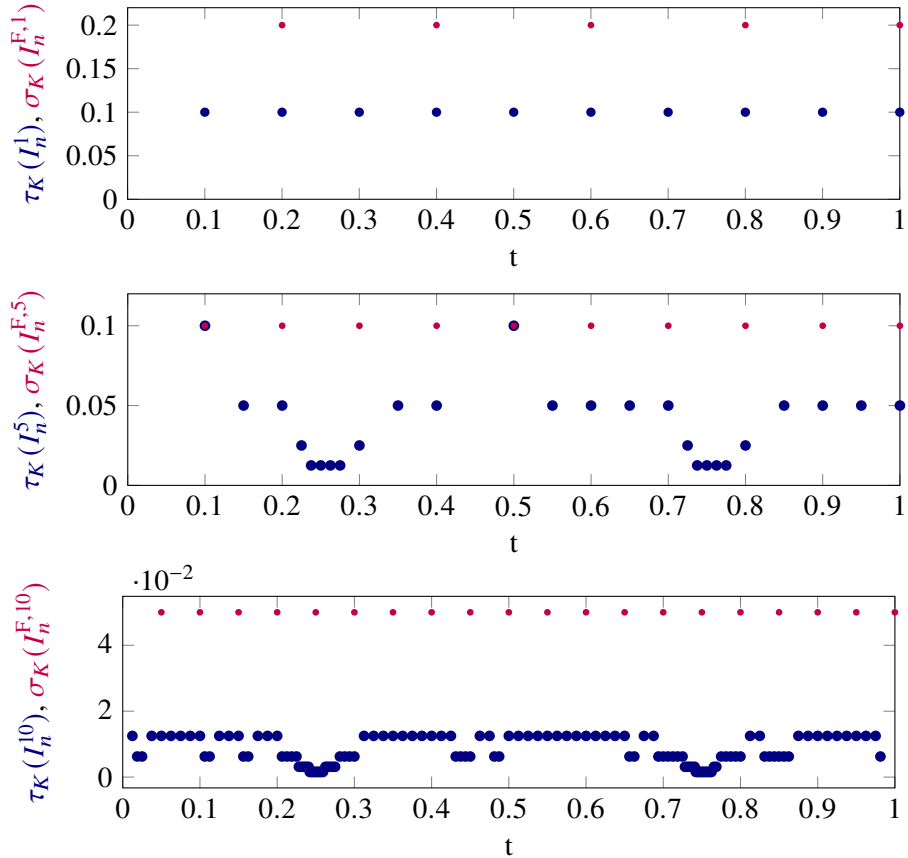

Finally, in Fig.~\ref{fig:7:DistributionTauSigmaEx2} we visualize exemplary the 
distribution of the  adaptively determined time cell lengths $\tau_K$ of 
$\mathcal{T}_{\tau,n}$, used for the transport problem, as well as the 
distribution of the globally determined time cell lengths $\sigma_K$ of 
$\mathcal{T}_{\sigma,n}$, used for the Stokes flow problem, over the whole time 
interval $I$ for different DWR refinement loops, corresponding to 
Tab.~\ref{tab:3:KB2BRdynamic:ada}.
The initial temporal meshes for the transport and Stokes flow problem are chosen
fulfilling the requirements presented in Sec.~\ref{sec:2:2:multirate} and 
Fig.~\ref{fig:1:multirate_time_scales}. 
While the time steps for the transport problem become smaller when the cone is
changing its orientation ($t=0.25$ and $t=0.75$), the time steps for the Stokes 
flow problem stay comparatively large in the course of the refinement process,
cf. the last two plots of Fig.~\ref{fig:7:DistributionTauSigmaEx2}.
Away from the time points of orientation change, the temporal mesh of the 
transport problem is almost equally decomposed.
This behavior seems natural for a global acting target quantity \eqref{eq:29:L2Goal}
and nicely confirms our approach of an efficient temporal approximation of a 
rapidly changing transport coupled with a slowly varying viscous flow.

%

\subsection{Example 3 (Transport in a channel)}
\label{sec:6:3}

\begin{figure}
\centering

\includegraphics[width=.44\linewidth]{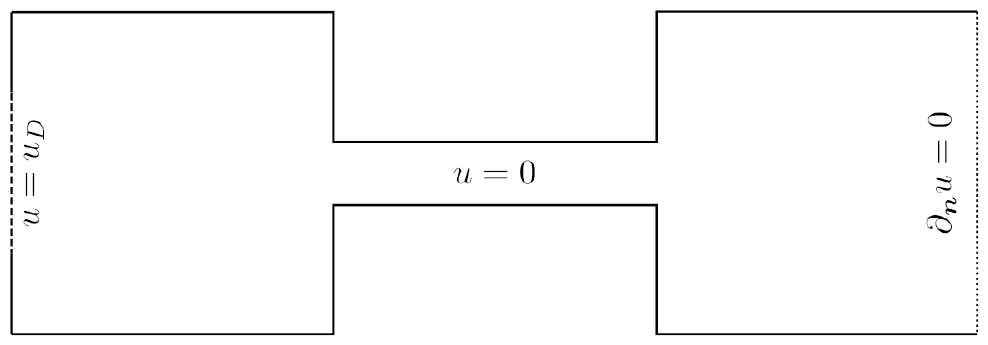}
~~
\includegraphics[width=.44\linewidth]{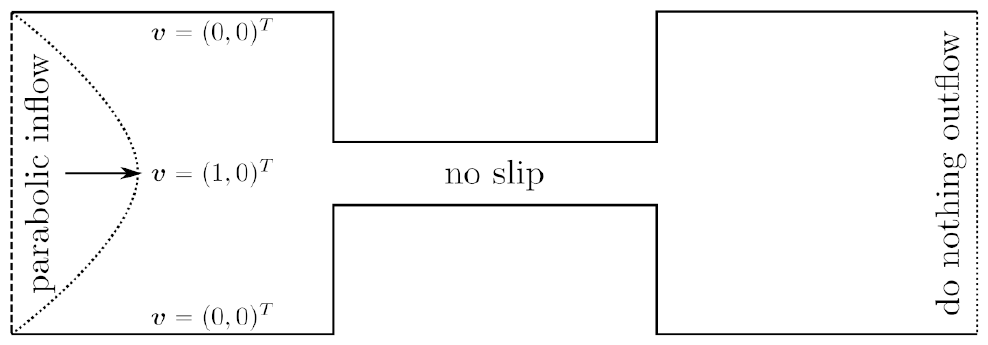}

\caption{Boundary colorization for the convection-diffusion-reaction problem (left)
and the coupled Stokes flow problem (right) for Sec.~\ref{sec:6:3}.}

\label{fig:8:boundary}
\end{figure}

In this example we simulate a convection-dominated transport with goal-oriented
adaptivity of a species through a channel with a constraint. The domain and its
boundary colorization are presented by Fig. \ref{fig:8:boundary}. Precisely,
the spatial domain is composed of two unit squares and a constraint in the middle
which restricts the channel height by a factor of 5.
Precisely, $\Omega = (-1,0)\times(-0.5,0.5) \cup (0,1)\times(-0.1,0.1) \cup
(1,2)\times(-0.5,0.5)$ with an initial cell diameter of $h=\sqrt{2 \cdot 0.025^2}$.
The time domain is set to $I=(0,2.5)$ with an initial $\tau=0.1$ for the transport
and $\sigma=2.5$ for the Stokes flow problem 
for the initialization of the slabs for the first loop $\ell=1$.
This choice has been made to compare the results to Example 2 in \cite{Bause2021}, 
where a quasi-stationary Stokes flow solution $\convection_h$ was used.
We approximate the primal solution $\concentration_{\tau h}^{1,0}$ with the
cG(1)-dG(0) method, the dual solution $\dualz_{\tau h}^{2,0}$ with the
cG(2)-dG(0) method and the Stokes flow solution $\convection_{\sigma h}^{2,0}$
with the cG(2)-dG(0) method.

\begin{figure}
\centering

\includegraphics[width=.6\linewidth]{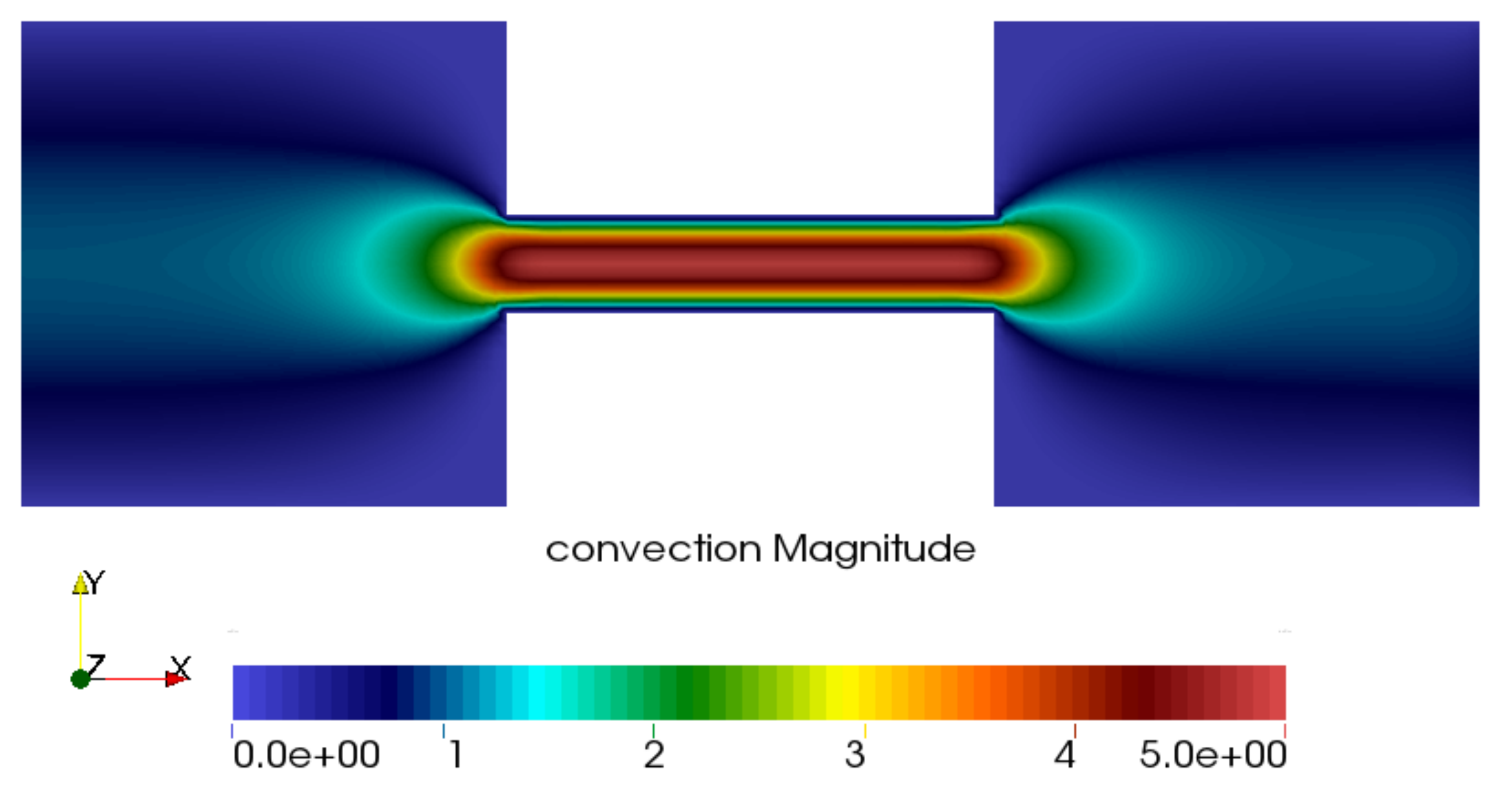}

\caption{Convection $\convection_{\sigma h}$ solution of the Stokes problem on 
one slab with a sufficiently globally refined spatial mesh with $Q_2$-$Q_1$ 
finite elements for Sec.~\ref{sec:6:3}.
On the left boundary a parabolic inflow profile in the positive x-direction
with unit magnitude is prescribed for the convection $\convection$.}

\label{fig:9:convection}
\end{figure}

The target quantity is
\begin{equation}
\label{eq:31:goal-mean}
J(u)= \frac{1}{|T|\cdot |\Omega|}
\displaystyle\int_I\int_\Omega u(\boldsymbol x, t)\, \mathrm{d}\boldsymbol{x}\mathrm{d}t\,.
\end{equation}
The transport of the species, which enters the domain on the left with an
inhomogeneous and time-dependent Dirichlet boundary condition and
leaves the domain on the right through a homogeneous Neumann boundary condition,
is driven by the convection with magnitudes between 0 and 5 as displayed in
Fig.~\ref{fig:9:convection}. The diffusion coefficient has the constant and small
value of $\varepsilon=10^{-4}$ and the reaction coefficient is chosen $\alpha=0.1$. 
%
The local SUPG stabilization coefficient is here set to
$\delta_K = \delta_0 \cdot h_K$, $\delta_0=0$,
i.e. a vanishing stabilization here.
The initial value function $\concentration_0=0$ as well as the forcing term
$\transportforce=0$ are homogeneous.
The Dirichlet boundary function value is homogeneous on $\Gamma_D$ except for the
line $(-1,-1) \times (-0.4,0.4)$ and time $0 \leq t \leq 0.2$ where the constant
value
\begin{displaymath}
\concentration(\boldsymbol{x},t)=1
\end{displaymath}
is prescribed on the solution.
The viscosity is set to $\viscosity=1$.
The tuning parameters of the goal-oriented adaptive Algorithm given in
Sec.~\ref{sec:5:algorithm} are chosen here in a way to balance
automatically the potential misfit of the spatial and temporal errors as
$\theta_h^\textnormal{bottom} = 0$,
\begin{equation}
\label{eq:32:balancing}
\theta_h^\textnormal{top} = \frac{1}{2} \cdot
\min\left\{ \left| \frac{\eta_h}{|\eta_h| + |\eta_\tau|} \right|\,, 1\right\}\quad
\textnormal{and}\quad
\theta_\tau^\textnormal{top} = \frac{1}{2} \cdot
\min\left\{ \left| \frac{\eta_\tau}{|\eta_h| + |\eta_\tau|} \right|\,, 1\right\}\,.
\end{equation}

\begin{figure}
\centering

\includegraphics[width=.44\linewidth]{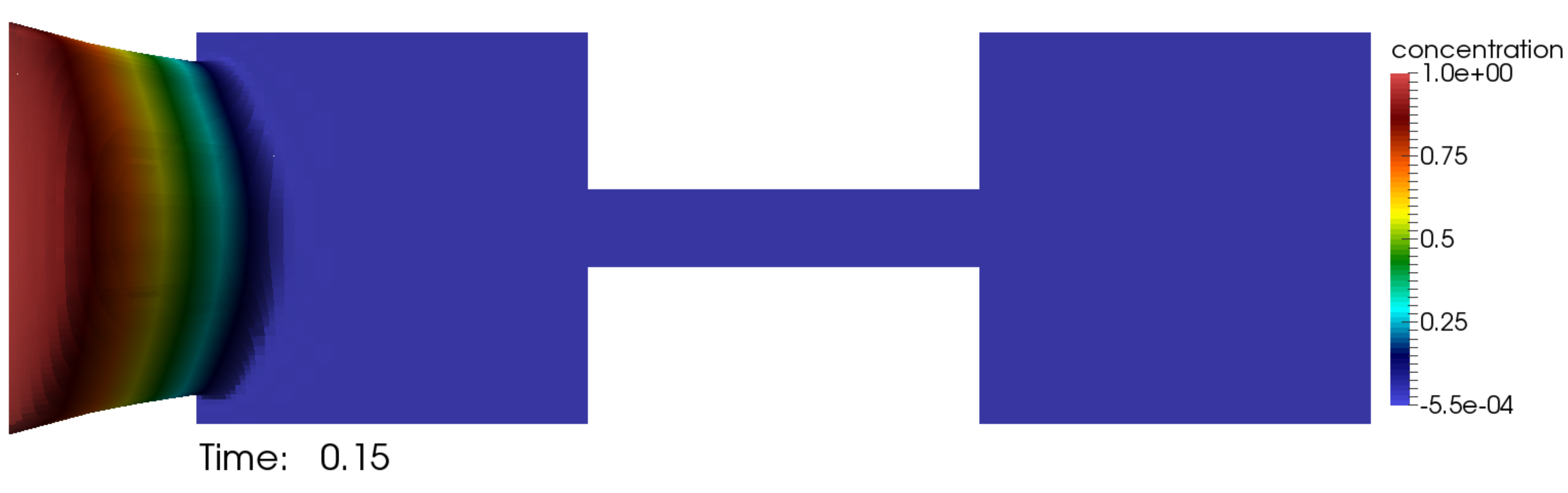}
~~
\includegraphics[width=.44\linewidth]{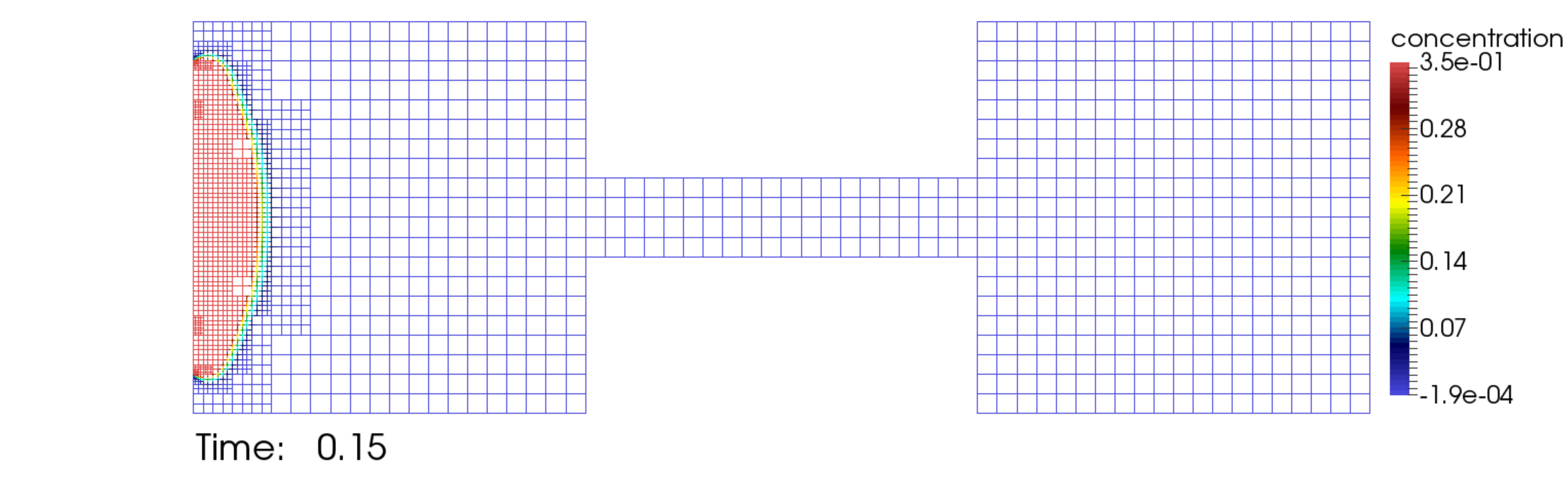}

\includegraphics[width=.44\linewidth]{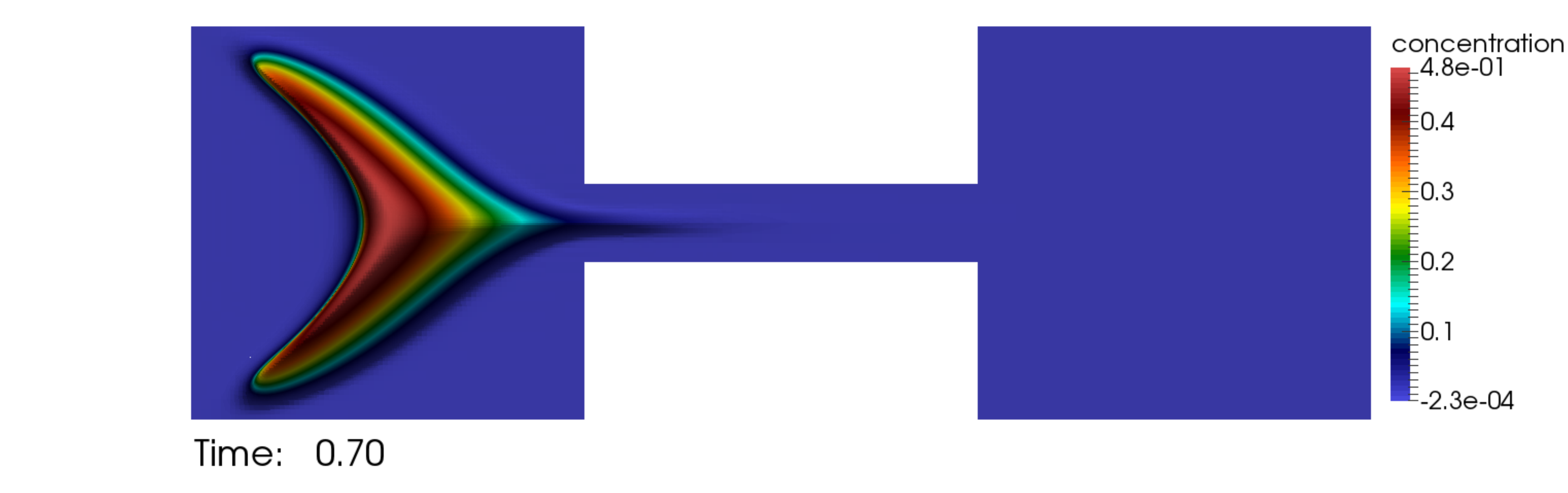}
~~
\includegraphics[width=.44\linewidth]{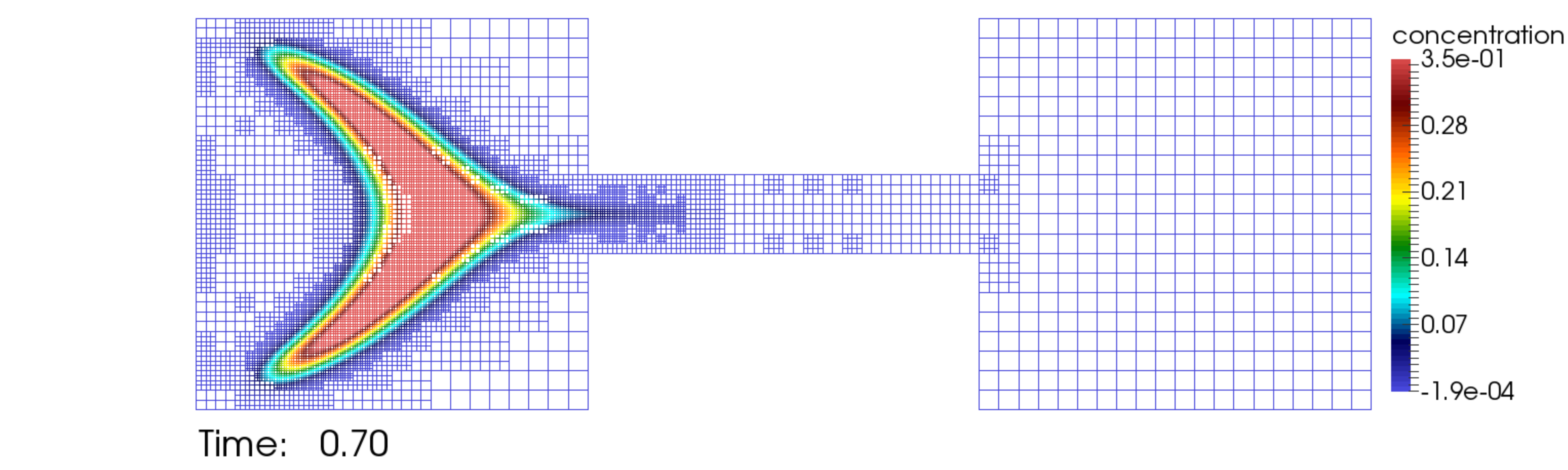}

\includegraphics[width=.44\linewidth]{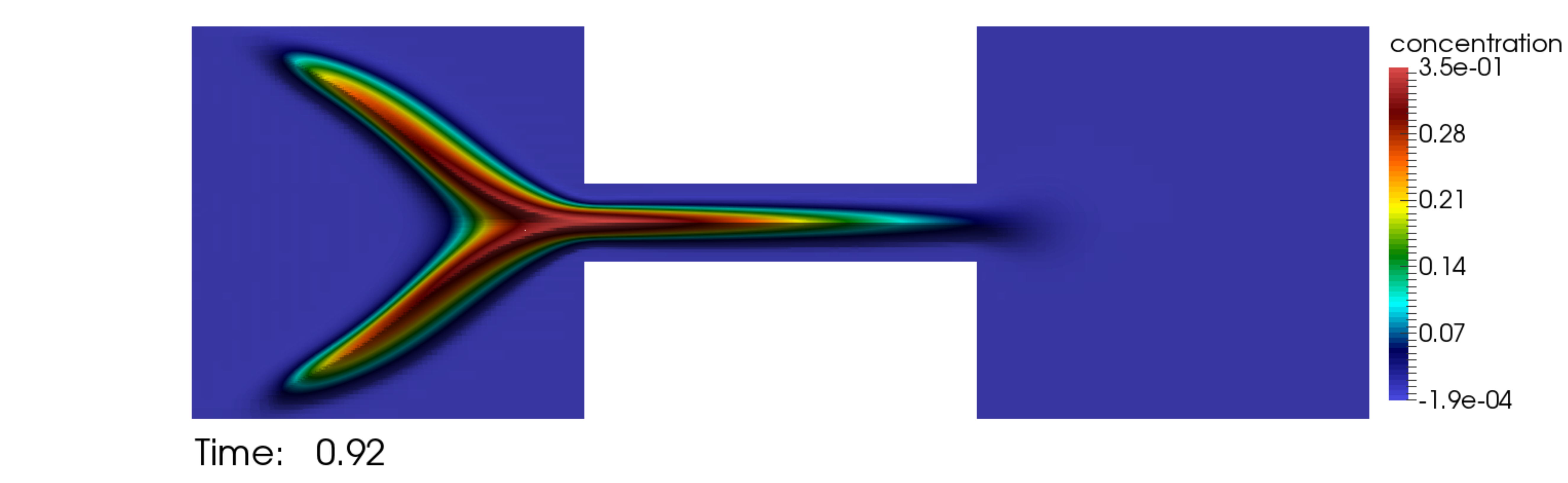}
~~
\includegraphics[width=.44\linewidth]{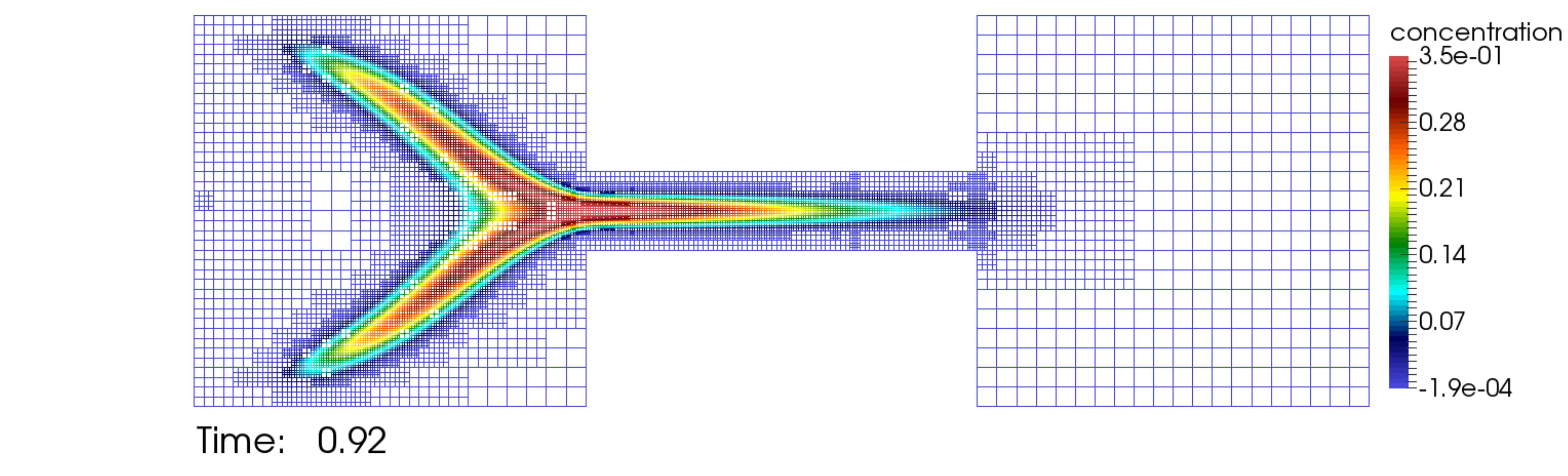}

\includegraphics[width=.44\linewidth]{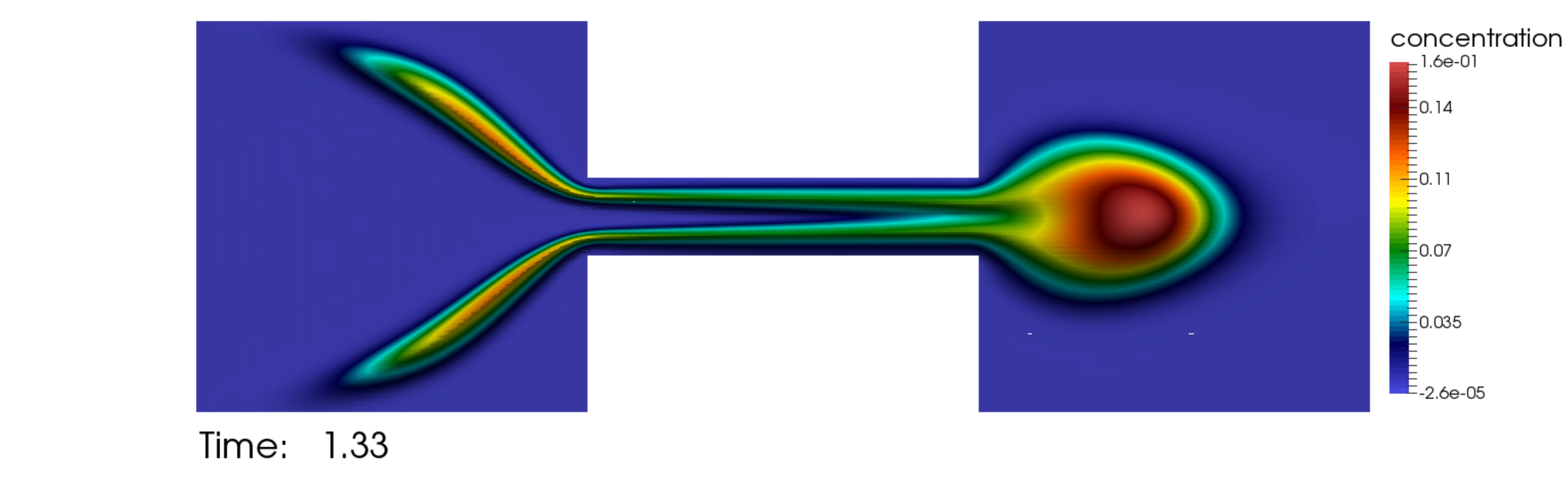}
~~
\includegraphics[width=.44\linewidth]{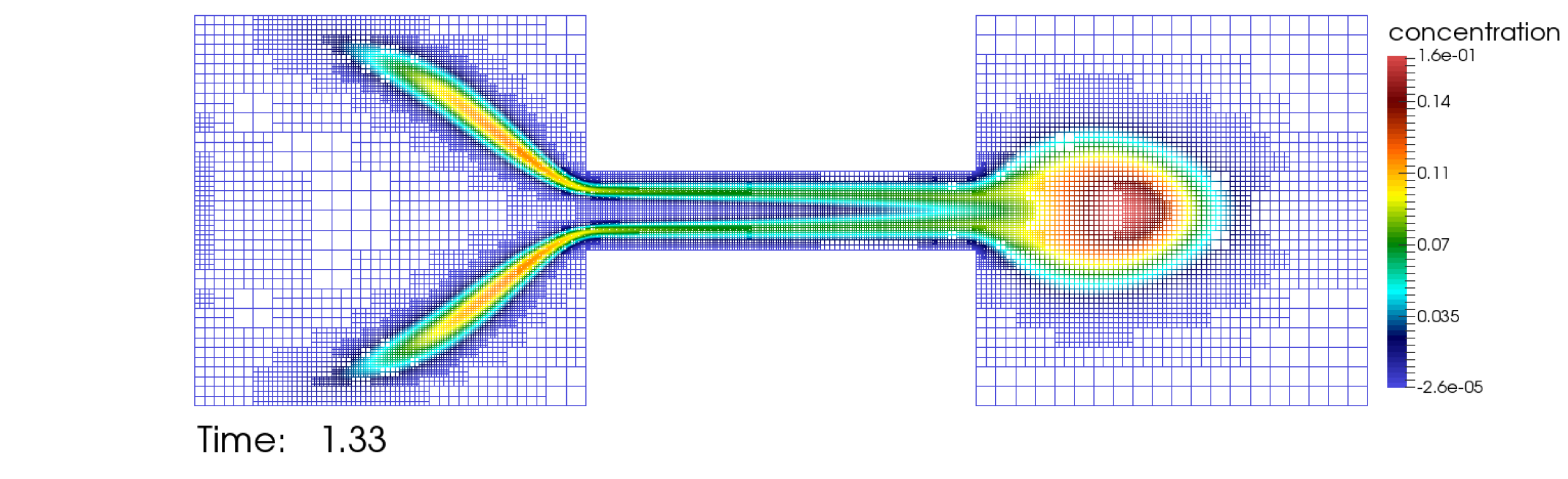}

\includegraphics[width=.44\linewidth]{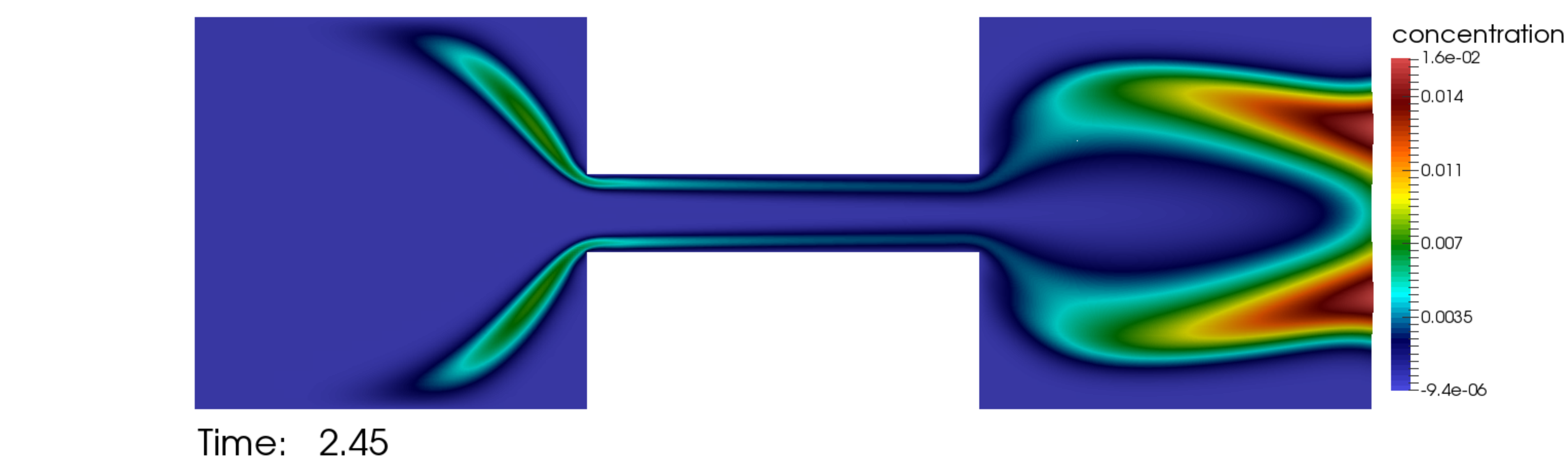}
~~
\includegraphics[width=.44\linewidth]{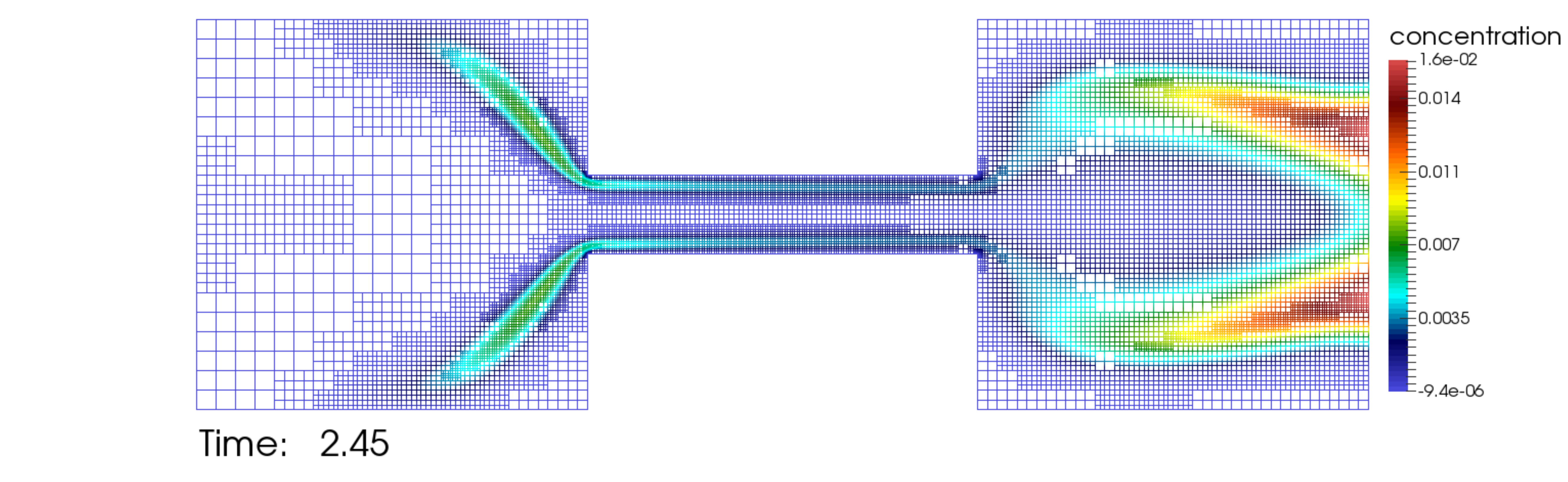}

\caption{Solution profiles and corresponding meshes of loop $\ell=8$
for Sec.~\ref{sec:6:3}.}
\label{fig:10:ex3:loop8}
\end{figure}
\begin{figure}
\centering

\includegraphics[width=.44\linewidth]{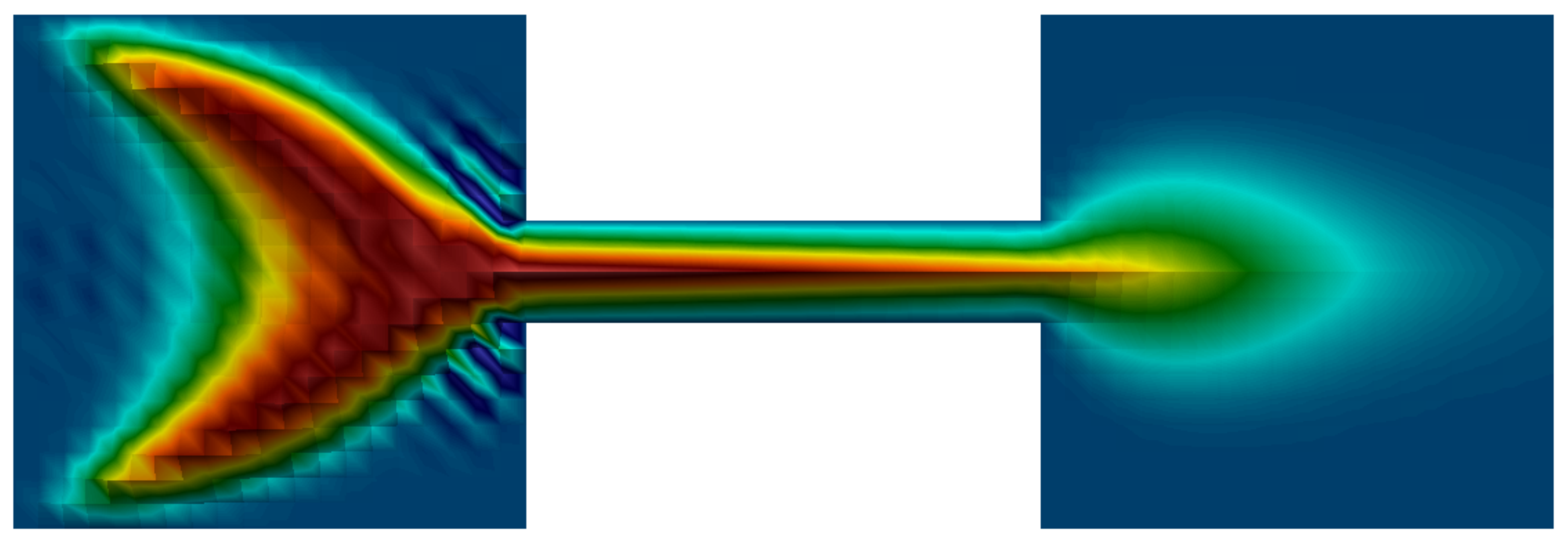}
~~
\includegraphics[width=.44\linewidth]{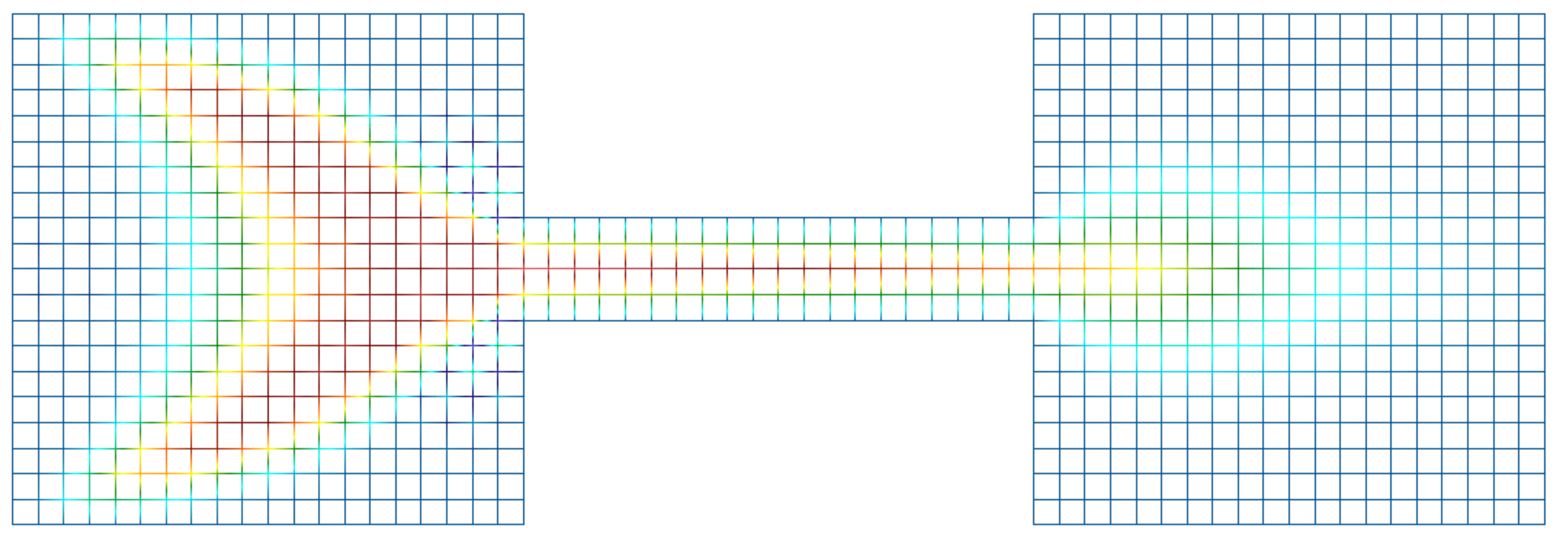}

\includegraphics[width=.44\linewidth]{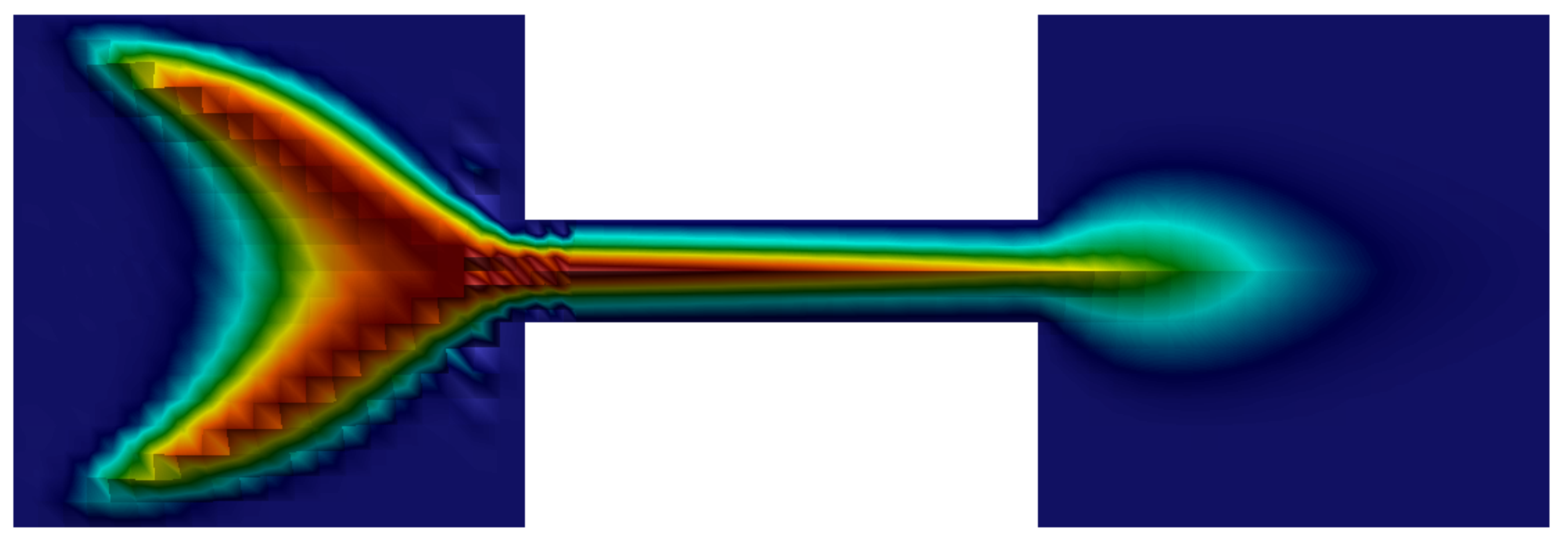}
~~
\includegraphics[width=.44\linewidth]{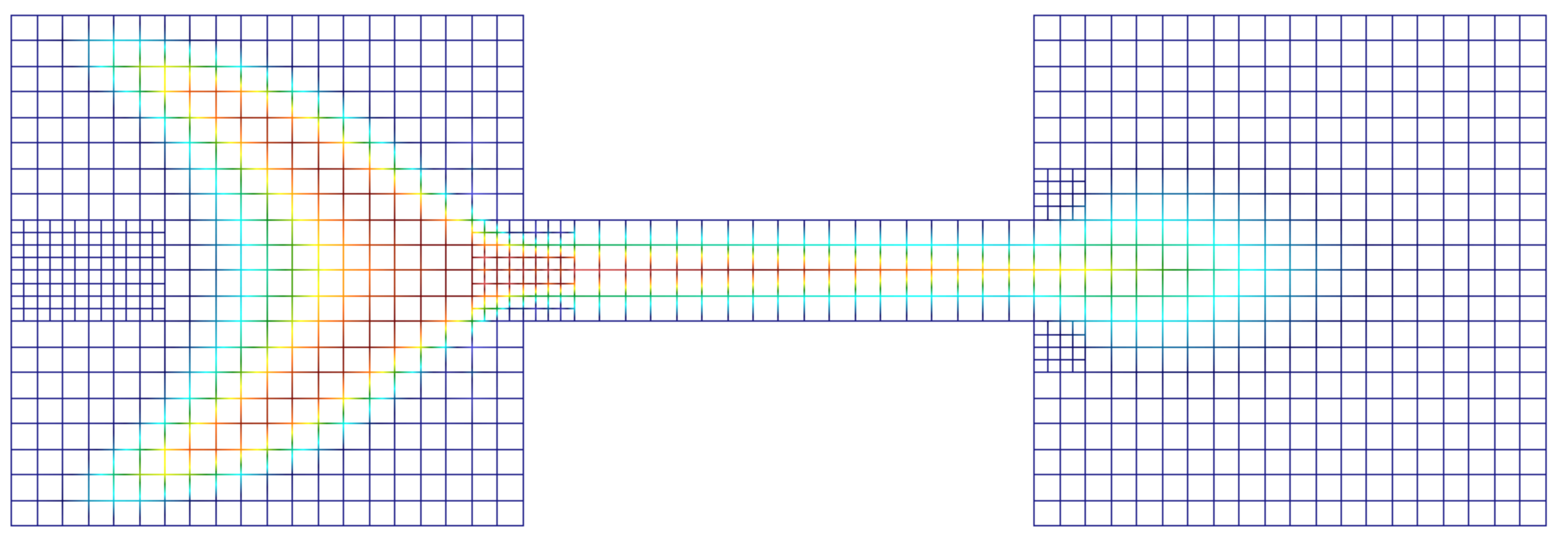}

\includegraphics[width=.44\linewidth]{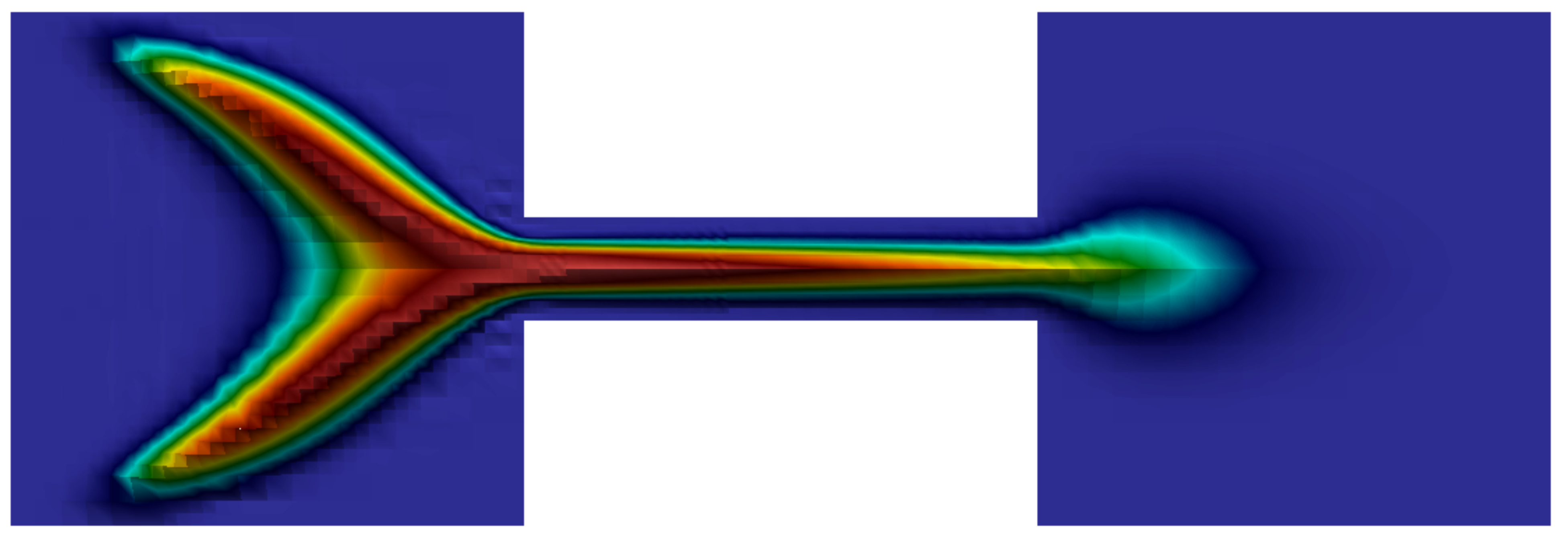}
~~
\includegraphics[width=.44\linewidth]{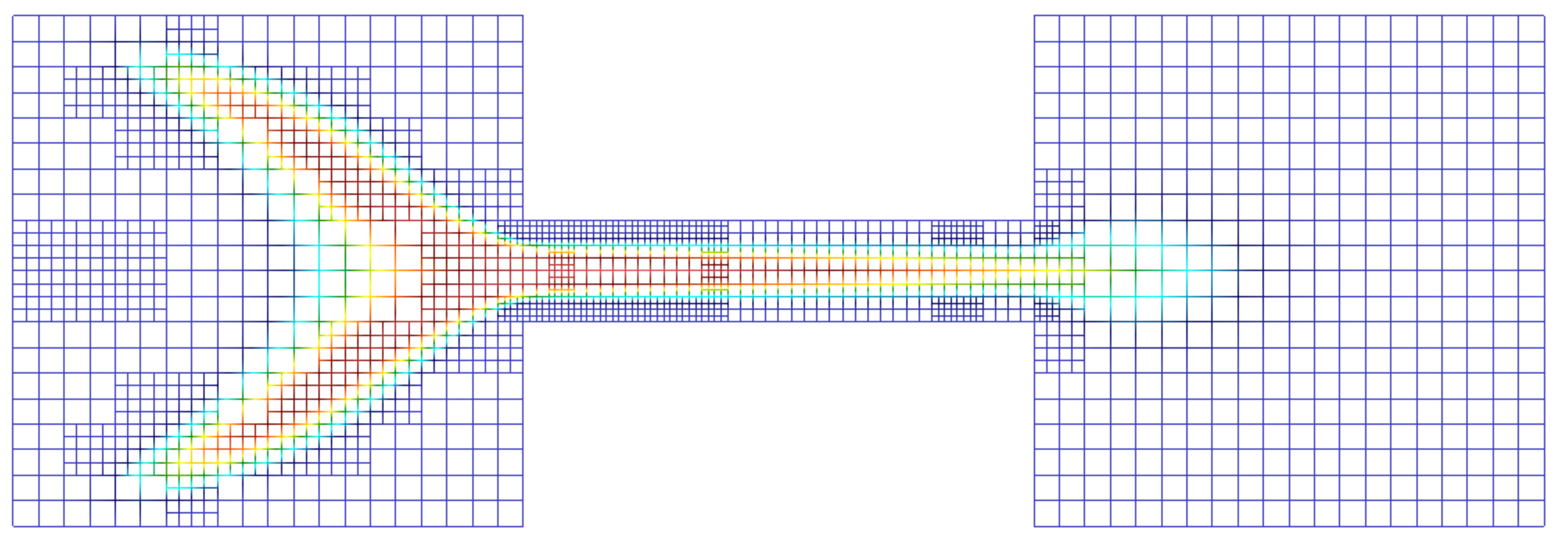}

\includegraphics[width=.44\linewidth]{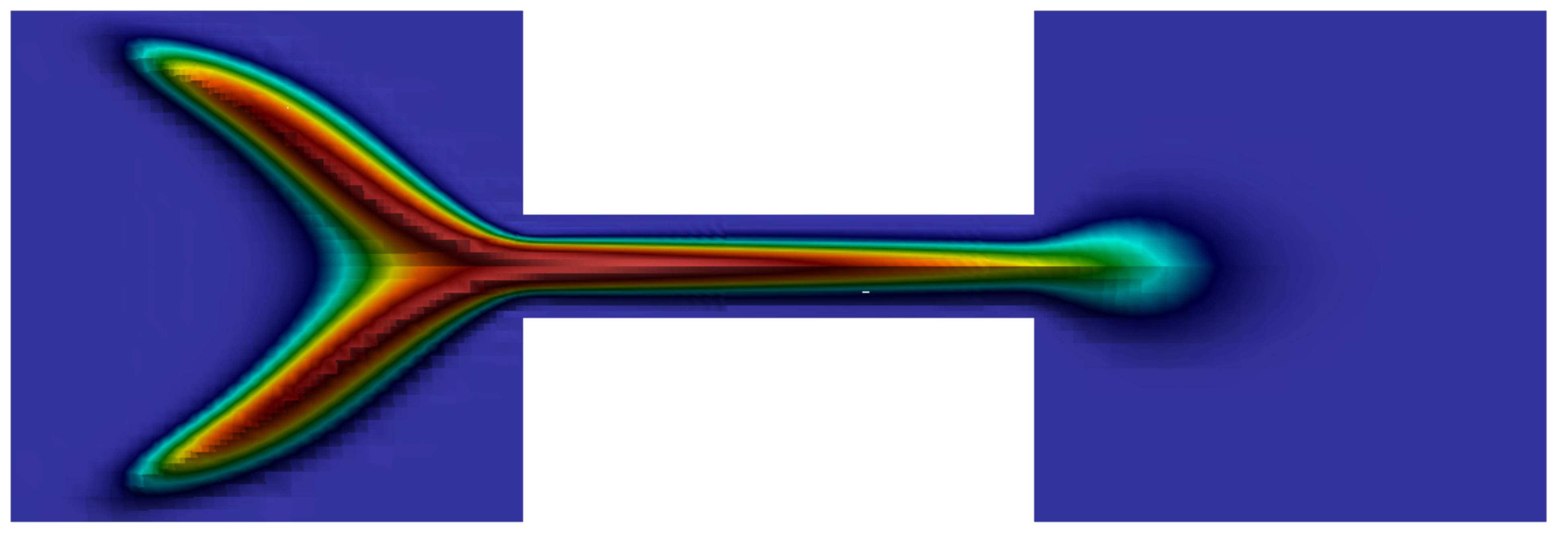}
~~
\includegraphics[width=.44\linewidth]{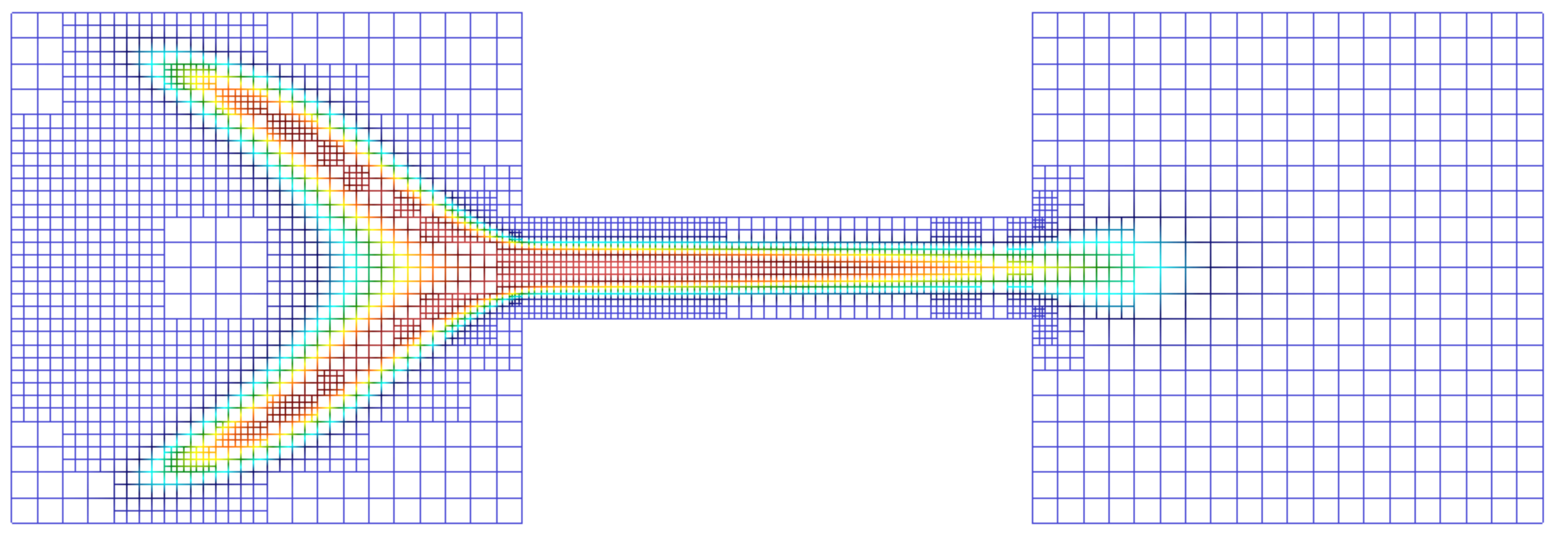}

\includegraphics[width=.44\linewidth]{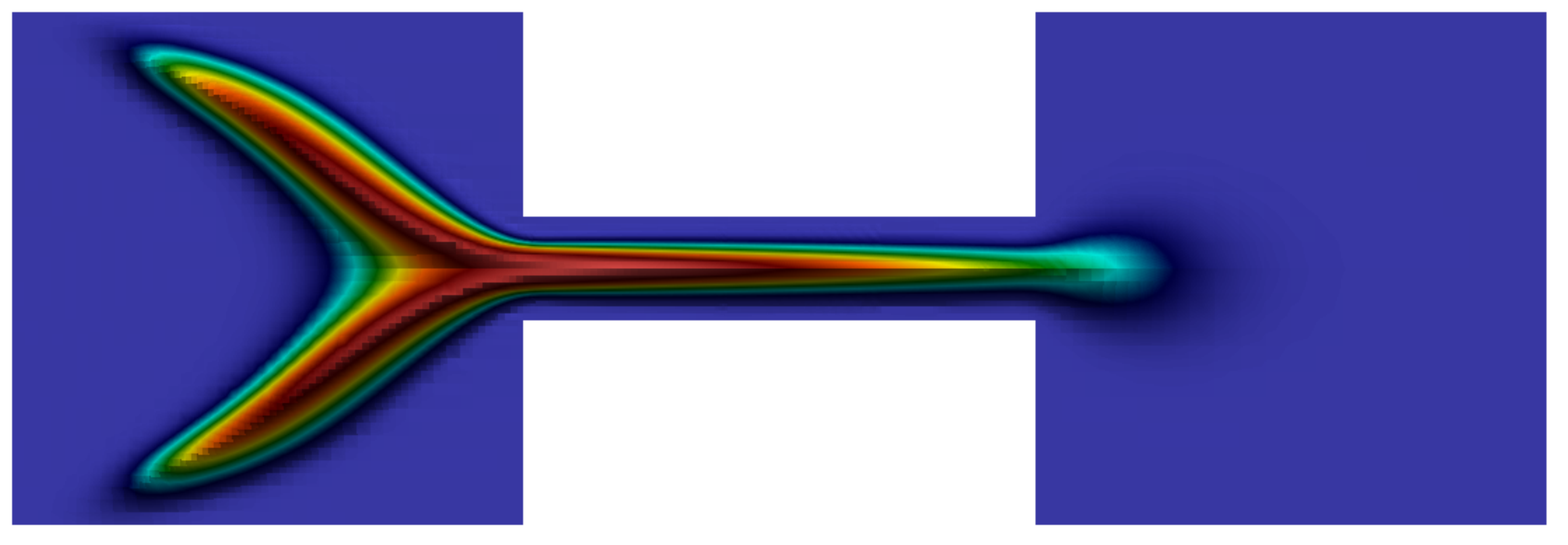}
~~
\includegraphics[width=.44\linewidth]{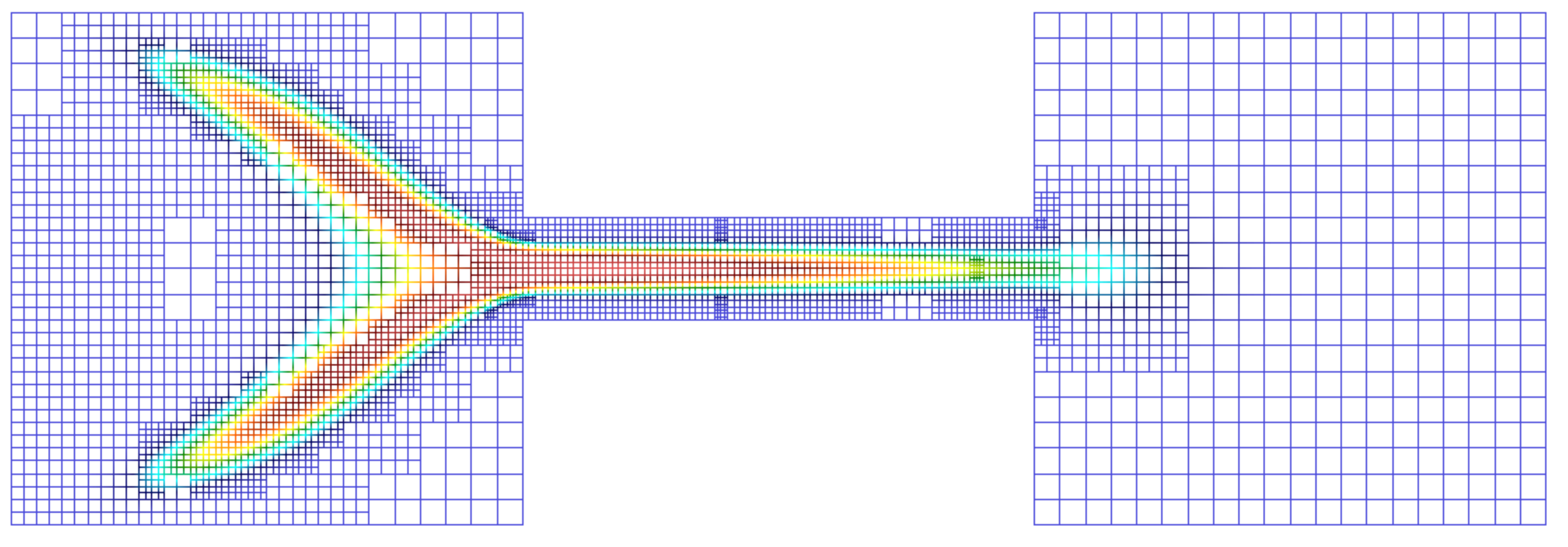}

\includegraphics[width=.44\linewidth]{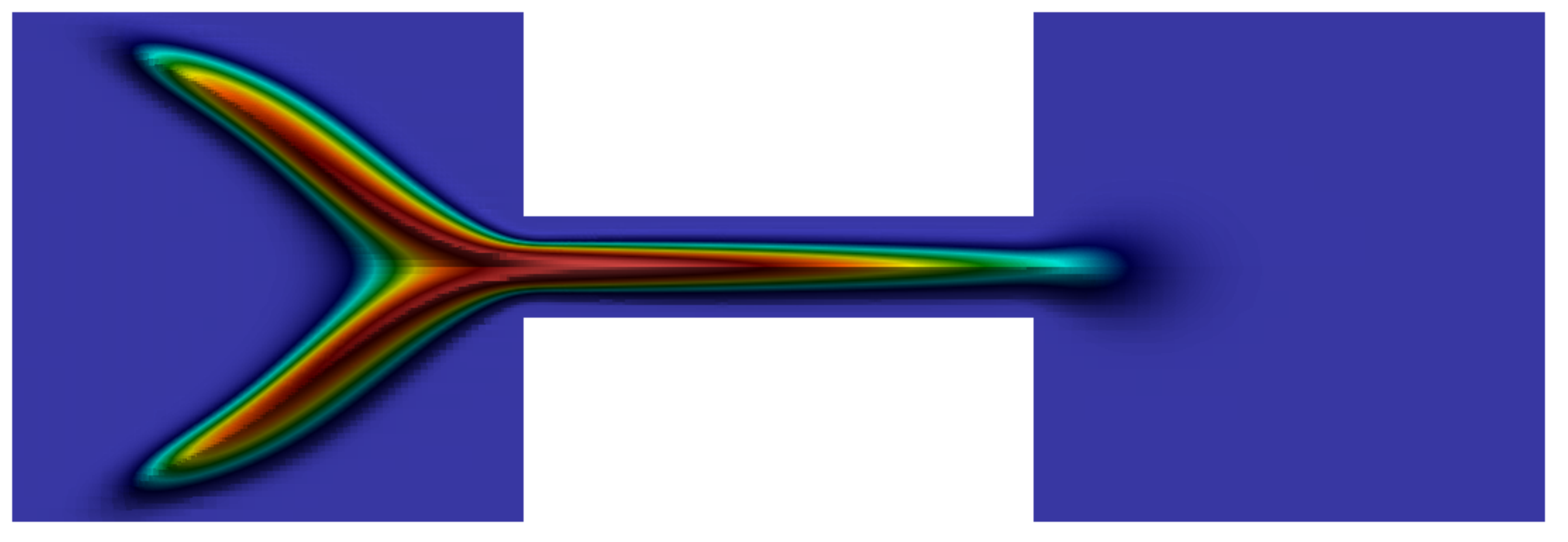}
~~
\includegraphics[width=.44\linewidth]{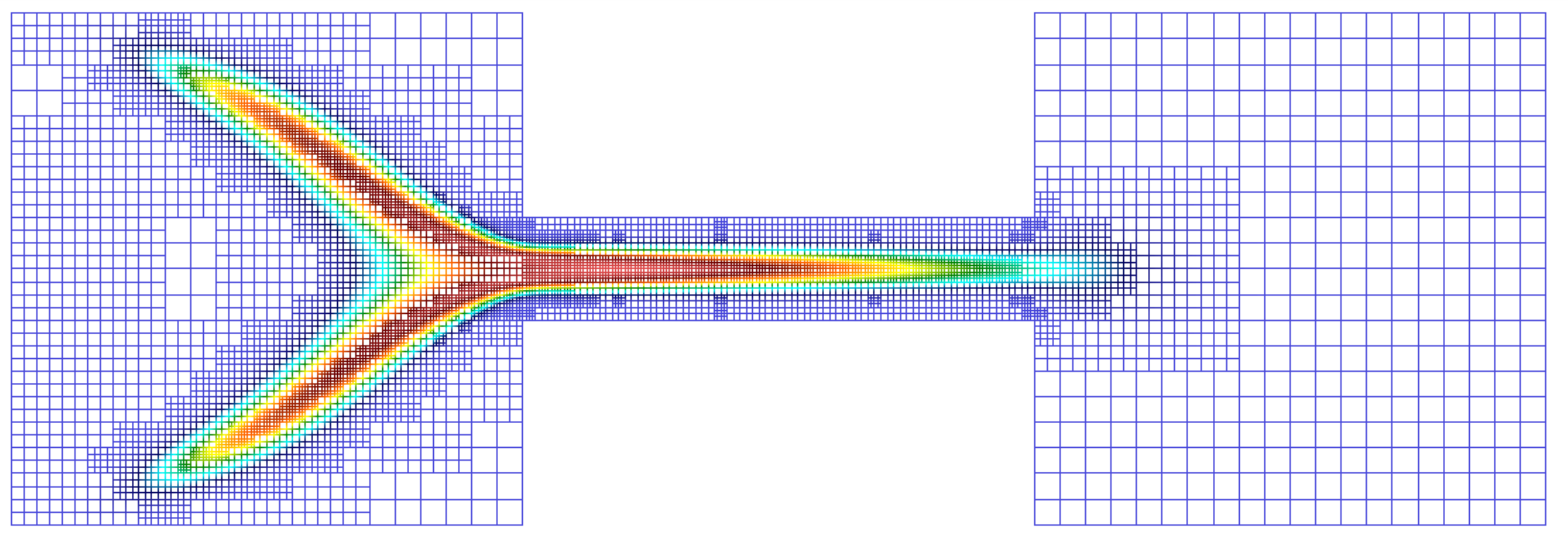}

\includegraphics[width=.44\linewidth]{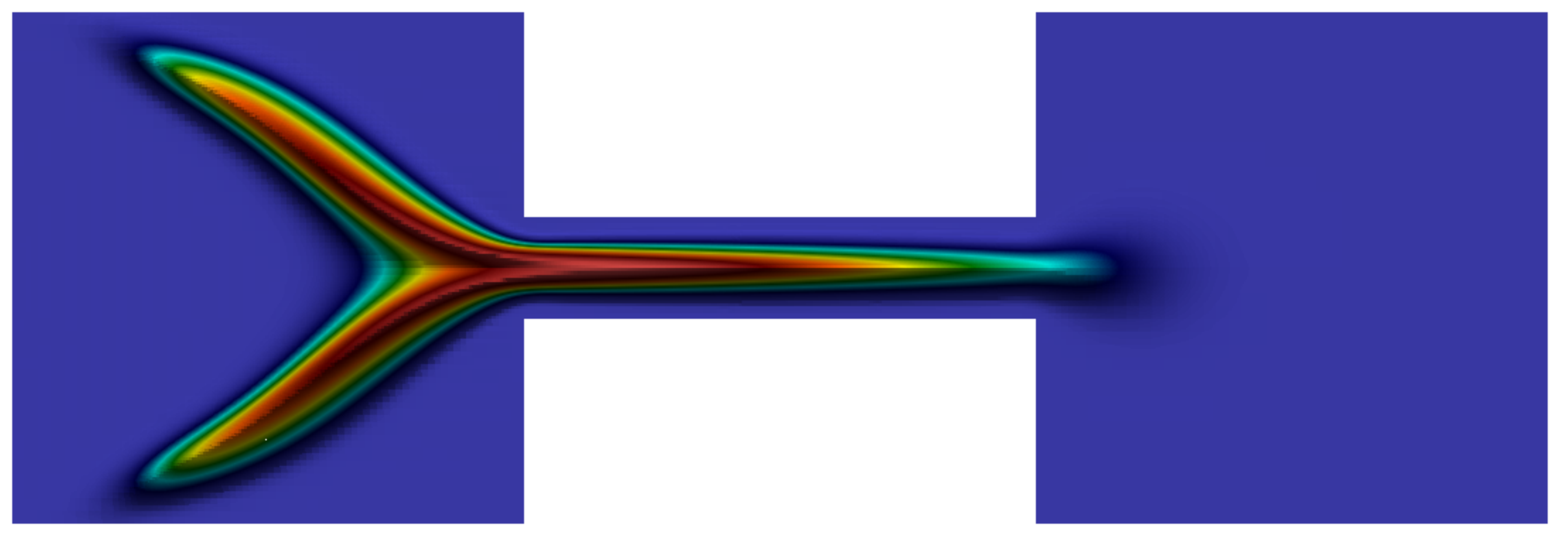}
~~
\includegraphics[width=.44\linewidth]{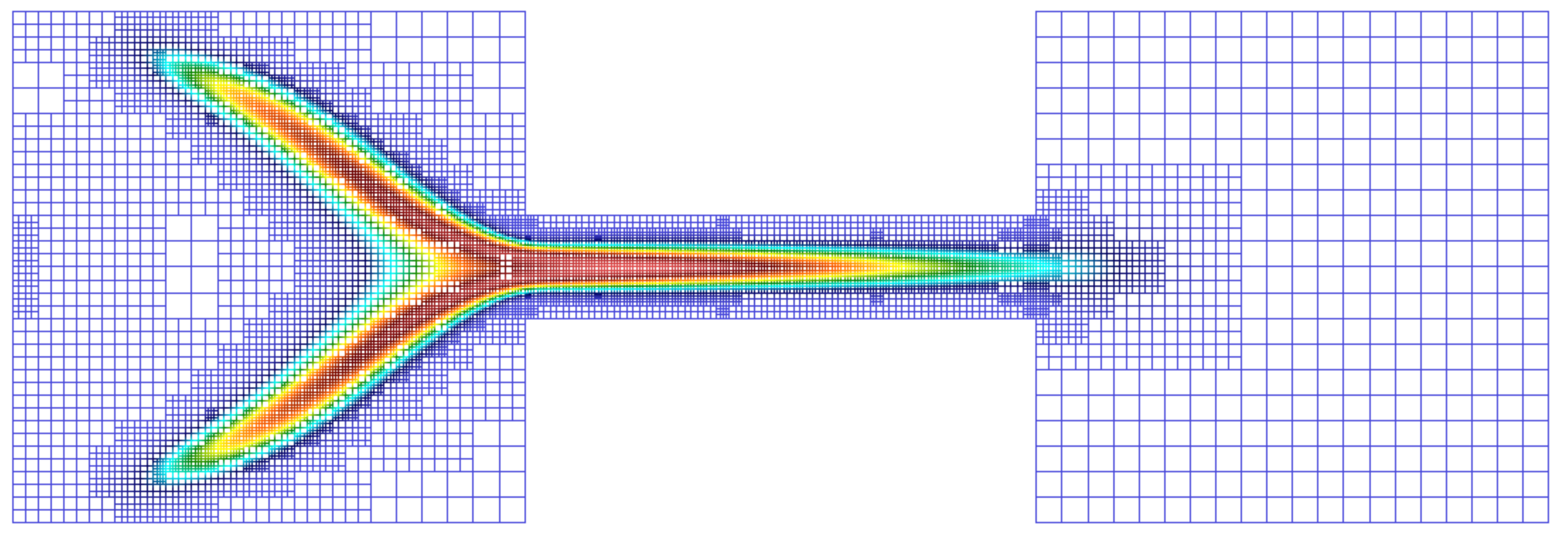}

\includegraphics[width=.44\linewidth]{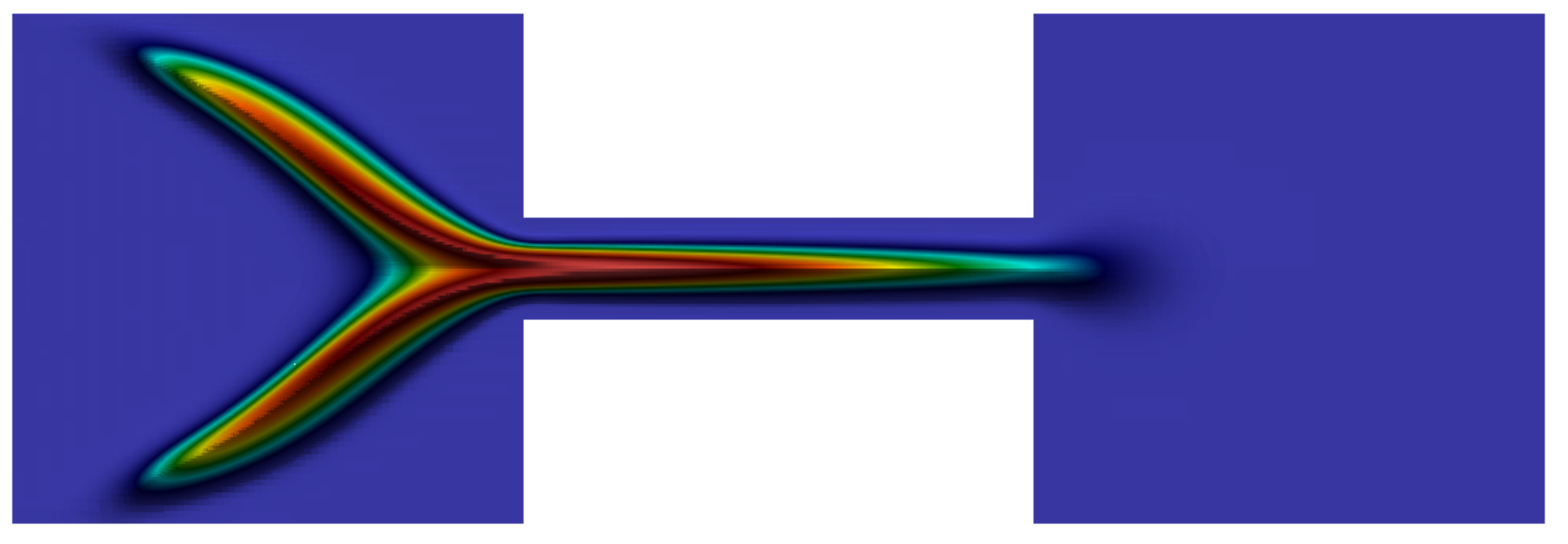}
~~
\includegraphics[width=.44\linewidth]{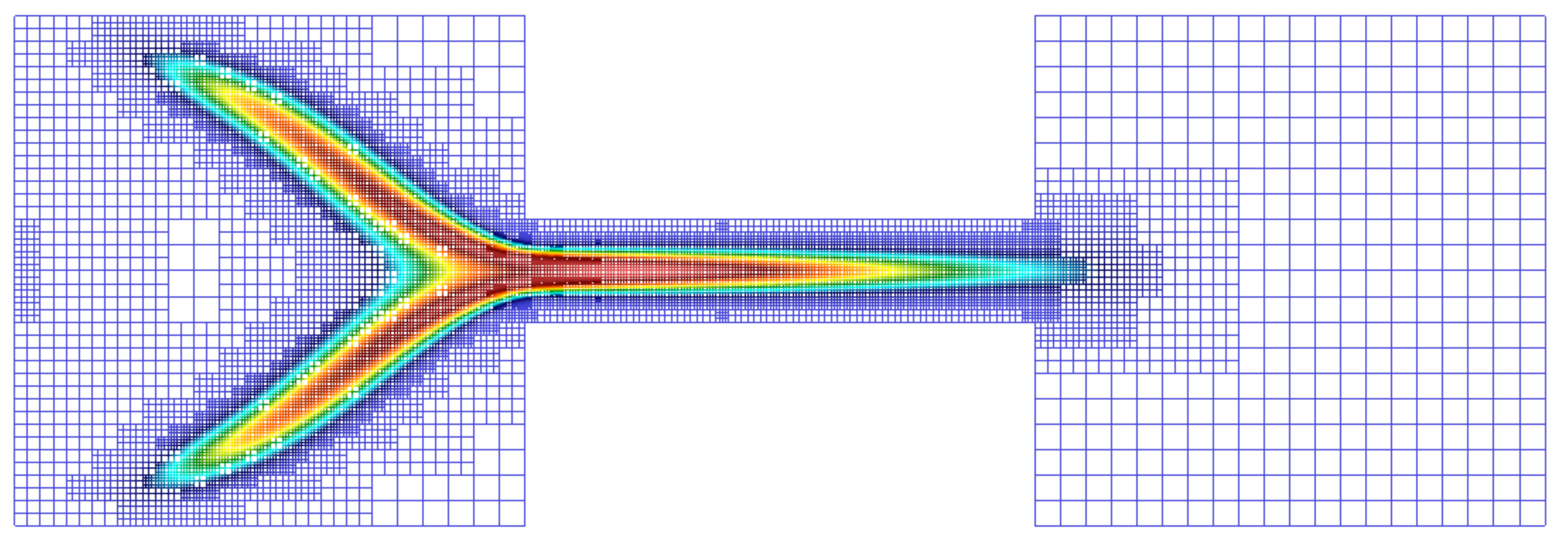}

\caption{Capturing of spurious oscillations with goal-oriented adaptivity
illustrated by comparative solution profiles and corresponding meshes of
the loops $\ell=1-8$ for Sec.~\ref{sec:6:3}.}
\label{fig:11:ex3:osc}
\end{figure}

\begin{figure}[h!]
\begin{minipage}{\linewidth}
\centering
\begin{tikzpicture}
\begin{axis}[%
width=3.9in,
height=1.0in,
scale only axis,
xlabel={t},
ylabel={\textcolor{navyblue}{$\tau_K(I_n^{1})$}
},
xmin=0.0,
xmax=2.5,
]
\addplot [
color=navyblue,
solid,
line width=1.5pt,
mark=*,
mark size = 1.,
only marks,
mark options={solid,navyblue}
]
table[row sep=crcr]{
0.1 0.1 \\
0.2 0.1 \\
0.3 0.1 \\
0.4 0.1 \\
0.5 0.1 \\
0.6 0.1 \\
0.7 0.1 \\
0.8 0.1 \\
0.9 0.1 \\
1 0.1 \\
1.1 0.1 \\
1.2 0.1 \\
1.3 0.1 \\
1.4 0.1 \\
1.5 0.1 \\
1.6 0.1 \\
1.7 0.1 \\
1.8 0.1 \\
1.9 0.1 \\
2 0.1 \\
2.1 0.1 \\
2.2 0.1 \\
2.3 0.1 \\
2.4 0.1 \\
2.5 0.1 \\
};
%
\end{axis}
\end{tikzpicture}
\end{minipage}

\begin{minipage}{\linewidth}
\centering
\begin{tikzpicture}
\begin{axis}[%
width=3.9in,
height=1.0in,
scale only axis,
/pgf/number format/.cd, 1000 sep={},
xlabel={t},
ylabel={\textcolor{navyblue}{$\tau_K(I_n^{4})$}
},
xmin=0.0,
xmax=2.5,
]

\addplot [
color=navyblue,
solid,
line width=1.0pt,
mark=*,
mark size = 1.5,
only marks,
mark options={solid,navyblue}
]
table[row sep=crcr]{
0.0125 0.0125 \\
0.025 0.0125 \\
0.05 0.025 \\
0.075 0.025 \\
0.1 0.025 \\
0.125 0.025 \\
0.15 0.025 \\
0.175 0.025 \\
0.2 0.025 \\
0.2125 0.0125 \\
0.225 0.0125 \\
0.25 0.025 \\
0.275 0.025 \\
0.3 0.025 \\
0.325 0.025 \\
0.35 0.025 \\
0.4 0.05 \\
0.45 0.05 \\
0.5 0.05 \\
0.55 0.05 \\
0.6 0.05 \\
0.65 0.05 \\
0.7 0.05 \\
0.75 0.05 \\
0.8 0.05 \\
0.9 0.1 \\
1 0.1 \\
1.1 0.1 \\
1.2 0.1 \\
1.3 0.1 \\
1.4 0.1 \\
1.5 0.1 \\
1.6 0.1 \\
1.7 0.1 \\
1.8 0.1 \\
1.9 0.1 \\
2 0.1 \\
2.1 0.1 \\
2.2 0.1 \\
2.3 0.1 \\
2.4 0.1 \\
2.5 0.1 \\
};
\end{axis}
\end{tikzpicture}
\end{minipage}
\begin{minipage}{\linewidth}
\centering
\begin{tikzpicture}
\begin{axis}[%
width=3.9in,
height=1.0in,
scale only axis,
/pgf/number format/.cd, 1000 sep={},
xlabel={t},
ylabel={\textcolor{navyblue}{$\tau_K(I_n^{8})$}
},
xmin=0.0,
xmax=2.5,
]

\addplot [
color=navyblue,
solid,
line width=1.0pt,
mark=*,
mark size = 1.5,
only marks,
mark options={solid,navyblue}
]
table[row sep=crcr]{
0.0015625 0.0015625 \\
0.003125 0.0015625 \\
0.00625 0.003125 \\
0.009375 0.003125 \\
0.0125 0.003125 \\
0.015625 0.003125 \\
0.01875 0.003125 \\
0.021875 0.003125 \\
0.025 0.003125 \\
0.03125 0.00625 \\
0.0375 0.00625 \\
0.04375 0.00625 \\
0.05 0.00625 \\
0.0625 0.0125 \\
0.075 0.0125 \\
0.0875 0.0125 \\
0.1 0.0125 \\
0.1125 0.0125 \\
0.125 0.0125 \\
0.1375 0.0125 \\
0.15 0.0125 \\
0.1625 0.0125 \\
0.175 0.0125 \\
0.1875 0.0125 \\
0.2 0.0125 \\
0.201563 0.0015625 \\
0.203125 0.0015625 \\
0.20625 0.003125 \\
0.209375 0.003125 \\
0.2125 0.003125 \\
0.21875 0.00625 \\
0.225 0.00625 \\
0.23125 0.00625 \\
0.2375 0.00625 \\
0.25 0.0125 \\
0.2625 0.0125 \\
0.275 0.0125 \\
0.2875 0.0125 \\
0.3 0.0125 \\
0.3125 0.0125 \\
0.325 0.0125 \\
0.3375 0.0125 \\
0.35 0.0125 \\
0.3625 0.0125 \\
0.375 0.0125 \\
0.3875 0.0125 \\
0.4 0.0125 \\
0.4125 0.0125 \\
0.425 0.0125 \\
0.4375 0.0125 \\
0.45 0.0125 \\
0.4625 0.0125 \\
0.475 0.0125 \\
0.4875 0.0125 \\
0.5 0.0125 \\
0.5125 0.0125 \\
0.525 0.0125 \\
0.55 0.025 \\
0.575 0.025 \\
0.6 0.025 \\
0.625 0.025 \\
0.65 0.025 \\
0.675 0.025 \\
0.7 0.025 \\
0.725 0.025 \\
0.75 0.025 \\
0.775 0.025 \\
0.8 0.025 \\
0.825 0.025 \\
0.85 0.025 \\
0.875 0.025 \\
0.9 0.025 \\
0.925 0.025 \\
0.95 0.025 \\
1 0.05 \\
1.05 0.05 \\
1.1 0.05 \\
1.15 0.05 \\
1.2 0.05 \\
1.25 0.05 \\
1.3 0.05 \\
1.35 0.05 \\
1.4 0.05 \\
1.5 0.1 \\
1.6 0.1 \\
1.7 0.1 \\
1.8 0.1 \\
1.9 0.1 \\
2 0.1 \\
2.1 0.1 \\
2.2 0.1 \\
2.3 0.1 \\
2.4 0.1 \\
2.5 0.1 \\
};
\end{axis}
\end{tikzpicture}
\end{minipage}

\caption{Distribution of the temporal step size $\tau_K$ of the transport problem
for a fixed $\sigma_K=2.5$ of the Stokes flow problem over the time 
interval $I=(0,T]$ for the initial (1) and after 4 and 8 DWR-loops.}
\label{fig:12:DistributionTauSigmaQS-Stokes}
\end{figure}

The solution profiles and corresponding adaptive meshes of the primal solution
$\concentration_{\tau h}^{1,0}$ of the loop $\ell=8$
for $t=0.15$, $t=0.70$, $t=0.92$, $t=1.33$ and $t=2.45$ are given by 
Fig.~\ref{fig:10:ex3:loop8}.
The refinement in space is adjusted to the position of the transported species
within the channel. It is located to the layers of the transported species, 
whereas the mesh stays coarse in the non-affected area.  
In Fig.~\ref{fig:11:ex3:osc} we present a comparative study of the solution
profile and corresponding meshes for $t=0.95$ over the adaptivity loops.
For $\ell=1,2,3$ obvious spurious oscillations in the left square and at the 
beginning of the constriction are existing, which are captured and resolved by 
the goal-oriented adaptivity by taking spatial mesh refinements along the 
layers of the transported species within the left square and within the 
constriction of the channel.
For $\ell>3$ the spatial refinements capture especially the solution profile
fronts with strong gradients with a focus on the high-convective middle of
the spatial domain.
In Fig.~\ref{fig:12:DistributionTauSigmaQS-Stokes} we visualize the temporal
distribution of the transport problem for several DWR-loops. The time cell 
lengths of the Stokes flow problem is kept fixed with value $\sigma_K=2.5$ for
all DWR-loops here and thus explicitly not displayed. We observe an adaptive
refinement in time at the beginning, consistent with the restriction in time of
the inflow boundary condition. The closer we get to the final time point $T$ the 
coarser the temporal mesh is chosen.

\begin{table}[ht]
\centering

\begin{tabular}{c | rrr | cc}
\hline
\hline
$\ell$ & $N$ & $N_K^{\text{max}}$ & $N_{\text{DoF}}^{\text{tot}}$ &
$\eta_h$ & $\eta_\tau$\\
\hline
1 & 25  &  880 &  24425 & 3.5795e-03 & 1.1452e-02 \\
2 & 29 &   880 &  28333 & 3.8619e-03 & 3.2318e-03 \\ 
3 & 32 &  1456 &  39456 & 2.9354e-03 & 5.4042e-03 \\
4 & 42 &  2116 &  62528 & 2.5532e-03 & 5.2001e-03 \\
5 & 51 &  4492 & 132483 & 2.3178e-03 & 5.2170e-03 \\
6 & 70 &  7072 & 239266 & 1.8934e-03 & 3.8571e-03 \\
7 & 79 & 10744 & 371015 & 1.7406e-03 & 2.3554e-03 \\
8 & 89 & 15376 & 619071 & 1.6069e-03 & 1.2974e-03 \\
\hline
\end{tabular}

\caption{Goal-oriented temporal and spatial refinements for the transport problem
in Sec.~\ref{sec:6:3}.
$\ell$ denotes the refinement level loop,
$N$ the accumulated total cells in time,
$N_K^{\text{max}}$ the maximal number of cells on a slab,
$N_{\text{DoF}}^{\text{tot}}$ the total space-time degrees of freedom 
and
$\eta_h$ and $\eta_\tau$ the computed error indicators in space and time, 
respectively.
}
\label{tab:5:ex3:eta}
\end{table}
The refinement in space and time is automatically balanced due to the dynamic
choice of $\theta_h^\textnormal{top}$ and $\theta_\tau^\textnormal{top}$
given by \eqref{eq:32:balancing} and is illustrated by Tab.~\ref{tab:5:ex3:eta}.
Regarding the spatial and temporal error indicators (cf. columns five and six of
Tab.~\ref{tab:5:ex3:eta}) a good equilibration can be observed within the final loop,
whereas in the first step a mismatch occurs resulting in a solely temporal 
refinement between $\ell=1$ and $\ell=2$.

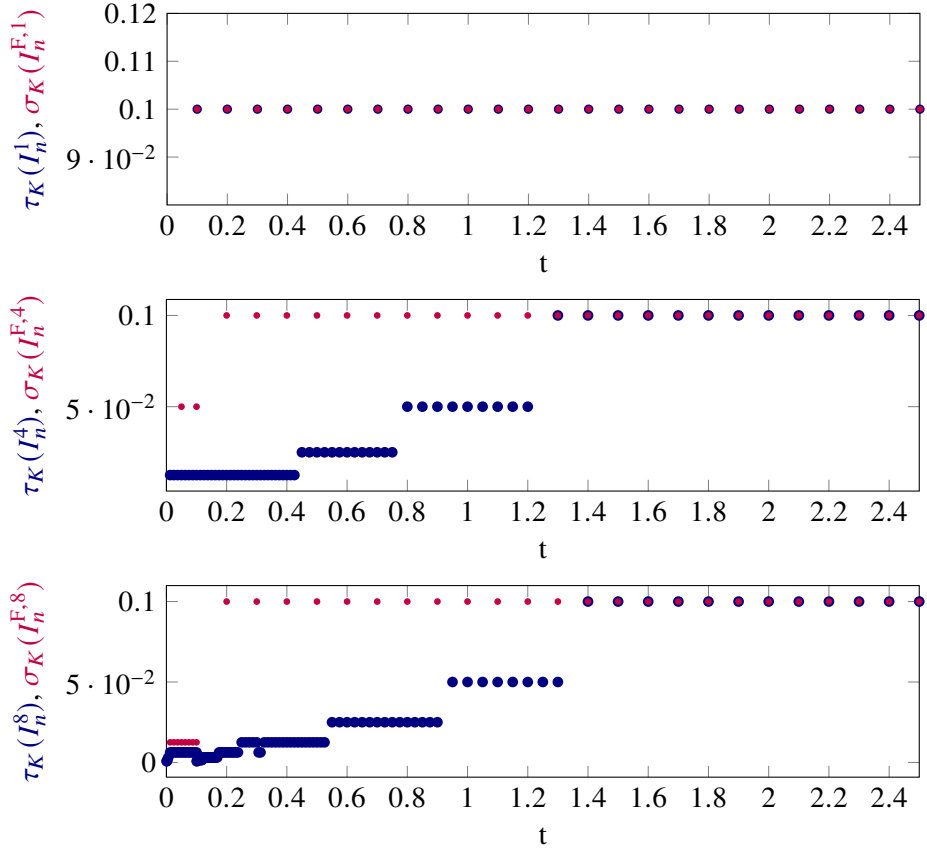
\begin{figure}[h!]
\begin{minipage}{\linewidth}
\centering
\begin{tikzpicture}
\begin{axis}[%
width=3.9in,
height=1.0in,
scale only axis,
xlabel={t},
ylabel={\textcolor{navyblue}{$\tau_K(I_n^{1})$}, 
\textcolor{HSUred}{$\sigma_K(I_n^{\textnormal{F},1})$}},
xmin=0.0,
xmax=2.5,
]
\addplot [
color=navyblue,
solid,
line width=1.5pt,
mark=*,
mark size = 1.,
only marks,
mark options={solid,navyblue}
]
table[row sep=crcr]{
0.1 0.1 \\
0.2 0.1 \\
0.3 0.1 \\
0.4 0.1 \\
0.5 0.1 \\
0.6 0.1 \\
0.7 0.1 \\
0.8 0.1 \\
0.9 0.1 \\
1 0.1 \\
1.1 0.1 \\
1.2 0.1 \\
1.3 0.1 \\
1.4 0.1 \\
1.5 0.1 \\
1.6 0.1 \\
1.7 0.1 \\
1.8 0.1 \\
1.9 0.1 \\
2 0.1 \\
2.1 0.1 \\
2.2 0.1 \\
2.3 0.1 \\
2.4 0.1 \\
2.5 0.1 \\
};
\addplot [
color=HSUred,
solid,
line width=0.5pt,
mark=*,
mark size = 1.,
only marks,
mark options={fill=HSUred}
]
table[row sep=crcr]{
0.1 0.1 \\
0.2 0.1 \\
0.3 0.1 \\
0.4 0.1 \\
0.5 0.1 \\
0.6 0.1 \\
0.7 0.1 \\
0.8 0.1 \\
0.9 0.1 \\
1 0.1 \\
1.1 0.1 \\
1.2 0.1 \\
1.3 0.1 \\
1.4 0.1 \\
1.5 0.1 \\
1.6 0.1 \\
1.7 0.1 \\
1.8 0.1 \\
1.9 0.1 \\
2 0.1 \\
2.1 0.1 \\
2.2 0.1 \\
2.3 0.1 \\
2.4 0.1 \\
2.5 0.1 \\
};
\end{axis}
\end{tikzpicture}
\end{minipage}

\begin{minipage}{\linewidth}
\centering
\begin{tikzpicture}
\begin{axis}[%
width=3.9in,
height=1.0in,
scale only axis,
/pgf/number format/.cd, 1000 sep={},
xlabel={t},
ylabel={\textcolor{navyblue}{$\tau_K(I_n^{4})$}, 
\textcolor{HSUred}{$\sigma_K(I_n^{\textnormal{F},4})$}},
xmin=0.0,
xmax=2.5,
]

\addplot [
color=navyblue,
solid,
line width=1.0pt,
mark=*,
mark size = 1.5,
only marks,
mark options={solid,navyblue}
]
table[row sep=crcr]{
0.0125 0.0125 \\
0.025 0.0125 \\
0.0375 0.0125 \\
0.05 0.0125 \\
0.0625 0.0125 \\
0.075 0.0125 \\
0.0875 0.0125 \\
0.1 0.0125 \\
0.1125 0.0125 \\
0.125 0.0125 \\
0.1375 0.0125 \\
0.15 0.0125 \\
0.1625 0.0125 \\
0.175 0.0125 \\
0.1875 0.0125 \\
0.2 0.0125 \\
0.2125 0.0125 \\
0.225 0.0125 \\
0.2375 0.0125 \\
0.25 0.0125 \\
0.2625 0.0125 \\
0.275 0.0125 \\
0.2875 0.0125 \\
0.3 0.0125 \\
0.3125 0.0125 \\
0.325 0.0125 \\
0.3375 0.0125 \\
0.35 0.0125 \\
0.3625 0.0125 \\
0.375 0.0125 \\
0.3875 0.0125 \\
0.4 0.0125 \\
0.4125 0.0125 \\
0.425 0.0125 \\
0.45 0.025 \\
0.475 0.025 \\
0.5 0.025 \\
0.525 0.025 \\
0.55 0.025 \\
0.575 0.025 \\
0.6 0.025 \\
0.625 0.025 \\
0.65 0.025 \\
0.675 0.025 \\
0.7 0.025 \\
0.725 0.025 \\
0.75 0.025 \\
0.8 0.05 \\
0.85 0.05 \\
0.9 0.05 \\
0.95 0.05 \\
1 0.05 \\
1.05 0.05 \\
1.1 0.05 \\
1.15 0.05 \\
1.2 0.05 \\
1.3 0.1 \\
1.4 0.1 \\
1.5 0.1 \\
1.6 0.1 \\
1.7 0.1 \\
1.8 0.1 \\
1.9 0.1 \\
2 0.1 \\
2.1 0.1 \\
2.2 0.1 \\
2.3 0.1 \\
2.4 0.1 \\
2.5 0.1 \\
};

\addplot [
color=HSUred,
solid,
line width=0.5pt,
mark=*,
mark size = 1.,
only marks,
mark options={fill=HSUred}
]
table[row sep=crcr]{
0.05 0.05 \\
0.1 0.05 \\
0.2 0.1 \\
0.3 0.1 \\
0.4 0.1 \\
0.5 0.1 \\
0.6 0.1 \\
0.7 0.1 \\
0.8 0.1 \\
0.9 0.1 \\
1 0.1 \\
1.1 0.1 \\
1.2 0.1 \\
1.3 0.1 \\
1.4 0.1 \\
1.5 0.1 \\
1.6 0.1 \\
1.7 0.1 \\
1.8 0.1 \\
1.9 0.1 \\
2 0.1 \\
2.1 0.1 \\
2.2 0.1 \\
2.3 0.1 \\
2.4 0.1 \\
2.5 0.1 \\
};
\end{axis}
\end{tikzpicture}
\end{minipage}
\begin{minipage}{\linewidth}
\centering
\begin{tikzpicture}
\begin{axis}[%
width=3.9in,
height=1.0in,
scale only axis,
/pgf/number format/.cd, 1000 sep={},
xlabel={t},
ylabel={\textcolor{navyblue}{$\tau_K(I_n^{8})$}, 
\textcolor{HSUred}{$\sigma_K(I_n^{\textnormal{F},8})$}},
xmin=0.0,
xmax=2.5,
]

\addplot [
color=navyblue,
solid,
line width=1.0pt,
mark=*,
mark size = 1.5,
only marks,
mark options={solid,navyblue}
]
table[row sep=crcr]{
0.00078125 0.00078125 \\
0.0015625 0.00078125 \\
0.003125 0.0015625 \\
0.00625 0.003125 \\
0.0125 0.00625 \\
0.01875 0.00625 \\
0.025 0.00625 \\
0.03125 0.00625 \\
0.0375 0.00625 \\
0.04375 0.00625 \\
0.05 0.00625 \\
0.05625 0.00625 \\
0.0625 0.00625 \\
0.06875 0.00625 \\
0.075 0.00625 \\
0.08125 0.00625 \\
0.0875 0.00625 \\
0.09375 0.00625 \\
0.1 0.00625 \\
0.100781 0.00078125 \\
0.101562 0.00078125 \\
0.103125 0.0015625 \\
0.104688 0.0015625 \\
0.10625 0.0015625 \\
0.107813 0.0015625 \\
0.109375 0.0015625 \\
0.110937 0.0015625 \\
0.1125 0.0015625 \\
0.114063 0.0015625 \\
0.115625 0.0015625 \\
0.117188 0.0015625 \\
0.11875 0.0015625 \\
0.121875 0.003125 \\
0.125 0.003125 \\
0.128125 0.003125 \\
0.13125 0.003125 \\
0.134375 0.003125 \\
0.1375 0.003125 \\
0.140625 0.003125 \\
0.14375 0.003125 \\
0.146875 0.003125 \\
0.15 0.003125 \\
0.153125 0.003125 \\
0.15625 0.003125 \\
0.159375 0.003125 \\
0.1625 0.003125 \\
0.165625 0.003125 \\
0.16875 0.003125 \\
0.175 0.00625 \\
0.18125 0.00625 \\
0.1875 0.00625 \\
0.19375 0.00625 \\
0.2 0.00625 \\
0.20625 0.00625 \\
0.2125 0.00625 \\
0.21875 0.00625 \\
0.225 0.00625 \\
0.23125 0.00625 \\
0.2375 0.00625 \\
0.25 0.0125 \\
0.2625 0.0125 \\
0.275 0.0125 \\
0.2875 0.0125 \\
0.3 0.0125 \\
0.30625 0.00625 \\
0.3125 0.00625 \\
0.325 0.0125 \\
0.3375 0.0125 \\
0.35 0.0125 \\
0.3625 0.0125 \\
0.375 0.0125 \\
0.3875 0.0125 \\
0.4 0.0125 \\
0.4125 0.0125 \\
0.425 0.0125 \\
0.4375 0.0125 \\
0.45 0.0125 \\
0.4625 0.0125 \\
0.475 0.0125 \\
0.4875 0.0125 \\
0.5 0.0125 \\
0.5125 0.0125 \\
0.525 0.0125 \\
0.55 0.025 \\
0.575 0.025 \\
0.6 0.025 \\
0.625 0.025 \\
0.65 0.025 \\
0.675 0.025 \\
0.7 0.025 \\
0.725 0.025 \\
0.75 0.025 \\
0.775 0.025 \\
0.8 0.025 \\
0.825 0.025 \\
0.85 0.025 \\
0.875 0.025 \\
0.9 0.025 \\
0.95 0.05 \\
1 0.05 \\
1.05 0.05 \\
1.1 0.05 \\
1.15 0.05 \\
1.2 0.05 \\
1.25 0.05 \\
1.3 0.05 \\
1.4 0.1 \\
1.5 0.1 \\
1.6 0.1 \\
1.7 0.1 \\
1.8 0.1 \\
1.9 0.1 \\
2 0.1 \\
2.1 0.1 \\
2.2 0.1 \\
2.3 0.1 \\
2.4 0.1 \\
2.5 0.1 \\
};

\addplot [
color=HSUred,
solid,
line width=0.5pt,
mark=*,
mark size = 1.,
only marks,
mark options={fill=HSUred}
]
table[row sep=crcr]{
0.0125 0.0125 \\
0.025 0.0125 \\
0.0375 0.0125 \\
0.05 0.0125 \\
0.0625 0.0125 \\
0.075 0.0125 \\
0.0875 0.0125 \\
0.1 0.0125 \\
0.2 0.1 \\
0.3 0.1 \\
0.4 0.1 \\
0.5 0.1 \\
0.6 0.1 \\
0.7 0.1 \\
0.8 0.1 \\
0.9 0.1 \\
1 0.1 \\
1.1 0.1 \\
1.2 0.1 \\
1.3 0.1 \\
1.4 0.1 \\
1.5 0.1 \\
1.6 0.1 \\
1.7 0.1 \\
1.8 0.1 \\
1.9 0.1 \\
2 0.1 \\
2.1 0.1 \\
2.2 0.1 \\
2.3 0.1 \\
2.4 0.1 \\
2.5 0.1 \\
};
\end{axis}
\end{tikzpicture}
\end{minipage}

\caption{Distribution of the temporal step size $\tau_K$ of the transport problem
and $\sigma_K$ of the Stokes flow problem over the time 
interval $I=(0,T]$ for the initial (1) and after 4 and 8 DWR-loops.}
\label{fig:13:DistributionTauSigmaEx2-instaStokes}
\end{figure}

Finally, we modify the parabolic inflow condition for the Stokes flow problem in 
order to investigate our multirate-in-time approach for the present example.
More precisely, on the left boundary $\Gamma_{\textnormal{inflow}}$ the inflow
condition $\convection_D$ is now given by
\begin{equation}
\label{eq:33:insta-inflow-condition}
\convection_D =
\begin{cases}
 \frac{\arctan(t)}{\pi/2}\cdot(1-4x_2^2,0)^\top & \textnormal{ for } 0 \leq t \leq 0.1\,,\\
 (1,0)^\top & \textnormal{ for } 0.1 < t \leq T\,.
\end{cases}
\end{equation}
Moreover, for the transport problem, the Dirichlet boundary function value is 
homogeneous on $\Gamma_D$ except for the line $(-1,-1) \times (-0.4,0.4)$ and 
time $0 \leq t \leq 0.1$ where the constant value
\begin{displaymath}
\concentration(\boldsymbol{x},t)=1
\end{displaymath}
is prescribed on the solution. 
Therefore, the time domain $I=(0,2.5)$ is now discretized with the same initial 
$\tau=\sigma=0.1$ for the transport and the Stokes flow problem for the first 
loop $\ell=1$. In Fig.~\ref{fig:13:DistributionTauSigmaEx2-instaStokes} we visualize 
the distribution of the  adaptively determined time cell lengths $\tau_K$ and
$\sigma_K$ used for the transport  and Stokes flow problem, respectively,
over the whole time interval $I$ for different DWR refinement loops.
We observe a similar behavior as displayed in Fig.~\ref{fig:12:DistributionTauSigmaQS-Stokes}.
The temporal mesh is refined close to the time conditions of the respective 
inflow boundaries for both problems, where the refinement in time for 
the Stokes flow problem is chosen to refine those slabs related to the inflow 
condition \eqref{eq:33:insta-inflow-condition} for each second DWR-loop. Away from
the temporal inflow condition both temporal meshes stay coarse.

\section{Conclusion}
\label{sec:7:conclusion}

In this work we presented a multirate-in-time approach regarding different time
scales for a rapidly changing transport coupled with a slowly creeping Stokes 
flow. The transport problem is represented by a convection-dominated 
convection-diffusion-reaction equation which is for this reason stabilized using the 
residual based SUPG method.
Both subproblems are discretized using a discontinuous Galerkin method dG($r$) 
with an arbitrary polynomial degree $r \geq 0$ in time and a continuous Galerkin
method cG($p$) with an arbitrary polynomial degree $p \geq 1$ in space.
A goal-oriented a posteriori error representation based on the Dual
Weighted Residual method was derived for the transport problem. 
This error representation is splitted into an amount in space and time whose 
localized forms serve as error indicators for the adaptive mesh refinement 
process in space and time.
The temporal weights of the DWR adaptivity process are approximated by a 
higher-order extrapolation approach whereas the spatial weights are approximated 
by higher-order finite elements.
The practical realization of the space-time slabs is based on tensor-product
spaces which enables for an efficient and flexible software implementation of
the underlying approach. 
In numerical experiments we verified expected experimental orders of convergence
of the underlying subproblems as well as the coupled problem. 
Furthermore, space-time adaptivity studies for the coupled problem were 
investigated for an academic test problem as well as a problem of practical 
interest, leading to high-efficient adaptively refined meshes in space and time. 
Effectivity indices close to one and well-balanced error indicators in space 
and time were obtained.
Spurious oscillations that typically arise in numerical approximations of 
convection-dominated problems could be reduced significantly. 
Finally, the here presented approach for coupled free flow and species transport is fairly 
general and can be easily adopted to other multi-physics systems coupling phenomena 
that are characterized by strongly differing time scales.

\end{document}